\documentclass[11pt]{amsart}
\usepackage{amsmath}
\usepackage{amsxtra}
\usepackage{amscd}
\usepackage{amsthm}
\usepackage{amsfonts}
\usepackage{amssymb}
\usepackage{eucal}

\newcommand{\nc}{\newcommand}

\numberwithin{equation}{section}

\swapnumbers

\newtheorem{thm}[subsection]{Theorem}
\newtheorem{prop}[subsection]{Proposition}
\newtheorem{lem}[subsection]{Lemma}
\newtheorem{cor}[subsection]{Corollary}

\newtheorem{conj}[subsection]{Conjecture}
\nc{\ssec}{\subsection}
\nc{\sssec}{\subsubsection}

\newcommand{\thmref}[1]{Theorem~\ref{#1}}
\newcommand{\secref}[1]{Sect.~\ref{#1}}
\newcommand{\lemref}[1]{Lemma~\ref{#1}}

\newcommand{\propref}[1]{Proposition~\ref{#1}}

\newcommand{\corref}[1]{Corollary~\ref{#1}}
\newcommand{\conjref}[1]{Conjecture~\ref{#1}}
\newcommand{\remref}[1]{Remark~\ref{#1}}

\nc{\on}{\operatorname}
\nc{\dmo}{\DeclareMathOperator}
\nc{\ZZ}{{\mathbb Z}}
\nc{\NN}{{\mathbb N}}
\nc{\FF}{{\mathbb F}}
\nc{\CC}{{\mathbb C}}

\nc{\GG}{{\mathbb G}}
\nc\one{{\mathbf 1}}
\nc{\OO}{{\mathcal O}}
\nc{\DD}{{\mathbb D}}
\renewcommand{\AA}{{\mathbb A}}
\nc{\Fq}{{\mathbb F}_q}
\nc{\Fqb}{\overline{\mathbb F}_q}
\nc{\Ql}{\ol{\mathbb Q}_\ell}

\nc{\Eis}{\on{Eis}}
\nc{\Aut}{\on{Aut}}
\nc{\Rep}{\on{Rep}}
\nc{\Hom}{\on{Hom}}
\nc{\Loc}{\on{Loc}}
\nc{\Bun}{\on{Bun}}
\nc{\Mod}{\on{Mod}}
\nc{\IC}{\on{IC}}
\nc{\Can}{\on{Can}}
\nc{\End}{\on{End}}
\nc{\rk}{\on{rk}}
\nc{\sym}{\on{sym}}
\nc{\Sh}{\on{Sh}}
\nc{\Perv}{\on{Perv}}
\nc{\Conv}{\on{Conv}}
\nc{\supp}{\on{supp}}
\nc{\BunBb}{\overline{\Bun}_B}
\nc{\BunNft}{\Bun_N^{\F_T}}
\nc{\BunNftb}{\overline{\Bun}_N^{\F_T}}
\nc{\BunNftw}{\widetilde{\Bun}_N^{\F_T}}
\nc{\BunPb}{\overline{\Bun}_P}
\nc{\BunPbw}{\widetilde{\Bun}_P}
\nc{\GUb}{\overline{G/U}}
\nc{\Psib}{\overline{\Psi}}
\nc{\Psio}{\overset{o}{\Psi}}
\nc{\GUPb}{\overline{G/U(P)}}
\nc{\GPPb}{\overline{G/[P,P]}}
\nc\hl{h^{\leftarrow}}
\nc\hr{h^{\rightarrow}}
\nc\hlb{{\bar h^{\leftarrow}}}
\nc\hrb{{\bar h^{\rightarrow}}}
\nc\hlt{{\tilde h^{\leftarrow}}}
\nc\hrt{{\tilde h^{\rightarrow}}}
\nc\Hl{H^{\leftarrow}}
\nc\Hr{H^{\rightarrow}}

\nc\I{{\mathcal I}}

\nc{\Q}{{\mathcal Q}}
\nc{\F}{{\mathcal F}}
\nc{\A}{{\mathcal A}}
\nc{\W}{{\mathcal W}}
\nc{\G}{{\mathcal G}}
\nc{\J}{{\mathcal J}}
\renewcommand{\P}{{\mathcal P}}
\nc{\Wb}{\overline{\W}}
\nc{\M}{{\mathcal M}}
\nc{\N}{{\mathcal N}}
\nc{\Y}{{\mathcal Y}}
\nc{\E}{{\mathcal E}}
\nc{\D}{{\mathcal D}}
\nc\Dh{\widehat{\D}}
\renewcommand{\O}{{\mathcal O}}
\nc{\C}{{\mathcal C}}
\nc{\K}{\mathcal K}
\renewcommand{\H}{{\mathcal H}}
\renewcommand{\S}{{\mathcal S}}
\nc{\T}{{\mathcal T}}
\renewcommand{\L}{{\mathcal L}}
\nc{\Gr}{\on{Gr}}

\nc{\p}{\mathfrak p}
\nc{\q}{\mathfrak q}

\nc{\s}{{\mathfrak s}}

\nc{\pw}{\widetilde{\mathfrak p}}
\nc{\qw}{\widetilde{\mathfrak q}}
\nc{\jw}{\widetilde j}

\nc{\grb}{\overline{\Gr}}

\nc{\lambdach}{\check\lambda}
\nc{\Lambdach}{\check\Lambda}
\nc{\much}{\check\mu}
\nc{\omegach}{\check\omega}
\nc{\nuch}{\check\nu}
\nc{\etach}{\check\eta}
\nc{\alphach}{\check\alpha}
\nc{\betach}{\check\beta}
\nc{\rhoch}{\check\rho}
\nc\xl{\overline{x}}
\nc\yl{\overline{y}}
\nc\nul{\overline{\nu}}
\nc\mul{\overline{\mu}}
\nc\lambdal{\overline{\lambda}}
\nc\zerol{\overline{0}}

\nc{\Hb}{\overline{\H}}

\nc{\mc}{\mathcal}
\nc{\ga}{\gamma}
\nc{\GL}{{}^L G}
\nc{\PL}{\Lambda}
\nc{\la}{\lambda}
\nc{\arr}{\rightarrow}
\nc{\ol}{\overline}
\nc{\al}{\alpha}
\nc{\De}{\Delta}
\nc{\Gaf}{{\mathbb G}_{a,\Fq}}
\nc{\wt}{\widetilde}
\nc{\bi}{\bibitem}
\nc{\WW}{{\mathbf W}}
\nc{\bn}{\BunNft}
\nc{\map}{\kappa}
\nc{\wh}{\widehat}
\nc{\bnx}{_{x,\infty}\BunNftb}
\nc{\canG}{\mathcal G}
\nc{\canB}{\mathcal B}
\nc{\canN}{{\mathcal N}^\epsilon}
\nc{\CanN}{\wt{\mathcal N}^\epsilon}
\nc{\canNnu}{{\mathcal N}^{\epsilon_\nu}}
\nc{\canNnup}{{\mathcal N}^{\epsilon_{\nu'}}}
\nc{\CanNnu}{\wt{\mathcal N}^{\epsilon_\nu}}
\nc{\CanNnup}{\wt{\mathcal N}^{\epsilon_{\nu'}}}
\nc{\out}{\on{out}}
\nc{\ft}{{\mathcal F}_T}
\nc{\zl}{\overline{z}}
\nc{\lal}{\overline{\la}}
\nc{\BunNftbm}{\ol{\Bun}_{N,\mu}^{\ft}}
\nc{\BunNftm}{\wt{\Bun}^{\F_T}_{N,\mu}}
\nc{\rhoc}{\check{\rho}}
\nc{\tboxtimes}{\widetilde{\boxtimes}}
\nc{\Coh}{{{\mathcal C}oh}}
\nc{\Cohn}{\Coh_n}
\nc{\bfd}{\Fq}
\nc{\bbfd}{\Fqb}

\nc{\V}{{\mc V}}
\nc{\canNK}{\ol{N(K_x)}^{\epsilon_\nu}}
\nc{\jj}{{\mathfrak j}}

\def\Homom#1#2{\operatorname{Hom}(#1,#2)}

\nc{\kk}{\Bbbk}
\nc{\He}{\on{H}}
\nc{\thl}{\wt{h}^{\leftarrow}}
\nc{\thr}{\wt{h}^{\rightarrow}}
\nc{\ovc}{\overset{\circ}}
\nc{\Ext}{\on{Ext}}
\nc{\HQ}{{\mc H}{\mc Q}}
\nc{\HL}{\on{HL}}
\nc{\oS}{\S^0}
\nc{\Hav}{\on{H}}
\nc{\HH}{H}
\nc{\Qp}{\ol\Q_+}
\nc{\Qpp}{\ol\Q_{++}}
\nc{\Qpm}{\ol\Q_{+-}}
\nc{\Mpp}{\Mod_{n,++}^d}
\nc{\Mpm}{\Mod_{n,+-}^d}
\nc{\Ppp}{{\mathcal P}_{E,++}^{d+1}}
\nc{\Ppm}{{\mathcal P}_{E,+-}^{d+1}}
\nc{\hlpm}{{}'\hl_{+-}{}}
\nc{\hrpm}{{}'\hr_{+-}{}}
\nc{\hlpp}{{}'\hl_{++}{}}
\nc{\hrpp}{{}'\hr_{++}{}}
\nc{\Wp}{{\mathcal W}_{E,+}}
\nc{\laa}{\la}
\nc{\Havv}{{\mathbf H}}

\nc{\bb}{\mathbf b}
\nc{\cc}{\mathbf c}

\title{On the geometric Langlands conjecture}

\author{E. Frenkel}

\address{Department of Mathematics, University of California,
Berkeley, CA 94720, USA}

\author{D. Gaitsgory}

\address{Department of Mathematics, Harvard University, Cambridge, MA
02138}

\author{K. Vilonen}

\address{Department of Mathematics, Northwestern University, Evanston, IL
60208, USA}

\date{December 2000; Revised: October 2001}

\begin{document}

\maketitle


\section*{Introduction}

\ssec{Background} Let $X$ be a smooth, complete, geometrically
connected curve over the finite field $\Fq$. Denote by $F$ the field of
rational functions on $X$ and by ${\mathbb A}$ the ring of ad\`eles of
$F$. The Langlands conjecture, recently proved by L.~Lafforgue \cite{Lf},
establishes a correspondence between cuspidal automorphic forms
on the group $GL_n({\mathbb A})$ and irreducible, almost everywhere
unramified, $n$--dimensional $\ell$--adic representations of the
Galois group of $\ol{F}$ over $F$ (more precisely, of the Weil group).

An unramified automorphic form on the group $GL_n({\mathbb A})$ can be
viewed as a function on the set $\Bun_n(\Fq)$ of isomorphism classes
of rank $n$ bundles on the curve $X$. The set $\Bun_n(\Fq)$ is the set
of $\Fq$--points of $\Bun_n$, the algebraic stack of rank $n$ bundles
on $X$. According to Grothendieck's ``faisceaux--fonctions''
correspondence, one can attach to an $\ell$--adic perverse sheaf on
$\Bun_n$ a function on $\Bun_n(\Fq)$ by taking the traces of the
Frobenius on the stalks. V.~Drinfeld's geometric proof \cite{Dr} of
the Langlands conjecture for $GL_2$ (and earlier geometric
interpretation of the abelian class field theory by P.~Deligne, see
\cite{La1}) opened the possibility that automorphic forms may be
constructed as the functions associated to perverse sheaves on
$\Bun_n$.

Thus, one is led to a geometric version of the Langlands conjecture
proposed by V.~Drinfeld and G.~Laumon: for each
geometrically irreducible rank $n$ local system $E$ on $X$ there
exists a perverse sheaf $\Aut_E$ on $\Bun_n$ (irreducible on each
component), which is a Hecke eigensheaf with respect to $E$, in an
appropriate sense (see \cite{La1} or Sect.~1 below for the precise
formulation). Moreover, the geometric Langlands conjecture can be made
over an arbitrary field $\kk$.

Building on the ideas of Drinfeld's work \cite{Dr}, G.~Laumon gave a
conjectural construction of $\Aut_E$ in \cite{La1,La2}. More
precisely, he attached to each rank $n$ local system $E$ on $X$ a
complex of perverse sheaves $\Aut'_E$ on the moduli stack $\Bun'_n$ of
pairs $\{\M,s\}$, where $\M \in \Bun_n$ is a rank $n$ bundle on $X$ and
$s$ is a regular non-zero section of $\M$. He  conjectured that if
$E$ is geometrically irreducible then this sheaf descends to a
perverse sheaf $\Aut_E$ on $\Bun_n$ (irreducible on each component),
which is a Hecke eigensheaf with respect to $E$.

In our previous work \cite{FGKV}, joint with D.~Kazhdan, we have shown
that the function on $\Bun'_n(\Fq)$ associated to $\Aut'_E$ agrees
with the function constructed previously by I.I.~Piatetskii-Shapiro
\cite{PS} and J.A.~Shalika \cite{Sha}, as anticipated by Laumon
\cite{La2}. This provided a consistency check for Laumon's
construction.

In this paper  we formulate a certain vanishing conjecture, and prove
that Laumon's construction indeed produces a perverse sheaf $\Aut_E$
on $\Bun_n$ with desired properties, when the vanishing conjecture holds.
In other words, the vanishing conjecture implies the geometric Langlands
conjecture, over any field $\kk$. For the sake of definiteness, we work
in this paper with a field $\kk$ of characteristic $p>0$, but our results
(with appropriate modifications, such as switching from perverse sheaves
to ${\mathcal D}$--modules) remain valid if $\on{char} \kk = 0$.

Moreover, in the case when $\kk$ is a finite field, we derive the
vanishing conjecture (and hence the geometric Langlands conjecture)
from the results of L.~Lafforgue \cite{Lf}.

To give the reader a feel for the vanishing conjecture, we give here one
of its formulations (see \secref{sect on van} for more details). Let $\M$
and $\M'$ be two vector bundles on $X$ of rank $k$, such that $\deg(\M')
- \deg(\M) = d$. Consider the space
$\Hom^0(\M,\M')$ of injective sheaf homomorphisms $\M \hookrightarrow
\M'$. Let $\Coh^d_0$ be the algebraic stack which classifies torsion
sheaves of length $d$, and $\pi:\Hom^0(\M,\M') \to \Coh^d_0$ the
natural morphism sending $\M \hookrightarrow \M'$ to $\M'/\M$.

G.~Laumon \cite{La1} has defined a remarkable perverse sheaf $\L^d_E$
on $\Coh^d_0$ for any local system $E$ of rank $n$ on $X$. The
vanishing conjecture states that if $E$ is irreducible and $n>k$, then

\hskip2cm $H^\bullet(\Hom^0(\M,\M'),\pi^*(\L^d_E)) = 0, \qquad \forall
d>kn(2g-2).$

\ssec{Contents}
The paper is organized as follows:

\medskip

In Sect. 1 we define Hecke functors and state the geometric Langlands
conjecture. We want to draw the reader's attention to the fact that our
formulation is different from that given in \cite{La1} in two respects.
The Hecke property is defined here using only the first Hecke functor;
according to \propref{thesis}, this implies
the Hecke property with respect to the other Hecke functors. We also do
not require the cuspidality property in the statement of the conjecture,
because we show in \secref{cuspidality} that the cuspidality of a Hecke
eigensheaf follows from the vanishing conjecture.

\smallskip

In Sect. 2 we recall the definition of Laumon's sheaf and state our
vanishing conjecture.

\smallskip

In Sect. 3 we present two constructions of $\Aut_E$ following Laumon
\cite{La1,La2} (see also \cite{FGKV}). A third construction, which
uses the Whittaker sheaves is given in Sect. 4. This construction is
the exact geometric analogue of the construction of Piatetskii-Shapiro
\cite{PS} and Shalika \cite{Sha} at the level of functions. The reader
is referred to \secref{sum three} for a summary of the relationship
between the three constructions and the strategy of our proof.

\smallskip

In Sects. 5--9 we derive the geometric Langlands conjecture assuming
that the vanishing conjecture is true. Sect. 5 is devoted to the proof
of the cleanness property in Laumon's construction. In Sect.~6 we
prove that the sheaf $\Aut'_E$ on $\Bun'_n$ descends to a perverse
sheaf $\Aut_E$ on $\Bun_n$. In Sects. 7 and 8 we give two
alternative proofs of the Hecke property of $\Aut_E$. We then show in
Sect. 9 that the perverse sheaf $\Aut_E$ is cuspidal.

\smallskip

In Sect.~10 we derive the vanishing conjecture from results of
L.~Lafforgue \cite{Lf} when $\kk$ is a finite field.

\smallskip

The Appendix contains proofs of some results concerning the Whittaker
sheaves, which are not necessary for our proof, but are conceptually
important.

\ssec{Notation and conventions}    \label{conventions}

Throughout this paper, $\kk$ will be a ground field of characteristic
$p>0$ and $X$ will be a smooth  complete geometrically connected curve
over $\kk$ of genus $g>1$.

This paper deals with $\Ql$-adic perverse sheaves and complexes of
perverse sheaves on various schemes over $\kk$, where $\ell$ is a
prime with $(\ell,p)=1$. In particular, by a local system on $X$ we
will understand a smooth $\ell$--adic sheaf over $X$. For brevity, we
will refer to a geometrically irreducible local system simply as an
irreducible local system.

When $\kk=\Fq$ we work with Weil sheaves (see \cite{De}), instead of
sheaves defined over $\Fq$.
We choose a square root of $q$ in $\Ql$, which defines a half-integral
Tate twist $\Ql(\frac{1}{2})$.

In addition to $\kk$-schemes, we will extensively use algebraic stacks
in the smooth topology (over $\kk$), see \cite{LMB}. If $G$ is an algebraic
group, we define $\Bun_G$ as a stack that classifies $G$--bundles on
$X$. This means that $\on{Hom}(S,\Bun_G)$ is the groupoid whose
objects are $H$--bundles on $X \times S$ and morphisms are
isomorphisms of these bundles. The pull--back functor for a morphism
$S_1\to S_2$ is defined in a natural way.

When $G=GL_n$, $\Bun_G$ coincides with $\Bun_n$, the moduli
stack of rank $n$ vector bundles on $X$. We write $\Bun_n^d$ for the
connected component of $\Bun_n$ corresponding to rank $n$ vector
bundles of degree $d$.

For an algebraic stack $\Y$ we will use the notation $\on{D}(\Y)$
for the derived category of $\Ql$-adic perverse sheaves on $\Y$.
We refer the reader to Sect.~1.4 of \cite{FGV} for our conventions
regarding this category. When we discuss objects of the derived category,
the cohomological grading should always be understood in the perverse
sense. In addition, for a morphism $f:\Y_1\to \Y_2$, the functors
$f_{!}$, $f_{*}$, $f^*$ and $f^{!}$ should be understood ``in the
derived sense''.

If $\Y$ is a stack over $\kk=\Fq$ and ${\mathbb F}_{q_1}$ is an
extension of $\Fq$, we denote by $\Y({\mathbb F}_{q_1})$ the set of
isomorphism classes of objects of the groupoid $\on{Hom}(\on{Spec}
{\mathbb F}_{q_1},\Y)$. If ${\mathcal S}$ is a perverse sheaf or a
complex of perverse sheaves on $\Y$, then $\Y({\mathbb F}_{q_1})$ is
endowed with the function ``alternating sum of traces of the Frobenius
on stalks'' (as in \cite{De}). We denote this function by $\text{\tt
f}_{q_1}(\S)$.

For the general definitions related to the Langlands correspondence
and the formulation of the Langlands conjecture we refer the reader to
\cite{La1}, Sect.~1 and \cite{FGKV}, Sect.~2. In particular, the
notions of cuspidal automorphic function or Hecke eigenfunction on
$GL_n(\AA)$ may be found there.

\ssec{Acknowledgments} We express our gratitude to D.~Kazhdan for his
collaboration in \cite{FGKV}, which has influenced this work.  We also
thank V.~Drinfeld, D.~Kazhdan, and I.~Mirkovi\'c for valuable
discussions.

\section{Hecke eigensheaves}

In this section we introduce the Hecke functors and state the
geometric Langlands conjecture.

\ssec{Hecke functors}

Consider the following correspondence:

\begin{equation}    \label{diagram H1}
\Bun_n \xleftarrow{\hl} \H^1_n \xrightarrow{\supp\times \hr} X\times
\Bun_n,
\end{equation}
where the stack $\H^1_n$ classifies quadruples
$(x,\M,\M',\beta:\M'\hookrightarrow\M)$, with $x\in X$,
$\M',\M\in\Bun_n$, such that
$\M/\M'$ is the simple skyscraper sheaf supported at $x$, i.e.,
$\M/\M'$ is (non-canonically) isomorphic to
$\OO_X(x)/\OO_X$. The morphisms $\hl$, $\hr$ and $\supp$
are given by $\hl(x,\M,\M')=\M$, $\hr(x,\M,\M')=\M'$, and
$\supp(x,\M,\M')=x$.

\medskip

The {\em Hecke functor} $\He^1_n:\on{D}(\Bun_n) \to
\on{D}(X\times\Bun_n)$ is defined by the formula
\begin{equation}    \label{formula H1}
\He_n^1(\K) = (\supp\times\hr)_! \hl{}^*(\K)\otimes
\Ql(\frac{n-1}{2})[n-1].
\end{equation}

Consider the $i$-th iteration of $\He^1_n$:
\begin{equation*}
(\He^1_n)^{\boxtimes i}:\on{D}(\Bun_n) \to \on{D}(X^i\times\Bun_n)
\end{equation*}
Let $\Delta$ denote the divisor in $X^i$ consisting of the
pairwise diagonals. Note that for any $\K\in \on{D}(\Bun_n)$ the
restriction $(H^1_n)^{\boxtimes i}(\K)|_{(X^i-\Delta)\times \Bun_n}$
is naturally equivariant with respect to the action of the symmetric
group $S_i$ on $X^i-\Delta$.

Consider a  rank $n$ local system $E$ on $X$. We say that $\K\in
\on{D}(\Bun_n)$ is a {\it Hecke eigensheaf}, or that it has a {\it Hecke
property} with respect to $E$, if $\K\neq 0$ and there exists an
isomorphism
\begin{equation} \label{eigen-property}
\He^1_n(\K)\simeq E\boxtimes \K,
\end{equation}
such that the resulting map
\begin{equation} \label{eigen-equivariance}
(\He^1_n)^{\boxtimes 2}(\K)|_{(X\times X-\Delta)\times \Bun_n}\to
E\boxtimes E\boxtimes \K|_{(X\times X-\Delta)\times \Bun_n}
\end{equation}
is $S_2$--equivariant.
This  implies that
$$(\He^1_n)^{\boxtimes i}(\K)|_{(X^i-\Delta)\times \Bun_n}\to
E^{\boxtimes i}\boxtimes \K|_{(X^i-\Delta)\times \Bun_n}$$ is
$S_i$-equivariant for any $i$.

\ssec{Statement of the Geometric Langlands conjecture}

We are now ready to  formulate  the unramified geometric Langlands
conjecture for $GL_n$:

\begin{conj} \label{geometric Langlands}
For each irreducible rank $n$ local system $E$ on $X$ there exists a
perverse sheaf $\Aut_E$ on $\Bun_n$, irreducible on each connected
component $\Bun^d_n$, which is a Hecke eigensheaf with respect to $E$.
\end{conj}

\conjref{geometric Langlands} has been proved by Drinfeld \cite{Dr} in
the case when $n=2$ (see also \cite{Ga}).

In this paper we reduce \conjref{geometric Langlands} to the Vanishing
\conjref{vanishing conjecture}. Then, in \secref{proof of vanishing}
we will show that when $\kk$ is a finite field $\Fq$, the Vanishing
Conjecture follows from recent results of Lafforgue \cite{Lf}.

\ssec{Other Hecke functors}

In addition to the functor $\He^1_n$, we also have Hecke functors
$\He^i_n:\on{D}(\Bun_n) \to \on{D}(X\times\Bun_n)$ for $i=2,...,n$.
To define them, consider the stack $\H^i_n$ which classifies quadruples
$$(x,\M,\M',\beta:\M'\hookrightarrow\M),$$ where $x\in X$,
$\M',\M\in\Bun_n$ such that $\M' \subset \M
\subset \M'(x)$, and $\on{length}(\M/\M')=i$.

As in the case of $\H^1_n$, we have a diagram
\begin{equation*}
\Bun_n \xleftarrow{\hl} \H^i_n \xrightarrow{\supp\times \hr} X\times
\Bun_n,
\end{equation*}
and the functor $\He^i_n$ is defined by the formula
\begin{equation*}
\He^i_n(\K) = (\supp\times\hr)_! \hl{}^*(\K)\otimes
\Ql(\frac{(n-i)i}{2})[(n-i)i].
\end{equation*}

The following result is borrowed from \cite{Ga}:

\begin{prop}  \label{thesis}
Let $\K$ be a Hecke eigensheaf with respect to $E$. Then for $i=1,...,n$
we have isomorphisms $\He^i_n(\K)\simeq \Lambda^i E \boxtimes \K$.
\end{prop}

\begin{proof}
Consider the stack $\Mod_n^{-i}$ of "lower modifications of length
$i$", which classifies the data of triples
$(\M,\M',\beta:\M'\hookrightarrow \M)$, where $\M,\M'\in \Bun_n$ and
$\beta$ is an embedding of coherent sheaves such that the quotient
$\M/\M'$ is a torsion sheaf of length $i$.

Let $X^{(i)}$ be the $i$-th symmetric power of $X$. We have a natural
morphism $\supp:\Mod^{-i}_n\to X^{(i)}$, which associates to
$(\M',\M,\beta)$ as above the divisor of zeros of the induced map
$\det \M' \to \det \M$.

Denote by $\H^{i,+}_n$ the preimage in $\Mod^{-i}_n$ of the main
diagonal $X\subset X^{(i)}$.  Note that $\H^i_n$ is naturally a closed
substack in $\H^{i,+}_n$.

Consider the stack $\wt\Mod{}^{-i}_n$, which classifies the data
$(\M'=\M_0\subset \M_1\subset...\subset\M_i=\M)$, where each $\M_j$ is
a rank $n$ vector bundle, and $\M_j/\M_{j-1}$ is a simple skyscraper
sheaf. There is a natural proper map $p:\wt\Mod^{-i}_n\to \Mod^{-i}_n$,
which "forgets" the middle terms of the filtration.

There is also a natural map $\wt\supp:\wt\Mod{}^{-i}_n\to X^i$
such that if $\on{sym}:X^i\to X^{(i)}$ denotes the symmetrization map,
we have $\sym\circ \,\wt\supp=\supp\circ\, p:\wt\Mod^{-i}_n\to X^{(i)}$.

The open substack $\wt\supp^{-1}(X^{i}-\Delta)$ of $\wt\Mod{}^{-i}_n$ is
isomorphic to the fiber product $\Mod^{-i}_n \underset{X^{(i)}}\times
(X^{i}-\Delta)$.

The map $p$ is known to be small (see, e.g., \cite{La1}).
This implies that the complex
$${\mc Spr}:=p_!(\Ql(\frac{i(n-1)}{2}))[i(n-1)]$$
on $\Mod^{-i}_n$ is perverse
(up to the cohomological shift by $n^2\cdot (g-1)=\dim(\Bun_n)$) and is
a Goresky-MacPherson extension of its
restriction to $\supp^{-1}(X^{(i)}-\Delta)$. In particular, ${\mc Spr}$
carries a canonical $S_i$-action and $({\mc Spr})^{S_i}\simeq
\Ql(\frac{i(n-1)}{2})[i(n-1)]$.

Let $\hl$ (resp., $\hr$) denote the morphism $\Mod_n^{-i}\to\Bun_n$,
which sends a triple $(\M,\M',\beta)$ to $\M$ (resp., $\M'$). By
construction, for any $\K\in \on{D}(\Bun_n)$,
\begin{equation*}
(\supp\times \hr)_!(\hl{}^*(\K)\otimes {\mc Spr})
\simeq (\on{sym}\times\on{id})_! (\He^1_n)^{\boxtimes i}(\K).
\end{equation*}
Thus, if $\K$ is a Hecke eigensheaf with respect to $E$, we obtain
an $S_i$-equivariant isomorphism
\begin{equation}   \label{Hecke fundamental}
(\supp\times \hr)_!(\hl{}^*(\K)\otimes {\mc Spr})
\simeq \on{sym}_!(E^{\boxtimes i})\boxtimes \K.
\end{equation}

To conclude the proof, we pass to the isotypic components of the
$sign$ representation of
$S_i$ on both sides of formula \eqref{Hecke fundamental} and  restrict
the resulting isotypic components to the main diagonal $X\subset X^{(i)}$.
By this process the RHS of \eqref{Hecke fundamental}  tautologically
yields $\Lambda^i E \boxtimes \K$. Thus, it remains to show that  the LHS
yields $\He^i_n(\K)$. To see this, it suffices to note  that,
$\Hom_{S_i}(sign,{\mc Spr})|_{\H^{i,+}_n}$ is isomorphic to the constant
sheaf on $\H^i_n$ tensored by $\Ql(\frac{(n-i)i}{2})[(n-i)i]$, by
the Springer theory \cite{BM,Sp}.

\end{proof}

The isomorphisms constructed in the above proposition have an
additional property. To state it, let $\sigma$ be the transposition acting
on $X\times X$ and let $i,j\in \{1,...,n\}$. Clearly, the functors
$$\K\mapsto (\He^i_n\times\on{id})\circ \He^j_n(\K)|_{(X\times
X-\Delta)\times \Bun_n}$$
and
$$\K\mapsto \sigma^*\circ (\He^j_n\times\on{id})\circ
\He^i_n(\K)|_{(X\times X-\Delta)\times \Bun_n}$$ from $\on{D}(\Bun_n)$
to $\on{D}((X\times X-\Delta)\times\Bun_n)$ are naturally isomorphic.
Hence, for a Hecke eigensheaf $\K$, the following diagram is
commutative:
\begin{equation}  \label{Hecke commute}
\begin{CD}
(\He^i_n\times\on{id})\circ \He^j_n(\K)|_{X\times X-\Delta}   @>>>
\sigma^*\circ (\He^j_n\times\on{id})\circ \He^i_n(\K)|_{X\times
X-\Delta}  \\
@VVV    @VVV   \\
\Lambda^i E \boxtimes \Lambda^j E \boxtimes \K|_{X\times X-\Delta}
@>>>\sigma^*(\Lambda^j E\boxtimes
\Lambda^i E)\boxtimes \K|_{X\times X-\Delta}.
\end{CD}
\end{equation}

Finally, let us consider the Hecke functor $\He^n_n$. By
definition, this is the pull-back under the morphism
$\on{mult}:X\times \Bun_n\to \Bun_n$ given by $(x,\M)\mapsto \M(x)$.
Hence if $\K$ is a Hecke eigensheaf with respect to $E$, then
\begin{equation*}
\on{mult}^*(\K)\simeq \Lambda^n E \boxtimes \K.
\end{equation*}

\section{The vanishing conjecture}    \label{sect on van}

Denote by $\Coh_n$ the stack classifying coherent sheaves on $X$ of
generic rank $n$. More precisely, for each $\kk$--scheme $S$,
$\on{Hom}(S,\Coh_n)$ is the groupoid, whose objects are coherent
sheaves $\M_S$ on $X \times S$, which are flat over $S$, and such that
over every geometric point $s\in S$, $\M_s$ is generically of rank
$n$. We write $\Coh_n^d$ for the substack corresponding to coherent
sheaves of generic rank $n$ and degree $d$.

\ssec{Laumon's sheaf}   \label{Laumon's sheaf}

In \cite{La2} Laumon associated to an arbitrary local system $E$ of
rank $n$ on $X$ a perverse sheaf ${\mathcal L}_E$ on $\Coh_0$. Let us
recall his construction. Denote by $\Coh_0^{\on{rss}}$ the open
substack of $\Coh_0$ corresponding to regular semisimple torsion
sheaves. Thus, a geometric point of $\Coh_0$ belongs to
$\Coh_0^{\on{rss}}$ if the corresponding coherent sheaf on $X$ is a
direct sum of skyscraper sheaves of length one supported at distinct
points of $X$. Let $\Coh_0^{\on{rss},d}=\Coh_0^{\on{rss}}\cap
\Coh^d$. We have a natural smooth map $(X^{(d)}-\Delta)\to
\Coh_0^{\on{rss},d}$.

Let $E^{(d)}$ be the $d$-th symmetric power of $E$, i.e., $E^{(d)} =
\on{sym}_!(E^{\boxtimes d})^{S_d}$, where $\on{sym}: X^d \to
X^{(d)}$. This is a perverse sheaf on $X^{(d)}$, and its restriction
$E^{(d)}|_{X^{(d)}-\Delta}$ is a local system, which descends to a
local system $\ovc\L{}_E^d$ on $\Coh_0^{\on{rss},d}$. The perverse
sheaf $\L^d_E$ on $\Coh^d_0$ is by definition the Goresky-MacPherson
extension of $\ovc\L{}_E^d$ from $\Coh_0^{\on{rss},d}$ to
$\Coh^d_0$. We denote by $\L_E$ the perverse sheaf on $\Coh_0$, whose
restriction to $\Coh^d_0$ equals $\L_E^d$.

\ssec{The averaging functor}    \label{averaging}

Using the perverse sheaf $\L_E^d$ we define the averaging functor
$\Hav^d_{k,E}:\on{D}(\Bun_k)\to \on{D}(\Bun_k)$. We stress that the
positive integer $k$ is independent of $n$, the rank of the local
system $E$.

For $d\geq 0$, introduce the stack $\Mod^d_k$, which classifies the
data of triples $(\M,\M',\beta:\M\hookrightarrow \M')$, where
$\M,\M'\in \Bun_k$ and $\beta$ is an embedding of coherent sheaves
such that the quotient $\M'/\M$ is a torsion sheaf of length $d$, and
the diagram
\begin{equation*}
\Bun_k \overset{\hl}\longleftarrow \Mod^d_k
\overset{\hr}\longrightarrow \Bun_k,
\end{equation*}
where $\hl$ (resp., $\hr$) denotes the morphism sending a triple
$(\M,\M',\beta)$ to $\M$ (resp., $\M'$). In addition, we have a
natural smooth morphism $\pi:\Mod^d_k\to \Coh^d_0$, which sends a
triple $(\M,\M',\beta)$ to the torsion sheaf $\M'/\M$.

Note that $\Mod^d_n$ is isomorphic to the stack $\Mod^{-d}_n$
which was used in the proof of \thmref{thesis}. Under this isomorphism
the maps $\hr$ and $\hl$ are reversed.

The {\em averaging functor} $\Hav^d_{k,E}: \on{D}(\Bun_k) \to
\on{D}(\Bun_k)$ is defined by the formula
\begin{equation*}
\K \mapsto \hr_!(\hl{}^*(\K)\otimes \pi^*(\L^d_E))\otimes
\Ql(\frac{d\cdot k}{2})[d\cdot k].
\end{equation*}

\ssec{Vanishing Conjecture} \label{vanishing conjecture} {\em Assume
that $E$ is an irreducible local system of rank $n$. Then for all
$k=1,\ldots,n-1$ and all $d$ satisfying $d > kn(2g-2)$, the functor
$\Hav^d_{k,E}$ is identically equal to $0$.}

The statement of the Vanishing Conjecture is known to be true for
$k=1$ (see below).

The goal of this paper is to show that if \conjref{vanishing
conjecture} holds for {\em any given} irreducible rank $n$ local
system $E$, then the geometric Langlands \conjref{geometric Langlands}
holds for $E$. In addition, in \secref{proof of vanishing} we will
prove \conjref{vanishing conjecture} in the case when $\kk$ is a
finite field $\Fq$ if the following statements are true (see
\cite{BBD} for the definition of a pure local system):

\smallskip

\noindent {(a)} $E$ is pure (up to a twist by a one-dimensional
representation of the Weil group of $\Fq$),

\smallskip

\noindent and either

\smallskip

\noindent {(b)} there exists a cuspidal Hecke eigenfunction associated to
the pull-back of $E$ to $X \underset{\Fq}\times {\mathbb F}_{q_1}$ for any
finite extension ${\mathbb F}_{q_1}$ of $\Fq$;

\smallskip

\noindent or

\smallskip

\noindent {(b')} the space of unramified cuspidal automorphic
functions on the group $GL_k$ over the ad\`eles is spanned by the
Hecke eigenfunctions associated to rank $k$ local systems on $X
\underset{\Fq}\times {\mathbb F}_{q_1}$, for all $k<n$.

The statements (a),(b),(b') follow from the recent work of Lafforgue
\cite{Lf} (note that (b) and (b') are specified by the Langlands
conjecture at the level of functions). Therefore, Lafforgue's results
together with the results of the present paper, imply
\conjref{vanishing conjecture} and hence \conjref{geometric Langlands}
over a finite field $\kk$.

\ssec{A reformulation of the Vanishing Conjecture}

Let $\M$ and $\M'$ be two rank $k$ vector bundles on $X$ and let us
write $\on{Hom}^0(\M,\M')$ for the open subset of injective maps in
the vector space $\Hom(\M,\M')$. There is a natural morphism $\pi:
\on{Hom}^0(\M',\M) \to \Coh^d_0$, where $d=\deg(\M')-\deg(\M)$, which
maps $(\M \hookrightarrow \M')$ to $\M'/\M$.

\begin{conj} \label{pointwise vanishing}
Under the assumptions on $E$ and $d$ given in \conjref{vanishing
conjecture},
$$\HH^{\bullet}(\on{Hom}^0(\M,\M'), \pi^*(\L^d_E))=0.$$
\end{conj}

Conjectures \ref{vanishing conjecture} and \ref{pointwise vanishing}
are equivalent. Indeed, consider the complex
\begin{equation*}
(\hl\times\hr)_! \pi^*(\L^d_E)\in \on{D}(\Bun_k\times \Bun_k).
\end{equation*}
\conjref{vanishing conjecture} is equivalent to the statement that
this complex equals $0$. But its fiber at $(\M,\M')\in
\Bun_k\times \Bun_k$ is precisely the cohomology
$\HH^{\bullet}(\on{Hom}^0(\M',\M), \pi^*(\L^d_E))$.

\ssec{Proof of the statement of \conjref{vanishing conjecture} in the
case $k=1$}
\label{deligne van thm}

Recall the Deligne vanishing theorem (see the Appendix of \cite{Dr}):

\smallskip

{\em Let $\on{AJ}: X^{(d)}\to \on{Pic}^d(X)$ be the Abel-Jacobi map,
and $E$ an irreducible local system of rank $n>1$. Then for}
$d>n(2g-2)$
\begin{equation*}
(\on{AJ})_!(E^{(d)})=0.
\end{equation*}

\smallskip

This theorem implies the case $k=1$ of \conjref{vanishing conjecture}.
Indeed, consider the morphism $X^{(d)}\to\Coh_0^d$ that associates to
a divisor $D$ the torsion sheaf $\O_X(D)/\O_X$. This morphism is
smooth and its image is the open substack $\Coh_0^{\on{r},d}$ of
$\Coh_0^d$ corresponding to those torsion sheaves $\T$ on $X$ for which
$\dim_\kk(\End(\T))=\on{length}(\T)$ (such torsion sheaves are called
regular).  Clearly, $\Coh_0^{\on{rss},d}\subset \Coh_0^{\on{r},d}$. Since
the Laumon sheaf $\L^d_E$ is an irreducible perverse sheaf on
$\Coh^d_0$, and $E^{(d)}$ is an irreducible perverse sheaf on
$X^{(d)}$, we obtain that the pull-back of $\L^d_E$ under the morphism
$X^{(d)}\to\Coh_0^d$ is isomorphic to $E^{(d)}$.  Observe now that the
diagram of stacks
\begin{equation*}
\Bun_1 \leftarrow \Mod_1^d \rightarrow \Bun_1
\end{equation*}
may be identified with
\begin{equation*}
\on{Pic}(X) \leftarrow \on{Pic}(X)\times X^{(d)} \rightarrow \on{Pic}(X),
\end{equation*}
where the left arrow is the projection on the first factor and the
right arrow is the composition
\begin{equation*}
\on{Pic}(X)\times X^{(d)} \xrightarrow{\on{id}\times \on{AJ}}
\on{Pic}(X)\times \on{Pic}^{d}(X)\xrightarrow{\on{mult}} \on{Pic}(X).
\end{equation*}

Therefore, for $\K\in \on{D}(\on{Pic}(X))$ we have:
\begin{equation*}
\hr_!(\hl{}^*(\K)\otimes \pi^*(\L^d_E))\simeq
\on{mult}_!(\K\boxtimes (\on{AJ})_!(E^{(d)}))=0,
\end{equation*}
by Deligne's theorem.

\section{The construction of $\Aut_E$} \label{construction}

Let $\Coh'_n$ denote the stack classifying pairs $(\M,s)$, where
$\M\in\Coh_n$ and $s$ is an injective map $\Omega^{n-1} \to \M$. Here
$\Omega$ stands for the canonical bundle of $X$ and we write
$\Omega^k$ for $\Omega^{\otimes k}$. We denote by $\Bun'_n$ the
preimage in $\Coh'_n$ of the open substack $\Bun_n\subset \Coh_n$.
Let $\varrho_n: \Coh'_n \to \Coh_n$ be the forgetful map; we use the same
notation for the forgetful map $\Bun'_n \to \Bun_n$.

In this section, starting with a local system $E$ on $X$ of an
arbitrary rank, we will construct a complex $\S'_E$ on
$\Coh'_n$. Later we will show that if $E$ is an irreducible local
system of rank $n$ which satisfies \conjref{vanishing conjecture},
then $\S'_E$ descends to a perverse sheaf $\S_E$
on $\Coh_n$. The restriction of $\S_E$ to $\Bun_n$ will then be the
Hecke eigensheaf $\Aut_E$. We present below three constructions of
$\S'_E$ (two of them in this section, and one more in the next
section).

\ssec{The first construction}    \label{the first}

The following is a version of the construction presented in
\cite{La2,FGKV}.

Define an algebraic stack $\wt{\mathcal Q}$ as follows. For a
$\kk$--scheme $S$, $\on{Hom}(S,\wt{\mathcal Q})$ is the groupoid, whose
objects are quadruples $(\M_S,\beta_S,(\M^0_{i,S}),(\wt
s_{i,S}))$, where $\M_S$ is a coherent sheaf on $X \times S$ of
generic rank $n$, $\M^0_S$ is a rank $n$ bundle on $X \times S$,
$\beta_S: \M^0_S \arr \M_S$ is an embedding of the corresponding
$\OO_{X \times S}$--modules, such that the quotient is $S$--flat,
$(\M^0_{i,S})$ is a full flag of subbundles
\begin{equation} \label{J}
0=\M^0_{0,S} \subset \M^0_{1,S} \subset \ldots \M^0_{n-1,S} \subset
\M^0_{n,s}=\M^0_S,
\end{equation}
and $\wt s_{i,S}$ is an isomorphism ${\Omega}^{n-i} \boxtimes \OO_S
\simeq \M^0_{i,S}/\M^0_{i-1,S}, i=1,\ldots,n$. The morphisms are the
isomorphisms of the corresponding $\OO_{X\times S}$--modules making
all diagrams commutative (we remark that in \cite{FGKV} we used the
notation $J$ instead of $\M^0$).

There is a representable morphism of stacks $\wt{\nu}:\wt{\mathcal
Q}\to\Coh'_n$, which for each $\kk$--scheme $S$ maps
$(\M_S,\beta_S,(\M^0_{i,S}),(\wt s_{i,S}))$ to the pair
$(\M_S,\beta_S \circ \wt s_{1,S})$, where $s_{1,S}$ is viewed as an
embedding of ${\Omega}^{n-1} \boxtimes \OO_S$ into $\M^0_S$.

We also define the morphism $\alpha: \wt{\mathcal Q}\to \Coh_0$
sending $(\M_S,\beta_S,(\M^0_{i,S}),(\wt s_{i,S}))$ to the
sheaf $\M_S/\on{Im} \beta_S$, and the morphism $\on{ev}: \wt{\mathcal
Q} \to \GG_a$ defined as follows.

Given two coherent sheaves $\L$ and $\L'$ on $X$, consider the stack
${\mathcal E}xt^1(\L',\L)$. The objects of the groupoid
$\on{Hom}(S,{\mathcal E}xt^1(\L',\L))$ are coherent sheaves $\L''$ on
$X \times S$ together with a short exact sequence
\begin{equation*}
0 \arr \L \boxtimes
\OO_S \arr \L'' \arr \L' \boxtimes \OO_S \arr 0,
\end{equation*}
and morphisms are maps between such exact sequences inducing the
identity isomorphisms at the ends. There is a canonical morphism
from the stack ${\mathcal E}xt^1(\L',\L)$ to the scheme $\on{Ext}^1(\L',\L)$.
We have for
each $i=1,\ldots,n-1$, a natural morphism $\on{ev}_i: \wt{\mathcal Q}
\arr {\mathcal E}xt^1(\Omega^{i},\Omega^{i-1})$, which sends
the data of $(\M_S,\beta_S,(\M^0_{i,S}),(\wt s_{i,S}))$ to
$$0\to \M^0_{i,S}/\M^0_{i-1,S}\to \M^0_{i+1,S}/\M^0_{i-1,S}\to
\M^0_{i+1,S}/\M^0_{i,S}\to 0.$$

Now $\on{ev}$ is the composition
\begin{equation} \label{ev}
\on{ev}:\wt{\mathcal Q} \arr \prod_{i=1}^{n-1} {\mathcal
E}xt^1(\Omega^{i},\Omega^{ i-1}) \arr \prod_{i=1}^{n-1}
\on{Ext}^1(\Omega^{ i},\Omega^{ i-1}) \arr \GG_a^{n-1}
\overset{\on{sum}}\arr \GG_a.
\end{equation}

We fix a non-trivial character $\psi:\Fq\to \Ql$, which gives rise to
the Artin-Shreier sheaf $\I_\psi$ on the additive group
$\GG_a$. Define the complex $\widetilde{\mathcal W}_E$ on $\wt\Q$ by
the formula
\begin{equation*}
\widetilde{\mathcal W}_E:=\alpha^*({\mathcal L}_E)\otimes
\on{ev}^*({\mathcal I}_\psi)\otimes \Ql(\frac{\dim}{2})[\dim],
\end{equation*}
where $\dim$ is the dimension of the corresponding connected component of
$\wt\Q$.

Since the morphism $\alpha$ is smooth, $\widetilde{\mathcal W}_E$ is a
perverse sheaf. Finally, we define
the complex $\S_E'$ on $\Coh'_n$ by
\begin{equation*}
\S'_E:=\wt{\nu}_!(\widetilde{\mathcal W}_E).
\end{equation*}

\ssec{The second construction, via Fourier transforms}
\label{seccon}  \label{notat}

This construction is due to Laumon \cite{La2}. It amounts to
expressing the first construction as a series of Fourier
transforms. Thus, we obtain an alternative construction of the
restriction of $\S'_E$ to the preimage in $\Coh'_n$ of a certain
open substack $\C_n$ of
$\Coh_n$. For technical reasons, which will become clear in the course
of the proof, we choose a slightly smaller open subset of $\Bun_n$
than in \cite{La2}.

\medskip

We will need the following result:

\medskip

We call a vector bundle $\M$ {\it very unstable} if $\M$ can be
decomposed into a direct sum $\M \simeq \M_1\oplus \M_2$, such that
$\M_i\neq 0$ and $\on{Ext}^1(\M_1,\M_2)=0$.  It is clear that very
unstable vector bundles form a constructible subset
$\Bun_n^{\on{vuns}}$ of $\Bun_n$.

Let $\L^{\on{est}}$ be any fixed line bundle.

\begin{lem} \label{estimate}
There exists an integer $c_{g,n}$ with the following
property: if $d\geq c_{g,n}$ and $\M \in
\Bun^d_n(\ol{\kk})$ is such that
$\Hom(\M,\L^{\on{est}}) \neq 0$, then $\M$ is very unstable.
\end{lem}

\begin{proof}
We will prove a slightly stronger statement. Namely, for each $n$ we
will find an integer $c_{g,n}$ such that for $\M\in \Bun^d_n$, $d\geq
c_{g,n}$, with $\Hom({\mathcal M},\L^{\on{est}}) \neq 0$ there exists a
decomposition $\M\simeq \M_1\oplus \M_2$ with $\M_i\neq 0$,
$\on{Ext}^1(\M_1,\M_2)=0$, and $\frac{\deg(\M_2)}{\on{rank}(\M_2)}\geq
\frac{\deg(\M)}{n}$.

Set
$c_{g,1}=\on{deg}(\L^{\on{est}})$. By induction, we can assume that $c_{g,i}$,
$i\leq n-1$, satisfying the above properties have been found. Let us show
that any integer $c_{g,n}$ such that
\begin{equation}  \label{c estimate}
(c_{g,n} -d^{\on{est}}-(n-1)(2g-2))\cdot \frac{i}{n-1} > c_{g,i},\qquad
\forall i=1,...,n-1,
\end{equation}
will do.

Indeed, let $c_{g,n}$ be such an integer. Suppose that for some $\M\in
\Bun^d_n$, $d\geq c_{g,n}$, we have $\Hom(\M,\L^{\on{est}}) \neq 0$.
Then there exists a short exact sequence
$$0\to \M'\to\M\to \L'\to 0$$ where $\L'$ is a line bundle such that
$\Hom(\L',\L^{\on{est}})\neq 0$.  By \eqref{c estimate} we have:
$\frac{\deg(\M')}{n-1}\geq \frac{\deg(\M)}{n}$. Hence, if
$\on{Ext}^1(\L',\M')$ vanishes,
the decomposition $\M\simeq \M'\oplus \L'$ satisfies our requirements
and we are done. Thus,
it remains to consider the case $\on{Ext}^1(\L',\M')\neq 0$. Then, by
Serre duality,
we obtain that $\Hom(\M'\otimes \Omega^{-1},\L')\neq 0$ and hence
$\Hom(\M'\otimes \Omega^{-1},\L^{\on{est}}) \neq 0$. From the
definition of $\M'$ and
\eqref{c estimate} we conclude that $\deg(\M'\otimes \Omega^{-1}) >
c_{g,n-1}$. We
first observe that this forces $n>2$. For $n=2$ we get
$\deg(\M'\otimes
\Omega^{-1})>\deg(\L^{\on{est}})$, forcing $\Hom(\M'\otimes
\Omega^{-1},\L^{\on{est}})$ to vanish, a contradiction.

Using our induction hypothesis, we can find a direct sum
decomposition $\M'\otimes
\Omega^{-1}\simeq \M'_1\oplus \M'_2$ with $\on{Ext}^1(\M'_1,\M'_2)=0$
and $\frac{\deg(\M'_2)}{i}\geq
\frac{\deg(\M'\otimes\Omega^{-1})}{n-1}$, where
$i=\on{rank}(\M'_2)$. Moreover, without loss of generality we can
assume that $\M'_2$ admits no further decomposition satisfying the
above condition (indeed, if it does, we simply split $\M'_2$
further). Since $\deg(\M'_2)\geq \deg(\M'\otimes \Omega^{-1})\cdot
\frac{i}{n-1}\geq c_{g,i}$, $\Hom(\M'_2,\L^{\on{est}})$ must vanish by
the induction hypothesis.

By  Serre duality, it follows that
$\on{Ext}^1(\L^{\on{est}},\M'_2\otimes\Omega)=0$, and hence
$\M_2:=\M'_2\otimes\Omega$ is a direct summand in $\M$. More
precisely, $\M \simeq \M_2 \oplus \M_1$, where $\M_1$ fits into a
short exact sequence
$$0\to \M'_1\otimes \Omega\to \M_1 \to \L'\to 0.$$ Therefore,
$\on{Ext}^1(\M_1,\M_2)=0$ and
$$\frac{\deg(\M_2)}{\on{rank}(\M_2)}\geq \frac{\deg(\M')}{n-1}\geq
\frac{\deg(\M)}{n}.$$
This completes the proof.
\end{proof}

\noindent{\bf Notational convention.} For notational convenience, in
what follows by {\em degree} of a coherent sheaf of generic rank $k$
we will understand its usual degree $- k(k-1)(g-1)$, so that the
bundle $\O \oplus \Omega \oplus \ldots \oplus \Omega^{k-1}$ is of
degree zero.

\bigskip

To define  $\C_n$, we choose the line bundle
$\L^{\on{est}}$ of a sufficiently large degree such that for any bundle
$\M$ on $X$ of rank $k\leq n$,
$\Hom(\M,\L^{\on{est}})=0$ implies that

\smallskip

\noindent{(a)} $\deg(\M)>n k (2g-2)$,

\smallskip

\noindent{(b)} $\Ext^1(\Omega^{k-1},\M)=0$.

\medskip

For example, any line bundle $\L^{\on{est}}$ of degree $>(2n+2)(g-1)$
will do.

\bigskip

Thus, let $c_{g,n}$ be an integer satisfying the requirements of
\lemref{estimate}. For $d \geq c_{g,n}$, let $\C^d_k$ be the open substack
of $\Coh^d_k$ consisting of $\M\in\Coh^d_k$ such that $\Hom({\mathcal
M},\L^{\on{est}})=0$. Finally, we set $\C_k = \cup_{d \geq c_{g,n}}
\C^d_k$.

Note that by construction any $\M \in
(\Bun^d_n-\C^d_n\cap\Bun^d_n)(\ol{\kk})$, for $d\geq c_{g,n}$, is very
unstable. This property of $\C_n$ will be crucial in \secref{descent}.

\ssec{The fundamental diagram}    \label{fund diag}

Let
\begin{align*}
{\mathcal E}_k&= \text{the stack classifying pairs} \ (\M_k,s_k),
{\mathcal M}_k\in\C_k, s_k\in\Homom{\Omega^{k-1}}{{\mathcal
M}_k} \\ {\mathcal E}_k^{\vee}&= \text{the stack classifying
extensions} \ \ 0\rightarrow\Omega^{ k}\rightarrow{\mathcal
M}_{k+1}\rightarrow{\mathcal M}_k\rightarrow 0, \text{with} \ {\mathcal
M}_k\in\C_k\,.
\end{align*}

We have natural projections $\rho_k: \E_k \to\C_k$ and $\rho_k^{\vee}:
\E_k^{\vee}\to \C_k$, which form dual vector bundles over $\C_k$, due
to the above conditions on $\L^{\on{est}}$. We have: $\rho_k =
\varrho_k|_{\E_k}$.

Next, we set:
\begin{align}
\E_k^0 &= \{(\M_k,s_k) \in \E_k \mid \, s_k \text{ is injective}\}
\subset\E_k \\ \E_k^{\vee 0}&= \{(0\rightarrow\Omega^{
k}\rightarrow\mathcal M_{k+1}\rightarrow\M_k\rightarrow 0) \in
\E_k^{\vee} \mid \M_{k+1}\in \C_{k+1} \} \subset \E_k^{\vee}\,.
\end{align}
Clearly, $\E_k^0 \simeq \E_{k-1}^{\vee 0}$.  Denote by $j_k$ the
embedding $\E^0_k \hookrightarrow \E_k$.
Note also that $\E^0_k$ is an open substack in $\Coh'_n$.

\bigskip

Consider the following diagram:


$$\begin{array}{cccccccccccccccccc}
&  & \E_n & \stackrel{j_n}{\hookleftarrow} & \E^0_n &
\simeq & \E^{\vee 0}_{n-1} & \hookrightarrow & \E^\vee_{n-1} &  &
&  & \E_{n-1} &  &  &  & \\
& \stackrel{\rho_n}{\swarrow} & & & & & & & &
\stackrel{\rho^\vee_{n-1}}{\searrow} & &
\stackrel{\rho_{n-1}}{\swarrow} & & \ldots & & & \\
\C_n &  &  &  &  &  &  &  &  &  & \C_{n-1} &  &  &  &  &  &
\end{array}
$$
$$
\begin{array}{cccccccccccccccccccccc}
& & & & & & & \E_1^\vee & & & & \E_1 & \stackrel{j_1}{\hookleftarrow}
& \E_1^0 & \simeq & \E^{\vee 0}_0 & \hookrightarrow & \E^\vee_0 & & \\
& & & & & & \ldots & & \stackrel{\rho^\vee_1}{\searrow} & &
\stackrel{\rho_1}{\swarrow} & & & & & & & &
\stackrel{\rho^\vee_0}{\searrow} & \\ & & & & & & & & & \C_1 & & & & &
& & & & & \C_0
\end{array}$$

\bigskip

We set $\F_{E,1}$ to be a complex on $\E_1^0$ equal to the pull-back
of Laumon's sheaf $\L_E$ under
$$\E_1^0 \simeq \E^{\vee 0}_0\hookrightarrow
\E^\vee_0\overset{\rho_0^\vee}\longrightarrow \C_0\simeq \Coh_0.$$

Since $\rho_0^\vee$ is a smooth morphism and $\E_1^0\to \E_0^\vee$ is
an open embedding, the restriction of $\F_{E,1}[d]$ to the connected
component of $\E_1^0$ corresponding
to coherent sheaves of degree $d$ is a perverse sheaf.

\medskip

Next, we define the complexes $\F_{E,k}$ on $\E_{k}^0$ by the formula:
\begin{equation*}
\F_{E,k+1}=\on{Four}\,(j_k{}_!(\F_{E,k}))\,|_{\E_{k+1}^0},
\end{equation*}
where $\on{Four}$ is the Fourier transform functor.

Unraveling the second construction we obtain (see \cite{La2}):

\begin{lem}  \label{coincide}
The complex $\F_{E,n}$ coincides (up to a cohomological shift and
Tate's twist) with the restriction of $\S_E'$ to $\E_n^0\subset
\Coh'_n$.
\end{lem}

\ssec{The cleanness property of $\F_{E,k}$}

Let us now assume that $E$ is an irreducible rank $n$ local system and
that \conjref{vanishing conjecture} holds for $E$.

In \secref{cleanness} we prove the following theorem, which was
conjectured by Laumon in \cite{La2}, Expos\'e I, Conjecture 3.2.

\begin{thm}  \label{clean}
For $k=1,...,n-1$, the canonical maps $j_{k!}(\F_{E,k}) \to
j_{k*}(\F_{E,k})$ are isomorphisms.
\end{thm}

Recall that a complex $\K$ on $Y$ is called {\it clean} with
respect to an embedding $Y\overset{j}\hookrightarrow \ol Y$ if
$j_!(\K)\to j_*(\K)$ is an isomorphism, i.e. $j_{*}(\K)|_{\ol
Y-Y}=0$. When $\K$ is a perverse sheaf, cleanness implies that
$j_!(\K) \simeq j_{!*}(\K) \simeq j_{*}(\K)$.
In this language \thmref{clean} states that the sheaf $\F_{E,k}$ on
$\E^0_k$ is clean with respect to $j_k: \E^0_k \hookrightarrow \E_k$.

By construction, Laumon's sheaf $\L_E^d$ is perverse and irreducible.
As was mentioned earlier, the restriction of $\F_{E,1}$ to each
connected component of $\E_1^0$ is, therefore, also an irreducible
perverse sheaf, up to a cohomological shift.  Since the Fourier
transform functor preserves perversity and irreducibility, we obtain
by induction:

\begin{cor}    \label{clean1}
The restriction of $\F_{E,n}$ to each connected component of $\E^0_n$
is an irreducible perverse sheaf, up to a cohomological shift.
\end{cor}

In \secref{descent} we will derive from \corref{clean1} the following
theorem, which was conjectured by Laumon in \cite{La2}, Expos\'e I,
Conjecture 3.1. Denote by $\rho_n^0$ the morphism $\E_{n}^0 \to \C_n$
obtained by restriction from $\rho_n$.

\begin{thm}  \label{descends}
The complex $\F_{E,n}$ descends to $\C_n$, i.e., there exists
a perverse sheaf $\oS_E$ on $\C_n$, such that
$$\F_{E,n}\otimes \Ql(\frac{n^2\cdot (g-1)}{2})[n^2\cdot (g-1)]\simeq
\rho^{0*}_n(\oS_E).$$ Moreover, the restriction of $\oS_E$ to each
connected component of $\C_n$ is a non-zero irreducible perverse
sheaf.
\end{thm}

\ssec{A summary}         \label{sum three}

Above we have described two constructions of a Hecke eigensheaf
associated to a local system $E$. In the next section we will describe
the third construction. Before doing that, we wish to summarize the
relations between the three constructions and to indicate the strategy
of the proof of the main result of this paper that the Vanishing
\conjref{vanishing conjecture} for $E$ implies the Geometric Langlands
\conjref{geometric Langlands}.

In the first construction given in \secref{the first} we constructed a
sheaf ${\mathcal S}'_E$ on the moduli stack $\Coh_n'$ of pairs
$(\M,\Omega^{n-1} \hookrightarrow \M)$, where $\M$ is a coherent sheaf
on $X$ of generic rank $n$.

In the second construction given in \secref{seccon} we defined a sheaf
${\mathcal F}_{E,n}$ on an open substack $\E^0_n$ of $\Coh_n'$ (recall
that this open substack is the preimage of $\C_n\subset \Coh_n$
under $\varrho_n:\Coh'_n\to \Coh_n$).

Finally, in the third construction given in the next section we
produce a sheaf $\Aut_E'$ on the stack $\Bun_n'$ of pairs
$(\M,\Omega^{n-1} \hookrightarrow \M)$, where $\M$ is a rank $n$
vector bundle on $X$ (obviously, $\Bun_n'$ is the preimage
$\varrho_n^{-1}(\Bun_n)\subset \Coh_n'$).

The relations between these sheaves are as follows:
$${\mathcal S}'_E|_{\E^0_n}\simeq {\mathcal F}_{E,n} \quad \text{ and
} \quad {\mathcal S}'_E|_{\Bun'_n}\simeq \Aut_E',$$ up to
cohomological shifts and Tate's twists. These isomorphisms are
established in \lemref{coincide} and \lemref{comp2&3},
respectively. In particular, $${\mathcal
F}_{E,n}|_{\varrho^{-1}(\Bun_n\cap \C_n)}\simeq
\Aut'_E|_{\varrho^{-1}(\Bun_n\cap \C_n)}.$$ Thus, the primary role of the
first construction is to establish a link between the second and the
third ones.

Our first goal is to prove \corref{clean1} that if $E$ is irreducible,
then ${\mathcal F}_{E,n}$ (obtained as the result of the second
construction) is irreducible and perverse when restricted to each
connected component (up to a cohomological shift and a twist). This
complex is obtained by iterating the operations of Fourier transform,
which are known to preserve irreducibility and perversity, and of
$!$--extensions with respect to the open embeddings $j_k$. Hence we
need to show that the $!$--extensions of the intermediate sheaves
$\F_{E,k}$ coincide with their Goresky-MacPherson extensions, i.e.,
that $\F_{E,k}$ are clean with respect to $j_k$. This is the content
of \thmref{clean}, which is derived in
\secref{cleanness} from the Vanishing \conjref{vanishing conjecture}.

Next, we will prove that ${\mathcal F}_{E,n}$ descends to a perverse
sheaf on the open subset $\C_n\cap\Bun_n\subset \C_n$ under a natural
smooth morphism $\rho^0_n: \E^0_n \to {\mathcal C}_n$ (see
\thmref{descends}). This will be done using the perversity and
irreducibility of ${\mathcal F}_{E,n}$ and some information about the
Euler characteristics of the stalks of the sheaf $\Aut'_E$. In order
to compute these Euler characteristics, we use the third construction
(and its relation to the second construction) in an essential way.

Having obtained the perverse sheaf $\oS_E$ on $\C_n\cap\Bun_n$, we
take its Goresky-MacPherson extension to the union of those connected
components of $\Bun_n$ which have a non-empty intersection with
$\C_n$, i.e., to $\underset{d\geq c_{g,n}}\cup\, \Bun_n^d$.  This
gives us a perverse sheaf $\Aut_E$ on this stack, irreducible on
each connected component.

Then we prove that $\Aut_E$ may
be extended to the entire stack $\Bun_n$ in such a way that it is a
Hecke eigensheaf. We give two independent proofs of this statement,
in \secref{Hecke property} and \secref{hecke laumon},
using the third and the second constructions, respectively.
Finally, in \secref{cuspidality} we prove that the Hecke property of
$\Aut_E$ and \conjref{vanishing conjecture} imply that $\Aut_E$ is
automatically cuspidal.

This summarizes the main steps in our proof of \conjref{geometric
Langlands}. In addition, we will prove the following result.  Let us
denote by $\S_E$ the Goresky-MacPherson extension of $\Aut_E$ from
$\Bun_n$ to $\Coh_n$. In \secref{hecke laumon} we will show that the
sheaves $\varrho_n^*(\S_E)$ and $\S'_E$ on $\Coh'_n$ are isomorphic,
up to a cohomological shift and Tate's twist. The same is true for the
sheaves $\varrho_n^*(\Aut_E)$ and $\Aut'_E$ on $\Bun'_n$.

Note, however, that neither ${\mathcal S}_E'$ nor
$\Aut'_E\simeq\S'_E|_{\Bun'_n}$ is perverse on the entire stack
$\Coh'_n$ and $\Bun'_E$, respectively. This does not contradict the
above assertions: although the morphism $\varrho_n$ is smooth over
$\C_n$, it is not smooth over the entire $\Coh_n$.

\section{The construction via the Whittaker sheaf}

In this section we will present another construction of the sheaf
$\S'_E$ (more precisely, of its restriction to
$\Bun'_n$).\footnote{This construction was independently found by
I.~Mirkovi\'c.} Conceptually, this construction should be viewed as a
geometric counterpart of the construction of automorphic functions for
$GL_n$ from the Whittaker functions due to Piatetskii-Shapiro \cite{PS}
and Shalika \cite{Sha} (see \cite{FGKV} and \secref{Whittaker
function} below for more details).

\ssec{Drinfeld's compactification}  \label{Drinfeld's compact}

We introduce the stack $\ol\Q$, which classifies the data
$(\M,(s_i))$, where $\M$ is a rank $n$ bundle and $s_i$, $i=1,...,n$,
are injective homomorphisms of coherent sheaves
\begin{equation}\label{Drinfeld space}
\begin{split}
\Omega^{n-1} &\xrightarrow {\ s_1 \ }\M
\\
\Omega^{(n-1) + (n-2)} &\xrightarrow  {\ s_2 \ }\Lambda ^2 \M
\\
&\ldots
\\
\Omega^{ \frac {n(n-1)} 2} &\xrightarrow  {\ s_{n-1} \ } \Lambda ^{n-1} \M
\\
\Omega^{ \frac {n(n-1)} 2} &\xrightarrow  {\ s_{n } \ } \Lambda ^n \M\,
\end{split}
\end{equation}
satisfying the requirement that at the generic point of $X$ the
collection $(s_i)$ defines a complete flag of subbundles in $\M$ (see
\cite{FGV}).

In concrete terms, this requirement may be phrased as follows: the
transposed maps $s_i^*: \Lambda^i \M^* \to
\Omega^{-(n-1)-\ldots-(n-i)}$ should satisfy the Pl\"ucker relations
\begin{equation}    \label{plucker}
(s^*_p \otimes s^*_q)(\Phi^k_{p,q}) = 0, \qquad 1\leq k\leq p\leq
q\leq n-1.
\end{equation}
Here $\Phi^k_{p,q}$ is the subsheaf of $\Lambda^p \M \otimes \Lambda^q
\M$ spanned by elements of the form
\begin{multline}    \label{Phi}
(v_1 \wedge \ldots \wedge v_q) \otimes (w_1 \wedge \ldots \wedge w_q)
\\ - \sum_{i_1 < \ldots < i_k} (v_1 \ldots w_1 \ldots w_k \ldots
v_p) \otimes (v_{i_1} \ldots v_{i_k} \wedge w_{k+1} \ldots w_q)
\end{multline}
(i.e., we exchange $k$--tuples of elements of the set $\{ v_i \}$ with
the first $k$ elements of the set $\{ w_j \}$ preserving the order).

To motivate this definition, denote by $\on{Fl}(V)$ the variety of
full flags in an $n$--dimensional vector space $V$ over $\kk$. We
have a natural embedding
$$(s_1,\ldots,s_{n-1}): \on{Fl}(V) \hookrightarrow \prod_{i=1}^{n-1}
{\mathbb P} \Lambda^i V.$$ According to the results of \cite{Tow} and
\cite{Ful}, Sect.~9.1, the ideal of the image of $\on{Fl}(V)$ under
this map is generated by the elements of the form
\eqref{Phi}.\footnote{If $\on{char} \kk = 0$, then the generators with
$k=1$ suffice, but not always so if $\on{char} \kk > 0$, see
\cite{Tow}.}

Thus, if all of the above homomorphisms
$$s_i: \Omega^{(n-1)+\ldots+(n-i)} \to \Lambda^i \M, \qquad
i=1,\ldots,n-1,$$ are maximal embeddings (i.e., bundle maps), then the
data of $(\M,(s_i))$ determine a full flag of subbundles of $\M$. We
denote the open substack of $\ol\Q$ classifying the data $(\M,(s_i))$
satisfying this condition by $\Q$ and the open embedding
$\Q\hookrightarrow \ol\Q$ by $j$.

\medskip

We will denote by $\ol\Q^d$ (resp., $\Q^d$, $\wt\Q^d$) the connected
component of $\ol\Q$ (resp., $\Q$, $\wt\Q$) corresponding to vector
bundles $\M$ of degree $d$ (recall our notational convention from
\secref{notat}). Note that $\dim(\ol\Q^d)=dn+\dim(\ol\Q^0)$.

There is a morphism
\begin{equation}    \label{tau}
\tau: \ol\Q^d \longrightarrow X^{(d)},
\end{equation}
sending $(\M,(s_i))$ to $D$, the divisor of zeros of the last map
$s_n: \Omega^{\frac {n(n-1)} 2} \longrightarrow \Lambda^n M$ in
\eqref{Drinfeld space}. Denote by $\ol\Q^D$ (resp., $\Q^D$) the
preimage of $D$ under $\tau$ in $\ol\Q^d$ (resp., $\Q^D$).

\ssec{Remark.}

The stack $\ol\Q^D$ is the stack $\ol\Bun_N^{\F_T}$ defined in
\cite{FGV}, Sect.~2.2.2, where $\F_T$ is the $T$--bundle on $X$, which
corresponds to the $n$--tuple of line bundles
$(\Omega^{n-1},\dots,\Omega$, $\O(D))$. The stack $\Q^D$
is the stack $\Bun_N^{\F_T}$ from \cite{FGV}, Sect.~2.2.1.

We recall from \cite{FGV,BG1} that the Drinfeld compactification
$\ol{\Bun}_N^{\F_T}$ classifies the data
$(\M,(\kappa^{\check{\la}}))$, where $\check{\la}$ runs over the set
of dominant weights of $GL_n$. Further, $\kappa^{\check{\la}}$ is a
homomorphism of coherent sheaves $\L^{\check{\la}}_{\F_T} \to
\V^{\check{\la}}_{\F_T}$, where $\L^{\check{\la}}_{\F_T}$ is the line
bundle $\F_T \underset{T}\times \check\la$, and
$\V^{\check{\la}}_{\F_T}$ is the vector bundle corresponding to $\M$
and the Weyl representation of $GL_n$ of highest weight $\lambda$. In
addition, the homomorphisms $\kappa^{\check{\la}}$ have to satisfy a
set of Pl\"ucker type relations described in Sect. 2.2.2 of
\cite{FGV}. These relations determine all $\kappa^{\check{\la}}$'s
from $s_i := \kappa^{\check{\omega}_i}, i=1,\ldots,n-1$. Equivalently,
these relations may be described in the form \eqref{plucker}.

\ssec{The Whittaker sheaf}    \label{Aut'}

Observe that the substack $\Q^d$ of $\ol\Q^d$ embeds as an open
substack into $\wt\Q^d$, which classifies those quadruples
$(\M,\beta,(\M^0_i),(s_i))$, for which $\M$ is torsion-free and
$\M^0_i\hookrightarrow \M$ are maximal embeddings (i.e., a bundle maps)
for $i=1,\ldots,n-1$. Recall the morphism $\on{ev}: \wt\Q \to
\GG_a$. Denote its restriction to $\Q^d \subset \wt\Q^d$ also by
$\on{ev}$. First, we define the complex $\Psi^0$ on
$\ol\Q^0$ as
$$\Psi^0 = j_! \on{ev}^*({\mathcal I}_\psi)\otimes
\Ql(\frac{\dim(\Q^0)}{2})[\dim(\Q^0)].$$

We have a natural morphism $q:\ol{\Q} \to \Bun_n$ taking $(\M,(s_i))$ to
$\M$. Recall the stack $\Mod_n^d$ from \secref{averaging} and consider
the Cartesian product
\begin{equation*}
Z^d:=\ol\Q^0\underset{\Bun_n}\times \Mod_n^d.
\end{equation*}
Let $'\hr:Z^d\to \ol\Q^d$ be the morphism that sends
$(\M,(s_i),\M',\beta:\M\to \M')$ to $(\M',(s'_i))$, where $s'_i$ is
the composition
\begin{equation*}
\Omega^{(n-1)+\dots+(n-i)}\overset{s_i}\longrightarrow
\Lambda^i \M \overset{\beta}\longrightarrow \Lambda^i \M'.
\end{equation*}
It is clear that $'\hr$ is a proper morphism of stacks, which makes
the following diagram commutative:
\begin{equation}   \label{diagram for whit}
\begin{CD}
\ol\Q^0 @<{'\hl}<< Z^d @>{'\hr}>> \ol\Q^d \\
@V{q}VV          @V{'q}VV         @V{q}VV \\
\Bun_n @<{\hl}<< \Mod_n^d @>{\hr}>> \Bun_n
\end{CD}
\end{equation}

The {\em Whittaker sheaf} ${\mathcal W}^d_E$ on $\ol\Q^d$ is defined
by the formula\footnote{In the case of $GL_2$ this sheaf was
studied in \cite{Lys:dis}.}
\begin{equation}  \label{def of Whit}
{\mathcal W}^d_E:={}'\hr_!({}'\hl{}^*(\Psi^0)\otimes
(\pi\circ{}'q)^*(\L^d_E))\otimes \Ql(\frac{d\cdot n}{2})[d\cdot n].
\end{equation}
In other words, ${\mathcal W}^0_E\simeq \Psi^0$, and in general
${\mathcal W}^d_E$ is obtained from ${\mathcal W}^0_E$ via a
$\ol\Q$--version of the averaging functor $\Hav^d_{n,E}$. Namely, define
the functor $'\Hav^d_{n,E}: \on{D}(\ol\Q^0) \to \on{D}(\ol\Q^d)$ by
the formula
\begin{equation}   \label{Hecke upstairs}
\K \mapsto {}'\hr_!({}'\hl{}^*(\K)\otimes
(\pi\circ{}'q)^*(\L^d_E))\otimes \Ql(\frac{d\cdot n}{2})[d\cdot n].
\end{equation}
Then ${\mathcal W}^d_E={}'\Hav^d_{n,E}(\Psi^0)$.

Let $\W_E$ be the complex on $\ol\Q$, whose restriction to $\ol\Q^d$
equals $\W^d_E$. Denote by $\nu: \ol\Q \to\Bun'_n$ the morphism which
sends $(\M,(s_i))$ to $(\M,s_1: \Omega^{n-1}\to\M)$. Set
\begin{equation*}
\Aut'_E:=\nu_!({\mathcal W}_E).
\end{equation*}

\begin{lem} \label{comp2&3}
The complex $\Aut'_E$ is canonically isomorphic to the restriction
of $\S'_E$ from $\Coh'_n$ to $\Bun_n'$.
\end{lem}

\begin{proof}
Let $\wt\Q^d_{tf}$ be the locus in $\wt\Q^d$ corresponding to those data
$(\M,\beta,(\M^0_i),(\wt{s}_i))$, for which $\M$ is
torsion-free. Observe that $\ovc{Z}{}^d:=({}'\hl{})^{-1}(\Q^0)\subset
Z^d$ is canonically identified with $\wt\Q^d_{tf}$. Since $\Psi^0$ is
extended by zero from $\Q^0$ to $\ol\Q^0$, the proposition follows
from the commutativity of the diagram
\begin{equation*}
\begin{CD}
\ovc{Z}{}^d @>\sim>> \wt\Q^d_{tf} \\
@V{'\hr}VV      @V{\wt\nu}VV  \\
\ol\Q^d @>\nu>> \Bun'_n
\end{CD}
\end{equation*}
which is verified directly from the definitions.
\end{proof}

\ssec{The structure of ${\mathcal W}^d_E$}

In the rest of this section we describe the structure of the Whittaker
sheaf using the results of our previous work \cite{FGV}. Strictly
speaking, these results are not necessary in our proof of the
geometric Langlands \conjref{geometric Langlands}. They are used only
in the proof of the Hecke eigensheaf property presented in
Sects. \ref{proof 1}--\ref{proof last}, for which we give an alternative
proof in \secref{hecke laumon}, which does not use the Whittaker
sheaf.

To state our results, we recall first that the substack $\Q^d$ of
$\ol\Q^d$ embeds as an open substack into $\wt\Q^d$ and the
composition $\Q^d\to \ol\Q^d\overset{\alpha}\to \Coh^d_0$ takes values
in $\Coh^{\on{r},d}$ (see \secref{deligne van thm} for this
notation). Hence,
\begin{equation}
{\mathcal W}^d_E|_{\Q^d}\simeq \wt{\mathcal
W}^d_E|_{\Q^d}\simeq \tau^*(E^{(d)})|_{\Q^d} \otimes
\on{ev}^*({\mathcal I}_\psi)\otimes \Ql(\frac{\dim(\Q^d)}{2})[\dim(\Q^d)]
\end{equation}
where $\tau$ is the morphism in \eqref{tau}.

We will prove the following statement:

\begin{prop}  \label{Whit is perverse}
The complex ${\mathcal W}^d_E$ on $\ol\Q^d$ is an irreducible perverse
sheaf  which is the Goresky-MacPherson extension of its restriction to
$\Q^d$.
\end{prop}

The first step in the proof of \propref{Whit is perverse} is provided
by the following result of \cite{FGV}:

\begin{lem} \label{Psi is clean}
The canonical map $\Psi^0 \to j_{!*} j^*(\Psi^0)$ is an isomorphism.
\end{lem}

In other words, the perverse sheaf $\Psi^0$ on
$\Q^0$ is clean with respect to $j:\Q^0\hookrightarrow\ol\Q^0$.
Therefore, since ${\mathcal W}^d_E\simeq {}'\Hav^d_{n,E}(j_{!*} j^*(\Psi^0))$
(cf. formula \eqref{Hecke upstairs}), we obtain that
$${\mathbb D}({\mathcal W}^d_E) \simeq {\mathcal W}^d_{E^*},$$
where ${\mathbb D}$ is the Verdier duality functor, and $E^*$ is
the dual local system to $E$.

Thus, ${\mathcal W}^d_E$ is Verdier self-dual up to replacing
$E$ by $E^*$ and to prove \propref{Whit is perverse} it
suffices to introduce a stratification of $\ol\Q^d$ and to show that
the $*$--restriction of ${\mathcal W}^d_E$ to each stratum in $\ol\Q^d
- \Q^d$ appears in strictly negative cohomological degrees (with
respect to the perverse $t$--structure). This stratification $\{
\Q^{\ol\mu} \}$ is introduced in \secref{strat of Q}. The description
of the restriction of ${\mathcal W}^d_E$ to each stratum given in
\propref{cass-shal} will imply that it appears in strictly negative
cohomological degrees (except for the open stratum).

\ssec{Substacks of $\ol\Q$ defined by orders of zeros of sections}
\label{explicit description}

The Langlands dual group to $GL_n$ is $GL_n(\Ql)$. In what follows we
represent each weight of $GL_n(\Ql)$ as a string of integers
$(d_1,...,d_n)$, so that dominant weights satisfy $d_1\geq d_2\geq
\ldots\geq d_n$. We denote the set of dominant weights by $P^+_n$. The
irreducible finite-dimensional representation of $GL_n(\Ql)$ with
highest weight $\mu \in P^+_n$ will be denoted by $V^\mu$. We denote
by $w_0$ is the longest element of the Weyl group of $GL_n$, which
acts on the weights by the formula $w_0 \cdot (d_1,...,d_n) =
(d_n,\ldots,d_1)$. For an anti-dominant weight $\mu$, we denote by
$V_\mu$ the irreducible finite-dimensional representation of $GL_n$
with lowest weight $\mu$, i.e., $V_\mu \simeq V^{-w_0(\mu)}$.

Let $\ol\mu=\{\mu^1,...,\mu^m\}$ be
a collection of weights of $GL_n(\Ql)$, where some of the $\mu^j$'s
may coincide. We will denote by $X^{\ol\mu}$ the corresponding partially
symmetrized power of $X$ with the all the diagonals removed. In other
words, if $m=m_1+\ldots+m_s$ is such that a given weight $\mu^r$ appears
in the collection exactly $m_r$ times, then
$$X^{\ol\mu}= X^{(m_1)}\times\ldots\times X^{(m_s)}-\Delta.$$

We will think of a point $\ol x$ of $X^{\ol\mu}$ as of a collection
of pairwise distinct points $x^j$, $j=1,...,m$, to each of which
there is an assigned weight $\mu^j=(d_1^j,\ldots,d_n^j)$.

\medskip

We associate to $\ol\mu$ a stack $\ol\Q^{\ol\mu}$, which classifies
the data $(\M,(s_i),\ol x)$, where $\M$ is a vector bundle of rank
$n$, $\ol x$ is a point of $X^{\ol\mu}$ represented by a collection of
distinct points $x^j\in X$, and
\begin{equation*}
s_i: \Omega^{(n-1)+...+(n-i)}\left(\underset{j}\Sigma\,
\,(d_1^j+...+d^j_i)\cdot x^j \right)\hookrightarrow
\Lambda^i \M
\end{equation*}
are injective homomorphisms of coherent sheaves satisfying the
Pl\"ucker relations from \cite{FGV,BG1}.

The locus where all the maps $s_i$ are maximal embeddings (i.e., are
bundle maps) is an open substack $\Q^{\ol\mu}$ of $\ol\Q^{\ol\mu}$. In
other words, $\Q^{\ol\mu}$ classifies the data $(\M,(\M_i),(\wt
s_i),\ol x)$, where $\M$ is a rank $n$ bundle, $0=\M_0\subset
\M_1\subset...\subset \M_n=\M$ is a full flag of {\it subbundles} of
$\M$ and $\wt s_i, i=1,\ldots,n$, is an isomorphism
$$\Lambda^i \M_i \simeq \Omega^{(n-1)+...+(n-i)}\left(\underset{j}\Sigma \,
(d_1^j+...+d^j_i)\cdot x^j
\right).$$

For a fixed point $\ol x\in X^{\ol\mu}$, we will denote by $\ol x$, by
$\ol\Q^{\ol\mu,\ol x}$ and $\Q^{\ol\mu, \ol x}$ the corresponding
closed substacks of $\ol\Q^{\ol\mu}$ and $\Q^{\ol\mu}$, respectively.

For reader's convenience, we identify the above stacks with those
studied in \cite{FGV}: $\Q^{\ol\mu,\ol{x}}\simeq \Bun_N^{\F_T}$ and
$\ol\Q^{\ol\mu,\ol x} \simeq \overline\Bun_N^{\F_T}$, where $\F_T$ is
the $T$--bundle on $X$, which corresponds to the $n$--tuple of line
bundles
\begin{equation*}
\left(\Omega^{n-1}\left(\underset{j}\Sigma\, d_1^j\cdot x^j
\right),\dots,\Omega\left(\underset{j}\Sigma\, d_{n-1}^j\cdot
x^j\right),\O\left(\underset{j}\Sigma\, d_n^j\cdot x^j \right) \right).
\end{equation*}

\subsection{Stratification of $\ol\Q^d$}    \label{strat of Q}

If the collection $\ol\mu$ satisfies the conditions
\begin{equation}    \label{cond for d-pos}
d_1^j + \ldots + d_i^j \geq 0, \quad \forall i,j, \qquad \sum_{i,j}
d_i^j = d,
\end{equation}
then we have a natural closed embedding
$\ol\Q^{\ol\mu}\hookrightarrow
\ol\Q^d$, and so $\Q^{\ol\mu}$ is a locally closed substack of
$\ol\Q^d$. In particular, for a divisor $D = \sum_{j=1}^m d^j\cdot
x^j \in X^{(d)}$ we have isomorphisms $\ol\Q^D \simeq \ol\Q^{\ol\mu,\ol x}$
and $\Q^D \simeq \Q^{\ol\mu,\ol x}$, where $\mu^j=(0,...,d^j)$.

The following statement follows from \cite{FGV},  Corollary 2.2.9.

\begin{lem}    \label{strata}
The locally closed substacks $\Q^{\ol\mu}$ with $\ol\mu$ satisfying
the condition \eqref{cond for d-pos} form a stratification of
$\ol\Q^d$.
\end{lem}

\subsection{Restrictions of $\W_E$ to the strata}

The collection $\ol\mu$ is called {\em anti-dominant} if all the weights
$\mu^j$ are anti-dominant. For an anti-dominant $\ol\mu$, we have a
map $\on{ev}^{\ol\mu}: \Q^{\ol\mu}\to \GG_a$ defined as in
\cite{FGV}. Namely,
\begin{equation*}
\begin{split}
\on{ev}: {\mathcal Q}^{\ol\mu} \arr \prod_{i=1}^{n-1}
\Ext^1(\M_i/\M_{i+1},\M_{i-1}/\M_i)) \arr \prod_{i=1}^{n-1}
\HH^1\left(X,\Omega\left(\underset{j}\Sigma\, (d_i^j-d_{i+1}^j)\cdot
x^j\right)\right) \\ \arr \HH^1(X,\Omega)^{\oplus n-1} \simeq \GG_a^{n-1}
\overset{\on{sum}}\arr \GG_a
\end{split}
\end{equation*}
(compare with formula \eqref{ev}). We then set
\begin{equation}
\ovc\Psi{}^{\ol\mu}:=\on{ev}^{\ol\mu}{}^*({\mathcal I}_\psi) \otimes
\Ql(\frac{\dim(\Q^{\ol\mu})}{2})[\dim(\Q^{\ol\mu})].
\end{equation}
Denote by $j^{\ol\mu}$ the embedding $\Q^{\ol\mu} \rightarrow
\ol\Q^{\ol\mu}$. According to Theorem 2 of \cite{FGV}, the sheaf
$\ovc\Psi{}^{\ol\mu}$ is clean, i.e., $\Psi{}^{\ol\mu}:=
j^{\ol\mu}_!(\ovc\Psi{}^{\ol\mu})$ is an irreducible perverse sheaf
isomorphic to $j^{\ol\mu}_{!*}(\ovc\Psi{}^{\ol\mu})$. In a similar
way we define the perverse sheaves $\ovc\Psi{}^{\ol\mu,\ol x}$ and
$\Psi{}^{\ol\mu,\ol x}$ on $\Q^{\ol\mu,\ol x}$ and $\ol\Q^{\ol\mu,\ol
x}$, respectively.

Next, to a local system $E$ on $X$ and an anti-dominant collection
$\ol\mu$ we associate a local system $E_{\ol\mu}$ on $\Q^{\ol\mu}$ as
follows. Recall that we denote by $V_\mu$ the irreducible
representation of $GL_n$ with {\em lowest} weight $\mu \in
-P^+_n$. Let $E_{\mu}$ be the local system on $X$ associated to $E$
and $V_\mu$.

For $\ol\mu$ corresponding to the partition $m=m_1+\ldots+m_s$
consider the sheaf $(E_{\mu_1})^{(m_1)} \boxtimes \ldots \boxtimes
(E_{\mu_s})^{(m_s)}$ on $X^{(m_1)} \times \ldots \times
X^{(m_s)}$. Denote by $E_{\ol\mu}$ its restriction to the compliment
of all diagonals, i.e. to $X^{\ol\mu}$. Let us denote by $\tau_{\ol\mu}$
the natural morphism from $\Q^{\ol\mu}$ to $X^{\ol\mu}$ and set
$'E_{\ol\mu}:=\tau_{\ol\mu}^*(E_{\ol\mu})$.

Thus, $'E_{\ol\mu}$ depends only on the positions of the points
$x^1,...,x^m$, and its fiber at $(\M,(\wt s_i),(x^1,...,x^m))$ is
$\otimes_{j=1}^m E_{\mu^j,x^j}$, where $E_{\mu^j,x^j}$ denotes the
fiber of $E_{\mu^j}$ at $x^j$.

\begin{prop}  \label{cass-shal}
The $*$--restriction of ${\mathcal W}^d_E$ to $\Q^{\ol\mu}\subset
\ol\Q^d$ is zero unless all weights $\mu^j$ are anti-dominant. When
they are, this restriction is canonically identified with
$\ovc\Psi{}^{\ol\mu}\otimes {}'E_{\ol\mu} \otimes \Ql(\frac{d-m}{2})[d-m]$.
\end{prop}

As was explained earlier, this proposition implies \propref{Whit is
perverse}.  Indeed, we had to show that $*$--restriction of ${\mathcal
W}^d_E$ to each stratum in $\Q^{\ol\mu}\subset \ol\Q^d - \Q^d$ appears
in negative cohomological degrees.  According to \propref{cass-shal},
this restriction lives in the cohomological degree $-(d-m)$. However,
since each $\mu^j=(d^j_1,...,d^j_n)$ satisfies $d^j_1\leq...\leq
d^j_n$ and $d^j_1+...+d^j_n\geq 0$, we obtain that $d-m\geq 0$ and the
equality takes place only when every $\mu^j$ is of the form
$(0,...,0,1)$. However, the stratum corresponding to this $\ol\mu$ is
contained in $\Q$.

The proof of \propref{cass-shal} is given in the Appendix. A similar
calculation has also been performed in \cite{Lys}.

\ssec{The Whittaker function}    \label{Whittaker function}

In this subsection we will assume that $\kk$ is the finite field
$\Fq$. We will show that the function associated with the Whittaker
sheaf may be identified with the restriction of the Whittaker
function.

First we briefly recall the definition of the Whittaker function (see
\cite{FGKV}, Sect. 2, for more details). Consider group $GL_n(\AA)$
over the ring of ad\`eles of $F=\Fq(X)$, and let $N(\AA)$ be its upper
unipotent subgroup. Denote by $u_{i,i+1}$ the $i$-th component of the
image of $u \in N(\AA)$ in $N(\AA)/[N(\AA),N(\AA)]$ corresponding to
the $(i,i+1)$ entry of $u$. Recall that we have fixed a non-trivial
additive character $\psi: \Fq \arr \Ql^\times$. We define the
character $\Psi$ of $N(\AA)$ by the formula\footnote{For this formula
to be well-defined, we should consider a twisted version $GL_n^J(\AA)$
of $GL_n(\AA)$ introduced in \cite{FGKV}, Sect. 2; then $u_{i,i+1}$
is naturally a differential.}
$$\Psi\left( (u_x)_{x\in |X|} \right) = \prod_{i=1}^{n-1} \prod_{x \in
|X|} \psi(\on{Tr}_{\kk_x/\Fq}(\on{Res}_x u_{x,i,i+1})).$$ It follows
from the residue theorem that $\Psi(u)=1$ if $u \in N(F)$.

Now let $E$ be a rank $n$ local system on $X$. Then there exists a
unique (up to a non-zero scalar multiple) function $W_{E}$ on
$GL_n(\AA)$, which is right $GL_n(\OO)$--invariant, left
$(N(\AA),\Psi)$--equivariant, and is a Hecke eigenfunction associated
to $E$. This function is called the Whittaker function corresponding
to $E$. Casselman--Shalika \cite{CS} and Shintani \cite{Shi} have
given an explicit formula for $W_E$ (see, e.g., Theorem 2.1 of
\cite{FGKV}).  The left $(N(\AA),\Psi)$--equivariance of $W_E$ implies
that it is left $N(F)$--invariant, where $F=\Fq(X)$. Therefore we
obtain a function on the double quotient $Q = N(F)\backslash
GL_n(\AA)/GL_n(\OO)$. We denote this function also by $W_E$.

In the same way as in the the proof of Lemma 2.1 from \cite{FGKV}, we
identify the set of $\Fq$--points of $\Q^{\ol\mu,\ol x}$, where $\mu^j
= (d^j_1,\ldots,d^j_n), j=1,\ldots,m$, with the projection onto $Q$ of
the subset
$$
N(\AA) \cdot \left( \on{diag}(\pi_x^{d_n(x)},\ldots,\pi_{x}^{d_1(x)})
\right)_{x \in |X|} \cdot GL_n(\OO) \subset GL_n(\AA),
$$
where $d_i(x) = d^j_i$, if $x=x_j$, and $d_i(x) = 0$, otherwise. Thus,
we may embed the set $\ol\Q^d(\Fq)$ of isomorphism classes of
$\Fq$--points of $\ol\Q^d$ into $Q$ for all $d\geq 0$. Comparing
\propref{cass-shal} with the Casselman-Shalika-Shintani formula, we
obtain:

\begin{prop}
The function $\text{\tt f}_q(\W^d_E)$
on $\ol\Q^d(\Fq)$ corresponding to the sheaf $\W^d_E$
equals the restriction of the Whittaker function $W_E$ to
$\ol\Q^d(\Fq) \subset Q$.
\end{prop}

Furthermore, the geometric construction of the sheaf $\Aut'_E$
described in this section translates at the level of functions into
the construction due to Shalika \cite{Sha} and Piatetskii-Shapiro
\cite{PS} (see Sect. 2 of \cite{FGKV} for a review of this
construction).

\section{Cleanness in Laumon's construction} \label{cleanness}

Let $E$ be an irreducible rank $n$ local system for which
\conjref{vanishing conjecture} holds. In this section we derive
\thmref{clean} for $E$, i.e., we  prove that the complex
$\F_{E,k}$ on $\E^0_k$ is clean
with respect to $j_k: \E^0_k \hookrightarrow \E_k$.

\medskip

We will begin by stating a well-known lemma, which will serve as one
of the key ingredients in the proof. Consider the following situation:
let $\E$ be a vector bundle over a scheme (or stack) $Y$.  Let us
denote by $\rho:\E\to Y$ the projection, by $i:Y\to \E$ the $0$
section, and by $\E^0\overset{j}\hookrightarrow \E$ the complement of
the zero section. Assume that $\K$ is a complex on $\E^0$, equivariant
with respect to the $\GG_m$-action.

\begin{lem} \label{cleanness trick}
The complex $\K$ is clean with respect to $j$ if and only if
$(\rho\circ j)_!(\K)=0$.
\end{lem}

To prove the lemma, it suffices to note that cleanness of $\K$ is
equivalent to  the
statement that $i^! j_!(\K)=0$. But for any $\GG_m$--equivariant
complex $\K'$ on
$\E$, we have: $i^!(\K')\simeq \rho_!(\K')$.

Our proof of \thmref{clean} will proceed by induction on the length
of torsion in
$\C_k$.  Let us first consider the case where there is no torsion at
all, i.e., we will
show that $j_k{}_!(\F_{E,k})\to j_k{}_*(\F_{E,k})$ is an isomorphism
on the open set $\rho_k^{-1}(\C_k\cap \Bun_k) \subset \C_k$.

Recall that $\rho_k$ (resp., $\rho_k^0$) denotes the projection
$\E_k\to \C_k$ (resp., $\E^0_k\to \C_k$).  On $\rho_k^{-1}(\C_k\cap
\Bun_k)$, $\E^0_k$ is the complement of the zero section in $\E_k$ and
$\F_{E,k}$ is $\GG_m$--equivariant. Thus, by \lemref{cleanness trick},
we are reduced to showing that $\rho^0_k{}_!(\F_{E,k})|_{\C_k\cap
\Bun_k}=0$. By \lemref{comp2&3}, with $n$ replaced by $k$, we obtain
that up to a cohomological shift and Tate's twist
\begin{equation*}
\rho^0_k{}_!(\F_{E,k})|_{\C_k\cap \Bun^d_k}\simeq
\Hav^d_{k,E}(q_!(\Psi^0))|_{\C_k\cap \Bun^d_k}.
\end{equation*}
The definition of $\C_k$ in \secref{fund diag} implies that if
$\C_k\cap \Bun_k\neq\emptyset$ then $d>nk(2g-2)$. Thus, the Vanishing
\conjref{vanishing conjecture} implies
$\rho^0_k{}_!(\F_{E,k})|_{\C_k\cap \Bun_k}=0$.

\ssec{Induction on the length of torsion}   \label{induction on torsion}

To set up the induction, we fix some notation. For an integer $\ell$,
let us write $\C_{k,\leq \ell}$ for the open substack of $\C_k$
consisting of coherent sheaves whose torsion is of length $\leq
\ell$. Set $\C_{k,<\ell}=\C_{k,\leq \ell-1}$ to be the open substack
in $\C_{k,\leq \ell}$ that corresponds to the locus where the torsion
is of length $<\ell$. Finally, let $\C_{k,\ell}$ be the closed
substack of $\C_{k,\leq\ell}$ corresponding to coherent sheaves whose
torsion is precisely of length $\ell$.

Set $\E_{k,\leq \ell}=\rho_k^{-1}(\C_{k,\leq \ell})$,
$\E_{k,<\ell}=\rho_k^{-1}(\C_{k,<\ell})$, and
$\E_{k,\ell}=\rho_k^{-1}(\C_{k,\ell})$.  Furthermore, let us write
\begin{align*}
&\E^0_{k,\leq \ell}=\E_k^0\cap  \E_{k,\leq \ell},\,\,\,
\E^0_{k,<\ell}=\E_k^0\cap  \E_{k,<\ell}, \\
&\E_k^t = \E_k-\E_k^0,\,\,\, \E_{k,\ell}^t=\E_k^t \cap \E_{k,\ell},
\text{ etc.}
\end{align*}

We assume, by induction, that $\F_{E,k}$ is clean with respect to the
inclusion $\E^0_{k,<\ell}\hookrightarrow \E_{k,<\ell}$. To show
cleanness of $\F_{E,k}$ with respect to the inclusion
$\E^0_{k,\leq\ell}\hookrightarrow \E_{k,\leq\ell}$, it suffices to
prove cleanness of $\F_{E,k}$ with respect to
$(\E_{k,\leq\ell}-\E_{k,\ell}^t)\hookrightarrow \E_{k,\leq\ell}$. We
would like to argue in the same manner as we did above in the case
when $\ell$ was zero and using \lemref{cleanness trick}.
Unfortunately, we cannot apply this lemma directly, because
$\E_{k,\leq\ell}$ is not a vector bundle over
$\E_{k,\ell}^t$. However, it will become a vector bundle after a
smooth base change.

Consider the stack $\Coh_{k, \leq \ell}$ which classifies coherent
sheaves of generic rank $k$ with torsion of length $\leq \ell$, and
the stack $\widetilde\Coh_{k, \leq \ell}$, which classifies the
following data: $\M_0\in \Bun_k$, $\T\in\Coh_0^\ell$, and a short
exact sequence $$0\to \M_0 \to \M\to\T\to 0.$$ There is a canonical
morphism $r:\widetilde\Coh_{k, \leq \ell} \to \Coh_{k, \leq \ell}$
which associates to a triple as above the coherent sheaf
$\M\in\Coh_k$.

\begin{lem}    \label{r smooth}
The morphism $r$ is smooth.
\end{lem}

\begin{proof}

First, the stack $\Coh_k$ is known to be smooth. One proves this
simultaneously with the fact that $\Coh_k$ is indeed an algebraic stack
in the smooth topology by covering it by a Hilbert scheme.

Therefore both stacks $\Coh_{k, \leq \ell}$ and $\widetilde\Coh_{k,
\leq \ell}$ are smooth, the former being an open substack in $\Coh_k$,
and the latter being a vector bundle over $\Bun_k$.

Hence in order to show that $r$ is smooth, it suffices
to show that the fiber of $r$ is smooth over
any field-valued point $\M\in \Coh_{k, \leq \ell}$. By definition,
the tangent space
to the fiber or $r$ at the point $0\to \M_0
\to \M\to\T\to 0$  is $\Hom(\M_0,\T)$. The dimension of $\Hom(\M_0,\T)$ is
$k\cdot \ell$ because $\M_0$ is torsion--free.
As $r$ is separable and the dimensions of the tangent spaces to the
fibers are constant, we conclude that $r$ is smooth.

\end{proof}

Since $r$ is smooth by \lemref{r smooth}, it is sufficient to prove
cleanness after this base change $r:\widetilde\Coh_{k, \leq \ell} \to
\Coh_{k, \leq \ell}$. Consider the stack $\widetilde\Coh_{k, \leq
\ell} \underset{\Coh_{k, \leq \ell}}\times \E_{k,\leq\ell}$.
It classifies the following data:
\begin{equation}  \label{the space E}
\begin{split}
&\M\in \C_k, \ \ \T\in \Coh_0^\ell, \ \ \M_{tf}\in\Bun_k \\
&0 \to \M_{tf} \to \M \to \T \to 0, \ \ \ \Omega^{k-1} \overset{s}\to \M.
\end{split}
\end{equation}
Its substack $$\tilde\E^t_{k,\ell} = \widetilde\Coh_{k, \leq \ell}
\underset{\Coh_{k, \leq \ell}}\times \E^t_{k,\ell}$$ consists of the
data \eqref{the space E} such that the extension is split and the
image of $s$ belongs to the torsion part of $\M$. In other words,
$\tilde\E^t_{k,\ell}$ classifies the data
\begin{equation} \label {the space E torsion}
\M_{tf}\in \Bun_k, \ \ \T\in \Coh_0^\ell  \text{  (such that }
\M_{tf}\oplus\T\in \C_k),\ \ \ \ \Omega^{k-1}
\xrightarrow{\kappa} \T.
\end{equation}

Denote by $\tilde\E_{k,\leq\ell}$ the open substack of $\widetilde\Coh_{k,
\leq \ell} \underset{\Coh_{k, \leq \ell}}\times \E_{k,\leq\ell}$
defined by the condition
$$\Hom(\M_{tf},\L^{\on{est}})=0\,;$$
here $\L^{\on{est}}$ is the line bundle of  \secref{seccon}.
The above condition guarantees that $\M_{tf}\oplus \T\in \C_k$. Set
$$\tilde\E_{k,\leq\ell}^0=\tilde\E_{k,\leq\ell}\cap{} \left(
\widetilde\Coh_{k, \leq \ell} \underset{\Coh_{k, \leq \ell}}\times
\E^0_{k,\leq\ell} \right).$$ Obviously, $\tilde\E^t_{k,\ell}$ is
contained in $\tilde\E_{k,\leq\ell}$.

There is a natural morphism $\rho_{k,\ell}: \tilde\E_{k,\leq\ell} \to
\tilde\E^t_{k,\ell}$ which maps the quadruple in the definition of
$\tilde\E_{k,\leq\ell}$ to the data $(\M_{tf},\T,\kappa)$, where $\kappa$
is the composition $\Omega^{k-1} \overset{s}\to \M \to \T$.
In this way $\tilde\E_{k,\leq\ell}$ becomes a vector bundle over
$\tilde\E^t_{k,\ell}$; the fiber over $(\M_{tf},\T,\kappa)$ can be canonically
identified with the vector space
\begin{equation*}
\on{Ext}^1(\on{Cone}(\Omega^{k-1} \xrightarrow{\kappa} \T),\M_{tf}).
\end{equation*}
Let $\rho^0_{k,\ell}$ denote the composition $\tilde\E_{k,\leq
\ell}^0\hookrightarrow \tilde\E_{k,\leq
\ell}\overset{\rho_{k,\ell}}\longrightarrow \tilde\E^t_{k,\ell}$.

\medskip

Let $\wt\F_{E,k}$ denote the pull-back of $\F_{E,k}$ to
$\tilde\E_{k,\leq\ell}^0$.  One readily verifies that $\wt\F_{E,k}$ is
equivariant with respect to the $\GG_m$-action along the fibers of
$\rho_{k,\ell}:\tilde\E_{k,\leq\ell}\to \tilde\E^t_{k,\ell}$.

The assertion of the theorem reduces to the fact that $\wt\F_{E,k}$ is
clean with respect to the open embedding
$\tilde\E_{k,\leq\ell}^0\hookrightarrow \tilde\E_{k,\leq\ell}$ and we already
know this assertion on the open substack $\tilde\E_{k,\leq
\ell}-\tilde\E^t_{k,\ell}$.

By applying \lemref{cleanness trick} to $\wt\F_{E,k}$ extended by $0$
from $\tilde\E_{k,\leq\ell}^0$ to $\tilde\E_{k,\leq
\ell}-\tilde\E^t_{k,\ell}$, we reduce our assertion to showing
that
\begin{equation}   \label{rho vanishing}
\rho^0_{k,\ell}{}_!(\wt\F_{E,k})=0.
\end{equation}

\ssec{Proof of formula \eqref{rho vanishing}}   \label{section rho vanishing}

Recall the stack $\wt\Q$ of \secref{the first} and let us denote by
$\wt\Q_k$ (resp., $\Q_k$, $\Q^0_k$, $\ol\Q_k$) its version with $n$
replaced by $k$. In particular, $\Q^0_k$ classifies points of the form
$(\M^0,\beta,(\M^0_i),(\wt{s}_i))$, where the map $\beta: \M^0 \to \M^0
=\M^0_n$ is the identity. We will denote such a point simply by
$(\M^0,(\wt{s}_i))$. Denote by $\wt\W_{E,k}$ the perverse sheaf on
$\wt\Q_k$ defined as in \secref{the first}.

Consider the Cartesian product $\wt\Q_k\underset{\Coh'_k}\times
\tilde\E_{k,\leq\ell}^0$. This is the stack that classifies the data of
\begin{equation}
\begin{split}
&\M_{tf}\in \Bun_k \text{ with }\Hom(\M_{tf},\L^{\on{est}})=0,
\ \ \T\in \Coh_0^\ell, \\
&0 \to \M_{tf} \to \M \to \T \to 0, \ \ \ \M^0\hookrightarrow \M, \\
& (\M^0,\wt{s}_1,...,\wt{s}_k)\in \Q_k^0.
\end{split}
\end{equation}

By \lemref{coincide}, $\wt\F_{E,k}$ is the direct image (with compact
supports) under the
map
$\wt\Q_k\underset{\Coh'_k}\times \tilde\E_{k,\leq\ell}^0\to
\tilde\E_{k,\leq\ell}^0$ of the pull-back of $\wt\W_{E,k}$ from $\wt\Q_k$ to
$\wt\Q_k\underset{\Coh'_k}\times \tilde\E_{k,\leq\ell}^0$. Hence, in
order to prove
\eqref{rho vanishing}, it suffices to show that the compactly
supported cohomology of the
fiber of $\wt\Q_k\underset{\Coh'_k}\times \tilde\E_{k,\leq\ell}^0 \to
\tilde\E_{k,\leq\ell}^0
\to \tilde\E^t_{k,\ell}$ with coefficients in the pullback of
$\wt\W_{E,k}$ vanishes. To this
end, let us fix a point
$(\M_{tf}\in\Bun_k,\T\in\Coh_0^\ell,\kappa:\Omega^{k-1}\to \T)$ in
$\tilde\E^t_{k,\ell}$ and analyze the fiber  over this point.
Let us write $Y$ for the closed
substack of the the fiber which lies over a
fixed  point $(\M^0,\wt{s}_1,...,\wt{s}_k)\in \Q_k^0$ and where the
composition $\phi:\M^0\to \M\to \T$ is also fixed.

{}From the discussion above we conclude that \eqref{rho vanishing}
follows if we show that
  for all
$$(\M_{tf}\in\Bun_k,\ \ \T\in\Coh_0^\ell, \ \
(\M^0,\,\wt{s}_1,...,\wt{s}_k)\in \Q^0_k,\ \ \phi:\M^0\to \T)$$
as above, we have
\begin{equation}    \label{rho vanishing1}
\HH^{\bullet}_c(Y,\wt\W_{E,k}|_Y)=0.
\end{equation}

To prove this, we will first reduce to the case when $\phi:\M^0\to \T$ is
surjective. Let us denote by $\T'$ the image of $\phi$ and by $\T''$
the cokernel of $\kappa$; write $\ell'$ (resp., $\ell''$) for the
length of $\T'$ (resp., $\T''$). Let $Y'$ be the scheme defined in the
same way as $Y$, but for
$$(\M_{tf}, \ \ \T', \ \ (\M^0,\wt{s}_1,...,\wt{s}_k)\in \Q^0_k, \ \
\phi':\M^0\to \T').$$

We have a natural map $\on{v}:Y\to Y'$,
which associates to a point
\begin{equation}\label{Y}
(0 \to \M_{tf} \to \M \to \T \to 0,\,\,
\M^0\hookrightarrow\M)
\end{equation}
the point
\begin{equation}\label{Y'}
(0 \to \M_{tf} \to \M' \to \T' \to 0,\,\,
\phi':\M^0\hookrightarrow\M'),
\end{equation}
where $\M'$ is the preimage of
$\T'$ under $\M\to\T$.

\begin{lem}
The complexes
$\on{v}_!(\wt\W_{E,k}|_Y)$ and $\wt\W_{E,k}|_{Y'}\otimes
(\L_E^{\ell''})_{\T''}$ are isomorphic up to a cohomological shift
and Tate's twist; here $(\L_E^{\ell''})_{\T''}$ is the stalk
of Laumon's sheaf at $\T''\in\Coh_0^{\ell''}$.
\end{lem}

\begin{proof}

Let us recall the following basic property of Laumon's sheaf $\L_E$,
cf. \cite{La2}:

Consider the stack $\on{Fl}^{d',d''}_0$ that classifies
short exact sequences
$$0\to \wt\T'\to \wt\T\to \wt\T''\to 0,\ \ \wt\T'\in \Coh_0^{d'},\ \
\wt\T''\in\Coh_0^{d''}.$$
Let $\mathfrak p$ denote the natural projection
$\on{Fl}^{\ell',\ell''}_0\to \Coh_0^d$
(here $d=d'+d''$), that associates to a short exact sequence as above
its middle term,
and let $\mathfrak q:\on{Fl}^{d',d''}_0\to \Coh_0^{d'}\times \Coh_0^{d''}$
denote the other natural projection.

In \cite{La2} Laumon proved that
\begin{equation} \label{simple Laumon}
{\mathfrak q}_!\circ {\mathfrak p}^*(\L^d_E)\simeq \L_E^{d'}\boxtimes
\L_E^{d''}.
\end{equation}

We have a natural map $Y'\to \Coh_0^{d'}\times \Coh_0^{d''}$ that
sends the data of
$$(0\to \M_{tf}\to \M'\to \T'\to 0, \ \ \M^0\hookrightarrow \M')$$
to
$(\M'/\M^0,\T'')$ with $d'=\deg(\M')-\deg(\M^0),d''=\ell''$ and
a map $Y\to \on{Fl}^{d',d''}_0$ that sends
$$(0\to \M_{tf}\to \M\to \T\to 0, \ \ \M^0\hookrightarrow \M)$$
to
$$0\to \M'/\M^0\to \M/\M^0 \to \T''\to 0.$$
Note that since $\M^0$ is fixed, $\wt\W_{E,k}|_Y$ is isomorphic to the
pull-back of ${\mathfrak p}^*(\L^d_E)$ under this map. Similarly,
$\wt\W_{E,k}|_{Y'}\otimes (\L_E^{\ell''})_{\T''}$ is isomorphic to
the pull-back
of $\L_E^{d'}\boxtimes \L_E^{d''}$ under $Y'\to \Coh_0^{d'}\times
\Coh_0^{d''}$.

\medskip

We have the following diagram:
$$
\begin{CD}
Y'   @<<< Y'\underset{\Coh_0^{d'}\times
\Coh_0^{d''}}\times\on{Fl}^{d',d''}_0  @<<<
Y  \\
@VVV     @VVV   \\
\Coh_0^{d'}\times \Coh_0^{d''}  @<{\mathfrak q}<<  \on{Fl}^{d',d''}_0,
\end{CD}
$$
in which the composed upper horizontal map is $\on{v}$. Moreover, it
is easy to see that the map $Y'\underset{\Coh_0^{d'}\times
\Coh_0^{d''}}\times\on{Fl}^{d',d''}_0 \leftarrow Y$ is a fibration
with fibers being affine spaces of the same dimension. Therefore, up
to Tate's twist and a cohomological shift, $\on{v}_!(\wt\W_{E,k}|_Y)$
is isomorphic to the pull-back under $Y'\to \Coh_0^{d'}\times
\Coh_0^{d''}$ of ${\mathfrak q}_!\circ {\mathfrak p}^*(\L^d_E)$.

Hence, the assertion of the lemma follows from \eqref{simple Laumon}.

\end{proof}

\ssec{End of the proof of formula \eqref{rho vanishing}}

The above considerations show that it suffices to treat the case when
$\phi:\M^0\to \T$ is
surjective. Let $\M^1$ denote the kernel of $\phi$. Let us observe that
the scheme $Y$ can be identified with the scheme $\Hom^0(\M^1,\M_{tf})$
of injective maps $\M^1\to \M_{tf}$.

Indeed, to
$$0\to \M_{tf}\to \M\to \T\to 0,\ \  \M^0\hookrightarrow \M$$
we associate $\M^1\hookrightarrow \M^0\hookrightarrow \M$,
which maps into $\M_{tf}$ by assumption. And vice versa:
to an embedding $\M^1\to \M_{tf}$ we associate
$\M=\M_{tf}\underset{\M^1}\oplus \M^0$.

Moreover, the sheaf $\wt\W_{E,k}|_Y$ becomes isomorphic to $\pi^*(\L^d_E)$,
where $d=\deg(\M_{tf})-\deg(\M^1)$. Therefore, the cohomology
$H_c(Y,\wt\W_{E,k}|_Y)$ equals the cohomology appearing in
\conjref{pointwise vanishing}.

By assumption, the vector bundle $\M_{tf}$ satisfies:
$\Hom(\M_{tf},\L^{\on{est}})=0$. Hence, by the condition on
$\L^{\on{est}}$ (cf. \secref{seccon}), $\deg(\M_{tf})>n k (2g-2)$, and
so $d=\deg(\M_{tf})+\ell> n k (2g-2)$.

The required vanishing statement now follows from
\conjref{pointwise vanishing} applied to $E$.  This completes the
proof of formula \eqref{rho vanishing} and \thmref{clean}.

\section{Descent of the sheaf $\F_{E,n}$} \label{descent}

As in the previous section, we keep the assumption that the local
system $E$ is irreducible and that \conjref{vanishing conjecture}
holds for $E$. Our goal here is to prove \thmref{descends}.

Having established \thmref{clean} for all $k=1,...,n-1$ we know,
according to \corref{clean1}, that over $\C_n^d$ the complex
$\F_{E,n}\otimes \Ql(\frac{d}{2})[d]$ is an irreducible perverse
sheaf.




\ssec{Euler characteristics}

The morphism $\rho^0_n: \E^0_n \to \C_n$ is smooth of relative
dimension $d-n^2(g-1)$, and the sheaf $\F_{E,n}\otimes
\Ql(\frac{d}{2})[d]|_{\E_n^0\cap (\rho^0_n)^{-1}(\C_n^d)}$ is perverse
and irreducible. Hence in order to prove \thmref{descends} it suffices
to show that when $d\geq c_{g,n}$, the restriction
$\F_{E,n}|_{(\rho_n^0)^{-1}(\C_n\cap\Bun_n^d)}$ is non-zero and that
it descends to a perverse sheaf on $\C_n\cap\Bun_n^d$.


Recall from
\lemref{comp2&3} that
\begin{equation*}
\F_{E,n}|_{(\rho_n^0)^{-1}(\C_n\cap\Bun_n^d)}\simeq
\Aut'_E|_{(\rho_n^0)^{-1}(\C_n\cap\Bun_n^d)},
\end{equation*}
up to a cohomological shift.

The proof will be based on the following proposition:

\begin{prop}  \label{euler constant}
The Euler characteristics of the stalks of $\Aut'_E|_{\Bun'_n}$ are
constant along the fibers of the projection $\varrho_n:\Bun'_n\to
\Bun_n$. Moreover, they are not identically equal to zero over
$(\varrho_n)^{-1}(\C_n\cap\Bun^d_n)$ (for $d\geq c_{g,n}$).
\end{prop}

\ssec{Derivation of descent from \propref{euler constant}}

Suppose $d \geq c_{g,n}$. Then the perverse sheaf
$\F_{E,n}|_{\E^0_n \cap (\rho^0_n)^{-1}(\C_n^d)}\otimes
\Ql(\frac{d}{2})[d]$ is the Goresky-MacPherson extension of a local
system on a locally closed substack $U'$ of $\E^0_n$, contained inside
the open substack $(\rho_n^0)^{-1}(\C_n\cap\Bun_n^d) \subset \E^0_n$.

There exists a smooth locally closed substack $U_1 \subset
\C_n\cap\Bun^d_n$, such that if we set $U'_1=(\rho_n^0)^{-1}(U_1)$,
the intersection $U'_1\cap U'$ is open and dense in $U'$ and
$\rho^0_n: U'_1\cap U'\to U_1$ is surjective. Since $\rho_n^0$ is
smooth, in order to prove \thmref{descends}, it suffices to show that
$\F_{E,n}\otimes \Ql(\frac{d}{2})[d]|_{U'_1}$ is a pull-back of a
local system on $U_1$.

We have the following general result:

\begin{lem}
Let $Y$ be a smooth scheme (or stack) and let $\K$ be an irreducible
perverse sheaf on $Y$. If the Euler characteristics of the stalks of
$\K$ are the same at all $\ol\kk$--points of $Y$, then $\K$ is a local
system. If these Euler characteristics are not identically equal to
$0$, then $\K \neq 0$.
\end{lem}

\begin{proof}
Let $Y_0\subset Y$ be the maximal open subset over which $\K$ is a
local system. By the irreducibility assumption, $\K$ is the
Goresky-MacPherson extension of a local system on $Y_0$. Since $Y$ is
smooth, it is enough to show that $Y-Y_0$ is of codimension $\geq 2$.

Suppose this is not so. Then $Y-Y_0$ contains a divisor. Let $A$
denote the strict Henselization of the local ring at the generic point
of this divisor and let $\eta$ (resp., $s$) be the generic (resp.,
closed) point of $\on{Spec}(A)$.

By our assumptions, $\K|_{\on{Spec}(A)}$ is the Goresky-MacPherson
extension of a local system on $\eta$. The stalk $\K_\eta$ is a
representation of the Galois group $\Gamma$ of the field of fractions
of $A$, and we have $\K_s\simeq (\K_\eta)^\Gamma$.

By the assumption on the Euler characteristics,
$\on{dim}(\K_s)=\on{dim}(\K_\eta)$, i.e., the representation of
$\Gamma$ on $\K_\eta$ is trivial. But this means that $\K$ extends as
a local system to the entire $\on{Spec}(A)$, which is a contradiction.

\end{proof}

Let us show that this lemma is applicable for $Y=U'_1$ and
$\K=\F_{E,n}\otimes \Ql(\frac{d}{2})[d]|_{U'_1}$.  Indeed, the Euler
characteristics of stalks of $\F_{E,n}|_{U'_1}$ are constant along the
fibers of $U'_1\to U_1$, by \propref{euler constant}.  Moreover, they
are constant on $U'_1\cap U'$, since $\F_{E,n}$ is a local system
there. Hence, the Euler characteristics are constant on all of $U'_1$,
since $U'_1\cap U'\to U_1$ is surjective.

Thus, we obtain that $\F_{E,n}\otimes
\Ql(\frac{d}{2})[d]|_{U'_1}$ is a non-zero local system.

\medskip

By definition, $U'_1$ is a complement to the zero section in the vector
bundle $\E_n|_{U_1}$. By construction, $\F_{E,n}$ is equivariant with
respect to the natural $\GG_m$--action along the fibers of the
projection $\rho^0_n:U'_1\to U_1$. Therefore, the assertion of
\thmref{descends} follows from the next lemma applied to
$\E:=\E_n|_{U_1}$, $\K:=\F_{E,n}\otimes \Ql(\frac{d}{2})[d]|_{U'_1}$.

\begin{lem}
Let $\E\to Y$ be a vector bundle and let us denote by $\E^0$ the
complement to the zero section.  Let $\K$ be a local system on $\E^0$,
equivariant with respect to the $\GG_m$-action along the fibers.  Then
$\K$ descends to a local system on $Y$.
\end{lem}

\begin{proof}
This follows from the fact that any local system on a projective space
is isomorphic to the trivial local system.
\end{proof}

Now we prove \propref{euler constant}. The first step is the following
statement.

\begin{lem}    \label{two local syst}
Let $E'$ be another rank $n$ local system, not necessarily
irreducible. Then the Euler characteristics of $\Aut'_E$
and $\Aut'_{E'}$ are equal at any given $\ol\kk$--point of $\Bun'_n$.
\end{lem}

In order to prove the lemma, we will use the following corollary of a
theorem of Deligne from \cite{Ill}, Corollary 2.10:

\begin{thm}    \label{del euler}
Let $f:Y_1\to Y_2$ be a proper morphism of schemes (or a proper
representable map of stacks). Let $\K$ and $\K'$ be two complexes on
$Y_1$, which are locally isomorphic, by which we mean that they can be
represented as inverse limits of \'etale-locally isomorphic complexes
with torsion coefficients.  Then $f_!(\K)$ and $f_!(\K')$ have equal
Euler characteristics at every $\ol\kk$-point of $Y_2$.
\end{thm}

Recall the stack $\ol\Q$ and note that the group $(\GG_m)^n$ acts on
it by the rule
\begin{equation*}
(c_1,...,c_n)\cdot (\M,s_1,...,s_n)=(\M,c_1\cdot s_1,...,c_n\cdot s_n).
\end{equation*}
Consider the quotient $\ol\Q_r:=\ol\Q/(\GG_m)^{n-1}$, where
$(\GG_m)^{n-1}\subset (\GG_m)^n$ corresponds to the omission of the
first copy of $\GG_m$. Then the morphism $\nu:\ol\Q\to \Bun'_n$
factors as
\begin{equation*}
\ol\Q\to \ol\Q_r \overset{\nu_r}\to \Bun'_n.
\end{equation*}
The following is proved in \cite{BG1}, Proposition 1.2.2:

\begin{lem}    \label{represen}
The morphism $q_r: \ol\Q/(\GG_m)^n\to \Bun_n$ is representable and
proper.
\end{lem}

We obtain from this lemma that $\nu_r:\ol\Q_r\to \Bun'_n$ is also
proper.

\medskip

Now let us take the quotient of the diagram \eqref{diagram for whit} by
$(\GG_m)^{n-1}$:
\begin{equation*}
\begin{CD}
\ol\Q^0_r  @<{'\hl_r}<<  Z_r^d  @>{'\hr_r}>>  \ol\Q^d_r   \\
@V{q_r}VV      @V{'q_r}VV    @V{q_r}VV   \\
\Bun^0_n   @<{\hl}<<  \Mod_n^d  @>{\hr}>>  \Bun^d_n
\end{CD}
\end{equation*}

Denote by $\Psi^0_r$ the $!$--direct image of $\Psi^0$ under
$\ol\Q^0\to \ol\Q^0_r$.  It is clear that $\Aut'_E$ can be written as
\begin{equation} \label{Radon}
\Aut'_E = \nu_r{}_! \; {}'\hr_r{}_! ({}'\hl_r{}^*(\Psi^0_r)\otimes
{}'q_r^* \pi^*(\L_E^d))\otimes \Ql(\frac{d\cdot n}{2})[d\cdot n],
\end{equation}
and similarly for $E'$.

This formula and \thmref{del euler} readily imply the equality of
the Euler characteristics of $\Aut'_E$ and $\Aut'_{E'}$. Indeed, the
morphism $\nu_r\circ {}'\hr_r:Z_r^d\to \Bun_n'$ is proper and the
complexes
$$'\hl_r{}^*(\Psi^0_r)\otimes {}'q_r^* \pi^*(\L_E^d)), \qquad
'\hl_r{}^*(\Psi^0_r)\otimes {}'q_r^* \pi^*(\L_{E'}^d))$$ are locally
isomorphic, since so are the corresponding Laumon's sheaves $\L^d_E$
and $\L_{E'}^d$. This completes the proof of \lemref{two local syst}.

\medskip

We remark that formula \eqref{Radon} is the generalization of the
Radon transform construction (as opposed to the Fourier transform
construction from \secref{seccon}) in Drinfeld's original proof
\cite{Dr} of the Langlands conjecture in the case of $GL_2$.

\ssec{Conclusion of the proof of \propref{euler constant}}

According to \lemref{two local syst}, in order to prove \propref{euler
constant} it suffices to show that there exists at least one local
system $E$, for which the statement of \propref{euler constant} is
true. Hence it suffices to prove it for the trivial local
system. Using the reduction technique of \cite{BBD}, Sect. 6.1.7, we
obtain:

\begin{lem}
Suppose that \propref{euler constant} holds when $E$ is the trivial
local system in the case when the ground field $\kk$ is a finite field
of characteristic $p$. Then \propref{euler constant} holds when $E$ is
the trivial local system in the case of an arbitrary field $\kk$ of
the same characteristic.
\end{lem}

\begin{proof}

Let $s_1$ and $s_2$ be two $\overline\kk$-points of $\Bun_n'$, which
project to the same point of $\Bun_n$.  First, we can assume that all
our data are defined over an algebra $A$ finitely generated over a
finite field; $A\subset \overline\kk$. In other words, we have the
stacks $(\Bun_n)_A$ and $(\Bun_n')_A$ and sections
$s_i:\on{Spec} A \to (\Bun_n')_A$.

Consider $s_i^*(\Aut'_{E_0}), i=1,2$, where $E_0$ is the trivial
local system (so that it is defined over $A$), as $\ell$-adic
complexes on $\on{Spec} A$.  By localizing $A$ we may assume that it
is smooth over a finite field, and that the above complexes are
locally constant.

Let $\eta: \on{Spec} \overline\kk \to \on{Spec} A$ be the canonical
generic geometric point of $\on{Spec} A$ and let $a: \on{Spec} \Fq \to
\on{Spec} A$ be some closed geometric point of $\on{Spec} A$. We need
to compare the Euler characteristics of the stacks
$(s_i^*(\Aut'_{E_0}))_\eta$ for $i=1,2$. Since our complexes are locally
constant, we may instead compare the stalks $(s_i^*(\Aut'_E))_a$.  In
other words, we can make the comparison over the finite field, as
required.
\footnote{Note that if the ground field $\kk$ is of characteristic
$0$, we can use a similar argument by choosing $A$ to be a
finitely generated algebra over $\ZZ$.}

\end{proof}

Thus, it suffices to prove \propref{euler constant} for the trivial
local system in the case when $\kk$ is a finite field.

\smallskip

Let us apply \lemref{two local syst} again and obtain that it
suffices to find just one local system $E^s$ in the case when $\kk$ is
a finite field, for which \propref{euler constant} is true.
We will take as $E^s$ any irreducible local system, which satisfies
the following conditions: \newline
\noindent (a) $E^s$ is pure, and

\smallskip

\noindent (b) there exists a cuspidal Hecke eigenfunction associated
to the pull-back of $E^s$ to $X \underset{\Fq}\times {\mathbb F}_{q_1}$
for any finite extension ${\mathbb F}_{q_1}$ of $\Fq$.

\smallskip

For example, such a local system can be constructed as follows: pick a
cyclic $n$--sheeted \'etale cover $\wt X\to X$, and let $E^s$ be the
direct image of a generic rank one local system on $\wt X$ of
finite order. Then $E^s$ is pure, so condition (a) is
satisfied. Moreover, according to Theorem 6.2 of \cite{AC} (see also
\cite{Kaz}), this local system also satisfies condition
(b).\footnote{Actually, Lafforgue's results \cite{Lf} imply that any
irreducible local system $E$ satisfies conditions (a) and (b), up to a
twist with a rank one local system.}

Thus, we have at our disposal at least one irreducible rank $n$ local
system $E^s$, for which the above conditions (a), (b), as well as
\conjref{vanishing conjecture} are true. We now prove that then
\propref{euler constant} also holds for this $E^s$. To prove the first
assertion of \propref{euler constant}, it suffices to show that the
function $\text{\tt f}_{q_1}(\Aut'_{E^s})$ (obtained by taking the
traces of Frobenius on the stalks of $\Aut'_{E^s}$, see
\secref{conventions}) on $\Bun'_n({\mathbb F}_{q_1})$ is constant
along the fibers of the projection
\begin{equation*}
\Bun'_n({\mathbb F}_{q_1}) \to \Bun_n({\mathbb F}_{q_1})
\end{equation*}
for all finite extensions ${\mathbb F}_{q_1}$ of $\Fq$.





But Theorem 3.1 of \cite{FGKV} states that if a cuspidal Hecke
eigenfunction associated to any given rank $n$ local system
$\wt{E}$ exists on $\Bun_n({\mathbb F}_{q_1})$, then its pull-back to
$\Bun'_n({\mathbb F}_{q_1})$ equals $\text{\tt
f}_{q_1}(\Aut'_{\wt{E}})$ up to a non-zero scalar.

Applying these results to our local system $E^s$, we obtain that the
function $\text{\tt f}_{q_1}(\Aut'_{E^s})$ is constant along the
fibers of $\varrho_n: \Bun'_n({\mathbb F}_{q_1}) \to \Bun_n({\mathbb
F}_{q_1})$. This proves the
first assertion of \propref{euler constant} for $E^s$.

It remains to prove the non-vanishing assertion of \propref{euler
constant} for $E^s$. According to \propref{let E}, if $E$ is an
irreducible local system on a curve $X$ over a finite field, which
satisfies the above conditions (a) and (b), then \conjref{vanishing
conjecture} holds for $E$. Hence by our assumptions on $E^s$,
\conjref{vanishing conjecture} holds for $E^s$. Therefore by
\thmref{clean} the restriction of $\Aut'_{E^s}$ to the preimage of
$\C_n\cap \Bun^d_n$ in $\Bun_n'$ is a perverse sheaf, up to a
cohomological shift. Hence it suffices to show that this restriction
is non-zero (for if a perverse sheaf has zero Euler characteristics
everywhere, then this sheaf is zero).

For that, it is enough to show that the corresponding function does
not vanish identically on $(\C_n\cap \Bun^d_n)(\Fq)$, if $d\geq
c_{g,n}$ (note that in this case $\C_n\cap \Bun^d_n \neq
\emptyset$). However, by assumption, $$\Bun^d_n-(\C_n\cap
\Bun^d_n)\subset \Bun_n^{\on{vuns}}, \qquad d\geq c_{g,n}$$ (see
\secref{notat} for the definition of $\Bun_n^{\on{vuns}}$). The
definition of cuspidal function implies the following

\begin{lem}  \label{cusp support}
Let $f$ be a cuspidal function on $\Bun_n(\Fq)$.
Then its restriction to $\Bun_n^{\on{vuns}}(\Fq)$ is identically zero.
\end{lem}

\begin{proof}

Let $\M\in \Bun_n(\Fq)$ be a very unstable bundle, and let $\M \simeq
\M_1\oplus\M_2$ be the corresponding decomposition, with
$\rk(\M_i)=n_i$.

Let $r^n_{n_2,n_1}:\on{Funct}(\Bun_n(\Fq))\to
\on{Funct}(\Bun_{n_2}(\Fq)\times \Bun_{n_1}(\Fq))$ be the
corresponding constant term operator.  Since $f$ is cuspidal, we have
$r^n_{n_2,n_1}(f)=0$. However, by applying the definition of
$r^n_{n_2,n_1}$ and evaluating $r^n_{n_2,n_1}(f)$ at the point
$\M_2\times \M_1\in \Bun_{n_2}(\Fq)\times \Bun_{n_1}(\Fq)$, we obtain
that it is equal to the integral
$$\underset{0\to \M_2\to \M'\to \M_1\to 0}\int f(M'),$$ over the
finite set $\on{Ext}^1(\M_1,\M_2)(\Fq)$ (the measure on this set is a
non-zero multiple of the tautological measure).

However, by our assumption, $\on{Ext}^1(\M_1,\M_2)=0$, therefore
$r^n_{n_2,n_1}(f)(\M_2,\M_1)=f(M)$, up to a non-zero constant.
\end{proof}

Thus, we obtain the second assertion of \propref{euler constant} for
our local system $E^s$. This completes the proof of \thmref{descends}.

\section{The Hecke property of $\Aut_E$}    \label{Hecke property}

In the previous section we constructed a perverse sheaf $\oS_E$ on
$\C_n$, whose pull-back to $\E^0_n$ is
$\F_{E,n}\otimes \Ql(\frac{n^2\cdot (g-1)}{2})[n^2\cdot (g-1)]$.

Let $\S_E$ be the Goresky-MacPherson extension of $\oS_E$ to
$\underset{d\geq c_{g,n}}\cup\, \Coh^d_n$. Finally, set
$\Aut_E:=\S_E|_{\underset{d\geq c_{g,n}}\cup\,\Bun^d_n}$.

Our goal is to prove the following

\begin{thm}   \label{Hecke theorem}
The perverse sheaf $\Aut_E$ can be uniquely extended to the entire
stack $\Bun_n$, so that it becomes a Hecke eigensheaf with respect to $E$.
\end{thm}

\thmref{Hecke theorem} will follow from \propref{firstHecke}, as will be
explained in \secref{proof of hecke}. We will give two independent
proofs of \propref{firstHecke}. The first one, presented in
Sects. \ref{proof 1}--\ref{proof last}, uses the Whittaker sheaf
$\W_E$. The second proof, given in \secref{hecke laumon}, uses the
Hecke-Laumon property of the Laumon sheaf $\L_E$.

\ssec{The Hecke property on $\Bun_n'$}    \label{proof 1}

Consider the Cartesian product $\Bun'_n\underset{\Bun_n}\times
\H^1_n$, where the map $\H^1_n\to \Bun_n$ is $\hr$. We have a
commutative diagram, in which the right square is Cartesian:
\begin{equation}   \label{diag Bun'}
\begin{CD}
\Bun'_n @<{''\hl}<< \Bun'_n\underset{\Bun_n}\times \H^1_n
@>{\supp\times{}''\hr}>> X\times \Bun_n'  \\
@VVV  @VVV    @VVV  \\
\Bun_n  @<{\hl}<< \H^1_n  @>{\supp\times\hr}>> X\times \Bun_n,
\end{CD}
\end{equation}
where the morphisms $''\hl$ and $''\hr$ are given by
\begin{eqnarray*}
''\hl: & (x,\M,\M',\beta:\M'\hookrightarrow \M,s':\Omega^{n-1} \to
\M') &\mapsto (\M,s=\beta\circ s':\Omega^{n-1}\to \M),\\
''\hr: & (x,\M,\M',\beta:\M'\hookrightarrow \M,s':\Omega^{n-1} \to
\M') &\mapsto (\M',s').
\end{eqnarray*}

\begin{prop} \label{Hecke on Bun'}
For any local system $E$ of rank $n$,
\begin{equation}    \label{isom of Aut'}
(\supp\times{}''\hr)_! {}''\hl{}^*(\Aut'_E)\otimes
\Ql(\frac{n-2}{2})[n-2] \simeq E \boxtimes \Aut'_E.
\end{equation}
\end{prop}

First, we will reformulate this proposition in terms of the stack
$\ol\Q^d$, introduced in \secref{Drinfeld's compact}.

\ssec{A reformulation}

We need to introduce two more stacks $\Qp^d$ and ${}\Qpp^d$ closely
related to $\ol\Q^d$. The stack ${}\Qpp^d$ classifies the data of
$(x,\M,(s_i))$ as in the definition of $\ol\Q^d$, but with $$s_i:
\Omega^{(n-1)+...+(n-i)}\to (\Lambda^i \M)(x), \qquad i=1,\ldots,n.$$
The stack $\Qp^d$ classifies the same data with the additional
condition that the image of $s_1$ is contained in $\M$ (and not just
$\M(x)$). We have tautological closed embeddings
\begin{equation*}
X\times\ol\Q^d \hookrightarrow \Qp^d \hookrightarrow {}\Qpp^d.
\end{equation*}

Recall that we have a forgetful morphism $\ol\Q^{d+1} \to \Bun'_n$,
and the morphism $''\hl: \Bun'_n\underset{\Bun_n}\times \H^1_n \to
\Bun'_n$. Denote by $\ol{Q}\H^{d+1}$ the corresponding fiber
product. Consider the following commutative diagram:
\begin{equation}    \label{diag for Q}
\begin{CD}
\ol\Q^{d+1} @<{\thl}<< \ol\Q^{d+1} \underset{\Bun_n}\times \H^1_n
@>{\thr}>> \Qpp^d \\
@AAA  @AAA  @AAA  \\
\ol\Q^{d+1} @<<< \ol{Q}\H^{d+1}
@>>> \Qp^d \\
@V{\nu}VV  @VVV    @V{\on{id}\times \nu}VV  \\
\Bun'_n @<{''\hl}<< \Bun'_n\underset{\Bun_n}\times \H^1_n
@>{\supp\times{}''\hr}>> X\times \Bun'_n
\end{CD}
\end{equation}
The bottom left and the top right squares in this diagram are Cartesian.

By definition, $\Aut'_E = \nu_!({\mathcal W}_E)$, where $\nu: \ol\Q
\to \Bun'_n$ is the forgetful morphism defined in \secref{Aut'}. Note
that $\nu: \ol\Q^d \to \Bun'_n$ extends to a morphism $\Qp^d\to
\Bun'_n$, which we also denote by $\nu$. Using the diagram \eqref{diag
for Q}, we obtain that the LHS of formula \eqref{isom
of Aut'}, restricted to the degree $d$ connected component of
$\Bun_n'$, is isomorphic to the complex $(\on{id}\times\nu)_!(\Wp^d)$, where
\begin{equation}   \label{restr}
\Wp^d := \thr_! \, \thl{}^*(\W_E^{d+1})|_{\Qp^d}\otimes
\Ql(\frac{n-2}{2})[n-2].
\end{equation}

Therefore \propref{Hecke on Bun'} follows from \propref{Hecke on Q},
which is proved in the Appendix.

\begin{prop} \label{Hecke on Q}
The complex $\Wp^d$ is supported on $X\times \ol\Q^d
\subset \Qp^d$, and its restriction to $X\times \ol\Q^d$ is isomorphic
to $E\boxtimes {\mathcal W}^d_E$.
\end{prop}

\ssec{The Hecke property on $\Bun_n$}    \label{proof last}

Observe that in the diagram \eqref{diagram H1} defining the Hecke
functor $\He^1_n$ we have:
$$(\on{supp} \times \hr)^{-1}(X\times (\C_n\cap \Bun_n)) \subset
(\hl)^{-1}(\C_n\cap \Bun_n).$$ Therefore we can define a functor
$\on{D}(\C_n\cap \Bun_n)\to\on{D}(X\times (\C_n\cap \Bun_n))$ by
formula \eqref{formula H1}.  We denote this functor also by $\He^1_n$
and consider its iterations $(\He^1_n)^{\boxtimes i}$ and the
corresponding functors $\He^i_n$. The notion of Hecke eigensheaf also
makes sense in this context.

We now derive from \propref{Hecke on Bun'} the following

\begin{prop}  \label{firstHecke}
The perverse sheaf $\Aut_E|_{\C_n\cap \Bun_n}$ is a Hecke eigensheaf
with respect to $E$.
\end{prop}

\begin{proof}

Recall from \secref{seccon} and \lemref{comp2&3}
that over $\E^0_n\cap (\varrho_n)^{-1}(\Bun_n^d)$ we have an isomorphism
$$\F_{E,n}\simeq \Aut'_E\otimes \Ql(\frac{-d+c}{2})[-d+c],$$
where $c$ is a constant depending only on $g$ and $n$.

The isomorphism of $\He^1_n(\Aut_E)$ and $E \boxtimes \Aut_E$ over
$\C_n\cap \Bun_n$ now follows from
\propref{Hecke on Bun'} via diagram \eqref{diag Bun'}, using the fact
that the morphism $\rho^0_n: (\Bun'_n \cap \E^0_n) \to (\Bun_n \cap
\C_n)$ is smooth, representable and has connected fibers.

Moreover, it follows from the construction of the isomorphism of
\propref{Hecke on Q}, that this isomorphism
satisfies condition \eqref{eigen-equivariance}.
\end{proof}

Now we derive \thmref{Hecke theorem} from the above proposition.

\ssec{Proof of \thmref{Hecke theorem}}    \label{proof of hecke}

Recall the morphism $\on{mult}:X\times \Bun_n\to \Bun_n$ given by
$(x,\M)\mapsto \M(x)$.  In the same way as in the proof of
\propref{thesis}, we obtain from \propref{firstHecke} that there is an
isomorphism
\begin{equation} \label{n-th Hecke}
\on{mult}^*(\Aut_E)|_{X\times (\C_n\cap \Bun_n)}\simeq \Lambda^n E
\boxtimes \Aut_E|_{X\times (\C_n\cap \Bun_n)}.
\end{equation}
Since the morphism $\on{mult}:X\times \Bun_n\to \Bun_n$ is smooth, the
isomorphism of formula \eqref{n-th Hecke} holds over the entire component
$\Bun^d_n$ for $d\geq c_{g,n}$ (and not only over $\C_n\cap\Bun^d_n$).

Now we extend $\Aut_E$ to all other connected components of $\Bun_n$
as follows: for every open substack $U\subset \Bun_n^{d'}$ of finite
type, there exists an integer $d''$ such that for any $x\in X$, the
morphism $\on{mult}_{d''\cdot x}: \Bun_n \to \Bun_n$ sending $\M$ to
$\M(d''\cdot x)$, maps $U$ into $\C_n\cap \Bun_n$. We set $\Aut_E|_U$
to be
\begin{equation*}
\on{mult}_{d''\cdot x}^*(\Aut_E)\otimes (\Lambda^n E_x)^{\otimes -d''}
\end{equation*}
According to formula \eqref{n-th Hecke}, this gives a well-defined
sheaf $\Aut_E$ on the entire $\Bun_n$, together with an isomorphism
$\on{mult}^*(\Aut_E)\simeq \Lambda^n E \boxtimes \Aut_E$.

\propref{firstHecke} then implies that $\Aut_E$ is a Hecke
eigensheaf. Indeed, the existence and uniqueness of the isomorphism
\eqref{eigen-property} satisfying \eqref{eigen-equivariance} over the
entire $\Bun_n$ follow from the construction, using the fact that
formula \eqref{Hecke commute} holds over $\C_n\cap \Bun_n$.

This completes the proof of \thmref{Hecke theorem}.

\ssec{Lifting of $\Aut_E$ to $\Bun'_n$}     \label{lifting to Bun'}

We have the sheaves $\Aut_E$ on $\Bun_n$ and $\Aut'_E$ on
$\Bun'_n$. Consider the commutative diagram
\begin{equation*}
\begin{CD}
\Bun'_n \cap \E^0_n @>>> \Bun'_n \\
@V{\rho^0_n}VV    @V{\varrho_n}VV  \\
\Bun_n \cap \C_n @>>> \Bun_n,
\end{CD}
\end{equation*}
By construction, for $d\geq c_{n,g}$, the sheaves
$\varrho_n{}^*(\Aut_E)$ and $\Aut'_E\otimes (\frac{-d+c}{2})[-d+c]$
are isomorphic over $\E^0_n\cap (\varrho_n)^{-1}(\Bun_n^d)$, where $c$
is a constant independent of $d$.  In this subsection we will address
the following question, posed by V.~Drinfeld:

\medskip

{\it Are the sheaves $\varrho_n{}^*(\Aut_E)$ and $\Aut'_E\otimes
(\frac{-d+c}{2})[-d+c]$ isomorphic on the entire $\Bun'_n$?}

\medskip

The answer is affirmative. Indeed, consider the diagram
\begin{equation*}
\begin{CD}
\Bun'_n @<{''\hl}<< \Bun'_n\underset{\Bun_n}\times \H^i_n
@>{\supp\times{}''\hr}>> X\times \Bun_n'  \\
@VVV  @VVV    @VVV  \\
\Bun_n  @<{\hl}<< \H^i_n  @>{\supp\times\hr}>> \Bun_n,
\end{CD}
\end{equation*}
defined in the same way as diagram \eqref{diag Bun'}. From
\propref{Hecke on Bun'} we derive, in the same way as in
\propref{thesis}, that
\begin{equation*}
(\supp\times{}''\hr)_! {}''\hl{}^*(\Aut'_E)\otimes
\Ql(\frac{i(n-i-1)}{2})[i(n-i-1)]
\simeq
\Lambda^i E \boxtimes \Aut'_E.
\end{equation*}

In addition, from the Hecke property of $\Aut_E$ it follows that
\begin{equation*}
(\supp\times{}''\hr)_! {}''\hl{}^*(\varrho_n{}^*(\Aut_E))\otimes
\Ql(\frac{i(n-i)}{2})[i(n-i)]
\simeq \Lambda^i E \boxtimes \rho^0_n{}^*(\Aut_E).
\end{equation*}

As before, for $i=n$, the functor $\K\mapsto (\supp\times{}''\hr)_!
{}''\hl{}^*(\rho^0_n{}^*(\K))$ amounts to the pull-back under the map
$$\on{mult}':X\times \Bun_n'\to \Bun_n'$$ given by
$$
(x,\M,s:\Omega^{n-1}\to \M)\mapsto (\M(x),s':\Omega^{n-1}\to \M\to
\M(x)).
$$

Any open substack $U$ of finite type in
$(\varrho_n)^{-1}(\Bun_n^{d'})$ can be mapped into $\E^0_n\cap
(\varrho_n)^{-1}(\Bun^d_n))$ with $d\geq c_{g,n}$ by means of
$\on{mult}'_{d''\cdot x}:\Bun'_n\to \Bun_n'$. Hence, over
$(\varrho_n)^{-1}(U)$ we have:
\begin{align*}
&\Aut'_E\otimes (\Lambda^n E_x)^{\otimes d'}\otimes
\Ql(\frac{-d'+c}{2})[-d'+c]\simeq \\
&\simeq \on{mult}'_{d''\cdot x}{}^*(\Aut'_E)\otimes
\Ql(\frac{-d'-n\cdot d''+c}{2})[-d'-n\cdot d''+c] \simeq \\
&\simeq \on{mult}'_{d''\cdot x}{}^*(\rho^0_n{}^*(\Aut_E)) \simeq
\varrho_n{}^*(\Aut_E)\otimes (\Lambda^n E_x)^{\otimes d''}.
\end{align*}

The fact that the constructed isomorphism does not depend on the
choice of $x$ and $d'$ follows in the same way as the corresponding
assertion for $\Aut_E$ in the proof of \thmref{Hecke theorem}.

\section{The Hecke--Laumon property of $\S_E$}    \label{hecke laumon}

In this section we give an alternative proof of \propref{firstHecke},
and hence of \thmref{Hecke theorem}.

Consider the diagram
\begin{equation}    \label{diagram HL}
\Coh_n \xleftarrow{\hl_l} \H\L^d_n \xrightarrow{\hr_l} \Coh^d_0\times\Coh_n,
\end{equation}
where the stack $\H\L^d_n$ classifies short exact sequences $0\to \M'
\to \M\to \T\to 0$ with $\M' \in \Coh_n$, $\T\in \Coh_0^d$. The
projections $\hl_l$ and $\hr_l$ send such data to $\M$ and $(\M',\T)$,
respectively. (Recall that in \secref{section rho vanishing} we
encountered this stack
for $n=0$ and called it $\on{Fl}_d^{d',d''}$).

The {\em Hecke--Laumon} functor $\HL^d_n: \on{D}(\Coh_n)
\to \on{D}(\Coh^d_0\times\Coh_n)$ is defined by the formula
\begin{equation}
\begin{cases}
&\HL^d_n(\K) = \hr_l{}_! \hl_l{}^*(\K)\otimes \Ql(\frac{d\cdot
(n+1)}{2})[d\cdot (n+1)],\,\, n\geq 1  \\
&\HL^d_0(\K) = \hr_l{}_! \hl_l{}^*(\K),
\end{cases}
\end{equation}
(see \cite{La2}).

Note that for $d=d_1+d_2$ there is a natural isomorphism of functors
$$(\on{id}\times \HL_n^{d_2})\circ \HL_n^{d_1}\simeq
(\HL_0^{d_1}\times \on{id})\circ \HL_n^{d}.$$

Finally, let us note that \eqref{simple Laumon} stated in
\secref{section rho vanishing}
reads as
$$\HL_0^{d_1}(\L^d_E)\simeq \L^{d_1}_E\boxtimes \L^{d_2}_E.$$

\ssec{Definition} We say that a complex $\K\in\on{D}(\Coh_n)$ has a
{\it Hecke--Laumon property} (or is a Hecke--Laumon eigensheaf) with
respect to $E$ if for each $d$ we are given an isomorphism
\begin{equation} \label{eigen Laumon}
\HL^d_n(\K)\simeq \L^d_E\boxtimes \K
\end{equation}
such that for $d=d_1+d_2$ the diagram
\begin{equation}  \label{Laumon assoc}
\begin{CD} (\on{id}\times \HL_n^{d_2})\circ \HL_n^{d_1}(\K) @>\sim>>
\L_E^{d_1}\boxtimes \L_E^{d_2}\boxtimes \K \\ @V{\sim}VV @V{\sim}VV
\\ (\HL_0^{d_1}\times \on{id})\circ \HL_n^{d}(\K) @>\sim>>
\HL_0^{d_1}(\L^d_E)\boxtimes \K.  \end{CD}
\end{equation}
is commutative.

\ssec{Restriction to $\C_n$}

Note that the definition of the Hecke--Laumon property makes sense not
only on $\Coh_n$, but also on $\C_n$.  Consider the stack
$(\H\L_n^d)':=\H\L^d_n\underset{\Coh^d_0\times\Coh_n}\times
\Coh^d_0\times\E^0_n$
and note that it fits into the diagram
$$
\CD
\E^0_n @<{'\hl_l}<< (\H\L_n^d)' @>{'\hr_l}>> \Coh^d_0\times\E^0_n   \\
@V{\varrho_n}VV      @VVV   @V{\on{id}\times \varrho_n}VV  \\
\Coh_n @<{\hl_l}<< \H\L^d_n @>{\hr_l}>> \Coh^d_0\times\Coh_n.
\endCD
$$

It is shown in \cite{La2} (by induction on $k$) that
$$'\hr_l{}_!{}'\hl_l(\F_{E,n})\otimes \Ql(\frac{d\cdot
(n+1)}{2})[d\cdot (n+1)]\simeq \L^d_E\boxtimes \F_{E,n}.$$
Therefore, since $\rho^0_n: \E^0_n\to \C_n$ is smooth, representable
and with connected fibers, we obtain:

\begin{cor}
The perverse sheaf $\oS_E$ on $\C_n$ is a Hecke-Laumon eigensheaf with
respect to $E$.
\end{cor}

We will now prove the following result:

\begin{prop}  \label{Hecke-Laumon}
Let $\S$ be a perverse sheaf on $\Coh_n$ and $\DD(\S)$ its Verdier
dual sheaf. Suppose that $\S$ and $\DD(\S)$ satisfy the Hecke-Laumon
property with respect to local systems $E$ and $E^*$,
respectively. Then $\K:=\S|_{\Bun_n}$ is a Hecke eigensheaf with
respect to $E$.
\end{prop}

\begin{proof}

We start with the following general observation. Let $\rho:\E\to Y$ be
a vector bundle, $i:Y\to \E$ be the zero section and $j:\E^0\to \E$
its complement.  Let us denote by $\wt\rho:{\mathbb P}\E\to Y$ the
corresponding projectivized bundle.

Suppose that $\F$ is a $\GG_m$--equivariant perverse sheaf on $\E$,
and set $\F^0 := \F|_{\E^0}$. We will denote by $\wt\F$ the perverse
sheaf on ${\mathbb P}\E$ corresponding to $\F^0$, i.e. the pull-back of
$\wt\F$ to $\E^0$ is $\F^0\otimes \Ql(\frac{-1}{2})[-1]$. We have the
following assertion (see \cite{Ga}).

\begin{lem}   \label{GMextension}
Assume that $\rho_*(\F)[-1]$ and $\rho_!(\F)[1]$ are perverse
sheaves. Then $\wt\rho_!(\wt\F)$ is a perverse sheaf as well, and
$\rho_!(\F)\otimes \Ql(\frac{1}{2})[1]\simeq \wt\rho_!(\wt\F)\simeq
\rho_*(\F)\otimes \Ql(\frac{-1}{2})[-1]$.
\end{lem}

\begin{proof}
Since $\F$ is $\GG_m$--equivariant, $\rho_!(\F)\simeq i^!(\F)$ and
$\rho_*(\F)\simeq i^*(\F)$.  By applying $i^!$ to the triangle
$$j_! \F^0\to \F\to i_* i^*(\F),$$ we obtain that $\rho_!
j_!(\F^0)\simeq i^! j_!(\F^0)$ has perverse cohomology only in
cohomological degrees $0$ and $1$.

Using the Leray spectral sequence of the composition $\E^0\to {\mathbb
P}\E\to Y$, we obtain that $\wt\rho_!(\wt\F)$ must be perverse. In
addition, we obtain that $\wt\rho_!(\wt\F)\simeq h^0(\rho_!
j_!(\F^0))\otimes \Ql(\frac{-1}{2})$, which identifies
$\wt\rho_!(\wt\F)$ with $i^*(\F)\otimes
\Ql(\frac{-1}{2})[-1]$. Similarly, we obtain: $\wt\rho_!(\wt\F) \simeq
i^!(\F)\otimes \Ql(\frac{1}{2})[1]$.
\end{proof}

We will reduce the assertion of \propref{Hecke-Laumon} to the above
lemma. Set $Y=\Coh^1_0\times \Bun_n\subset \Coh^1_0\times \Coh_n$ and
take $\E$ to be
$$(\hr_l)^{-1}(\Coh^1_0\times \Bun_n) \subset \H\L^1_n.$$

(Note that in general the preimage $(\hr_l)^{-1}(\Coh^d_0\times \Bun_n)
\subset \H\L^d_n$ is the same as the stack $\widetilde\Coh_{n, \leq
d}$ introduced in \secref{induction on torsion} and the map
$\widetilde\Coh_{n, \leq d}\hookrightarrow \H\L^d_n\overset{\hl_l}\to
\Coh_n$ becomes the map $r:\widetilde\Coh_{n, \leq d} \to \Coh_{n,
\leq d}\subset \Coh_n$.)

Set $\F=\hl_l{}^*(\S)\otimes \Ql(\frac{n}{2})[n]$. Then $\F$ is
$\GG_m$--equivariant and perverse, according to \lemref{r smooth}.  In
addition, the image of $\E^0$ under $\hl_l$ lies in $\Bun_n\subset
\Coh_n$.

By the assumption of \propref{Hecke-Laumon}, $\S$ is a Hecke--Laumon
eigensheaf. This implies that $\rho_!(\F)\otimes
\Ql(\frac{1}{2})[1]\simeq \L^1_E\boxtimes \K$, and so $\rho_!(\F)[1]$
is a perverse sheaf. Applying Verdier duality and using the
assumptions of \propref{Hecke-Laumon} regarding $\DD(\S)$, we obtain
that $\rho_*(\F)[-1]$ is a perverse sheaf too. Hence, we can apply
\lemref{GMextension}.

Let us perform a base change with respect to $X\to \Coh^1_0$. Then
$X\underset{\Coh^1_0}\times {\mathbb P}\E$ identifies naturally with the
Hecke correspondence $\H^1_n$ in such a way that
$$\wt\rho: X\underset{\Coh^1_0}\times{\mathbb P}\E\to X\times \Bun_n$$
becomes the projection $\hr$. Therefore, \lemref{GMextension} implies that
$\He^1_n(\K)\simeq E\boxtimes \K$.

The fact that this
isomorphism indeed satisfies condition \eqref{eigen-equivariance}
follows from property \eqref{Laumon assoc} in the case $d=2$. This
completes the proof of \propref{Hecke-Laumon}.
\end{proof}

\ssec{Remark}

V.~Drinfeld has asked the following question about the possibility of
proving a theorem converse to \propref{Hecke-Laumon}:

{\it Let $\K$ be a perverse sheaf on $\Bun_n$, which is a Hecke
eigensheaf with respect to $E$.  Is it true that the
Goresky-MacPherson extension of $\K$ to $\Coh_n$ has the Hecke-Laumon
property with respect to $E$?}

We conjecture that the answer to this question is affirmative.

\ssec{Second proof of \propref{firstHecke}} It is clear from the above
proof that \propref{Hecke-Laumon} is still valid if we replace the
stacks $\Coh_n$ and $\Bun_n$ by their substacks $\C_n$ and $\C_n\cap
\Bun_n$, respectively. Now we apply this modification of
\propref{Hecke-Laumon} in the situation when $\S=\oS_E$, and
$\K=\Aut_E|_{\C_n\cap \Bun_n}$. It
follows from the definitions that all conditions of
\propref{Hecke-Laumon} are satisfied (in particular, we have:
$\DD(\oS_E) \simeq \oS_{E^*}$). The statement of \propref{firstHecke}
(and hence \thmref{Hecke theorem}) now follows directly from
\propref{Hecke-Laumon}.

\ssec{Lifting of $\S_E$ to $\Coh'_n$}  \label{lifts everywhere}

We have the diagram
\begin{equation*}
\begin{CD}
\E^0_n @>>> \Coh'_n \\
@V{\rho^0_n}VV    @V{\varrho_n}VV  \\
\C_n @>>> \Coh_n,
\end{CD}
\end{equation*}

Recall the definition of the perverse sheaf $\S_E$ on $\underset{d\geq
c_{n,g}}\cup\,\Coh^d_n$ given in the beginning of \secref{Hecke
property}.  By \thmref{descends} and \lemref{coincide}, over
$(\varrho_n)^{-1}(\C_n\cap (\underset{d\geq c_{n,g}}\cup\,\Coh^d_n))$
the complexes $\varrho_n^*(\S_E)$ and $\S'_E\otimes
\Ql(\frac{-d+c}{2})[-d+c]$ are isomorphic, where $c$ is a constant
independent of $d$.

Consider the Goresky-MacPherson extension of $\Aut_E$ (which by now
is defined on the whole of $\Bun_n$) to
$\Coh_n$. By abuse of notation we still denote this extension by
$\S_E$. Now we prove the following assertion:

\medskip

{\em The sheaf $\S_E$ has the Hecke-Laumon property with respect to
$E$, and $\varrho_n{}^*(\S_E)\simeq \S'_E\otimes
\Ql(\frac{-d+c}{2})[-d+c]$.}

\medskip

Let us denote by $(X\times \Coh_n)^0$ the open substack of $X\times
\Coh_n$, corresponding to those pairs $(x,\M)$, for which $\M$ has no
torsion supported at $x$. In a similar way we define the substack
$(X^i\times \Coh_n)^0$ of $X^i\times \Coh_n$.  We define the functors
\begin{equation}
\He^i_n:\on{D}(\Coh_n)\to \on{D}((X\times \Coh_n)^0)
\end{equation}
in the same way as before. Since we already know $\Aut_E$ is a Hecke
eigensheaf, we obtain that
\begin{equation}  \label{central Hecke-Laumon}
\He^n_n(\S_E)\simeq \Lambda^n(E)\boxtimes \S_E|_{(X^i\times \Coh_n)^0}.
\end{equation}

By arguing as in \secref{proof of hecke}, we deduce
the Hecke-Laumon property of $\S_E$ on the entire $\Coh_n$ from
\eqref{central Hecke-Laumon} and the fact that $\S_E|_{\C_n}$ is a
Hecke-Laumon sheaf.

Similarly, the isomorphism $\varrho_n{}^*(\S_E)\simeq
\S'_E\otimes \Ql(\frac{-d+c}{2})[-d+c]$ follows in the same way as
in \secref{lifting to Bun'}.

\section{Cuspidality}     \label{cuspidality}

\ssec{Constant term functors}

Let $P\subset GL_n$ be the standard (upper) parabolic subgroup
corresponding to a partition $(n_1,\ldots,n_k)$ of $n$, with the Levi
quotient $M \simeq GL_{n_1} \times \ldots \times GL_{n_k}$. The
embedding of $P$ in $GL_n$ and the projection $P\to M$ induce
morphisms $p$ and $q$ in the diagram
\begin{equation}    \label{cus diagram}
\Bun_n \xleftarrow {\ \ \p\ \ }\Bun_P \xrightarrow{\ \ \q \ \ }\Bun_M\,.
\end{equation}

The {\it constant term} functor $\on{R}^G_M: \on{D}(\Bun_n) \to
\on{D}(\Bun_M)$ is defined by the formula $\on{R}^G_M(\K)=q_!
p^*(\K)$. We say that $\K \in \on{D}(\Bun_n)$ is {\em cuspidal} if
$\on{R}^G_M(\F) = 0$ for all proper parabolic subgroups $P$ of $G$.

For a partition $n=n_1+n_2$ let $P(n_1,n_2)$ be the corresponding parabolic
subgroup in $GL_n$ with the Levi factor $GL_{n_1}\times GL_{n_2}$. In
this case diagram \eqref{cus diagram} is
\begin{equation*}
\Bun_n \overset{\p}\leftarrow \Bun_{P(n_1,n_2)}
\overset{\q}\longrightarrow \Bun_{n_1}\times \Bun_{n_2}.
\end{equation*}
We denote the corresponding constant term functor $\on{D}(\Bun_n) \to
\on{D}(\Bun_{n_1}\times \Bun_{n_2})$ by
$\on{R}^n_{n_1,n_2}(\K)$.

It is easy to see that a complex $\K$ is cuspidal if and only if
$\on{R}^n_{n_1,n_2}(\K)=0$ for all
partitions $n=n_1+n_2$, with $n_1,n_2>0$.

\medskip

In this section we prove the following

\begin{thm}   \label{Aut is cusp}
Let $\Aut_E$ be a Hecke eigensheaf on $\Bun_n$ with respect to an
irreducible rank $n$ local system $E$, which satisfies
\conjref{vanishing conjecture}. Then $\Aut_E$ is cuspidal.
\end{thm}

As a corollary we obtain the following statement:

\begin{cor}    \label{finite type}
The perverse sheaf $\Aut_E$ is the extension by zero
from an open substack of finite type on every connected component of
$\Bun_n$.
\end{cor}

\smallskip

\noindent{\em Proof of the corollary} will rely on the following
well-known assertion:

\smallskip

\begin{lem}
For a fixed line bundle $\L$ and an integer $d$, consider the open
substack $U$ of $\Bun^d_n$ which classifies vector bundles $\M$ with
$\Hom(\M,\L)=0$. Then $U$ is of finite type.
\end{lem}

To prove the corollary, let us consider a connected component
$\Bun_n^d$. Without loss of generality we may assume that $d\geq
c_{g,n}$. Take $\L=\L^{\on{est}}$ and let $U$ be as in the lemma.  By
definition, $\Bun_n^d-U$ is contained in $\Bun_n^{\on{vuns}}$.

However, by arguing as in \lemref{cusp support} we obtain that
if a complex $\K$ is cuspidal, then it has zero
stalks at all very unstable bundles.

Therefore, $\K|_{\Bun_n^d}$ is extended by zero from $U$. This
completes the proof of \corref{finite type}.\qed

\smallskip


The proof of \thmref{Aut is cusp} as well as other results of this section
relies on the following computation:

\begin{prop}  \label{average Aut}
Let $\K$ be a Hecke eigensheaf with respect to some rank $n$ local
system $E'$, and let $E$ be another local system, of an arbitrary
rank. Then
$$\Hav^d_{n,E}(\K) \simeq \K\otimes \HH^{\bullet}(X^{(d)},(E\otimes
E'{}^*)^{(d)})\otimes \Ql(\frac{d}{2})[d].$$
\end{prop}

\begin{proof}

Consider the diagram
\begin{equation*}
\Bun_n \xleftarrow{\hl} \Mod^d_n \xrightarrow{\supp\times \hr}
X^{(d)} \times \Bun_n
\end{equation*}

We need to prove that
\begin{equation}   \label{strong av Aut}
(\supp\times\hr)_!(\hl{}^*(\K)\otimes \pi^*(\L^d_{E}))\otimes
(\Ql(\frac{1}{2})[1])^{\otimes d\cdot (n-1)}
\simeq (E\otimes E'{}^*)^{(d)} \boxtimes \K.
\end{equation}

Consider the stack $\wt\Mod{}^d_n$, which classifies the data
$(\M=\M_0\subset \M_1\subset...\subset\M_i=\M')$, where each $\M_j$ is
a rank $n$ vector bundle, and $\M_j/\M_{j-1}$ is a simple skyscraper
sheaf. (Note that there is a canonical isomorphism between
$\wt\Mod{}^d_n$ and the stack $\wt\Mod{}^{-d}_n$, introduced in the
proof of \thmref{thesis}, under which the roles of the projections
$\hl$ and $\hr$ get reversed.)

Let $p: \wt\Mod{}^d_n \to \Mod^d_n$ be the forgetful
map. We also have a natural morphism $\wt\supp:\wt\Mod{}^d_n\to
X^d$. Consider the corresponding diagram
\begin{equation*}
\begin{CD} \Bun_n @<{\wt\hl}<< \wt\Mod{}^d_n @>{\wt\supp\times \wt{\hr}}>>
X^d\times \Bun_n \\ @V{\on{id}}VV @VpVV @V{\sym\times \on{id}}VV \\
\Bun_n @<{\hl}<< \Mod^d_n @>{\supp\times \hr}>> X^{(d)}\times \Bun_n
\end{CD}
\end{equation*}

Since $p$ is small, the complex $p_! \; \wt\supp^*(E^{\boxtimes d})$ is
a perverse sheaf (up to a cohomological shift), which is the
Goresky-MacPherson extension of its own restriction to
$\supp^{-1}(X^{(d)}-\Delta)$. In particular, it carries a canonical
action of the symmetric group $S_d$ and
\begin{equation*}
(p_! \; \wt\supp^*(E^{\boxtimes d}))^{S_d}\simeq \pi^*(\L^d_{E}).
\end{equation*}

As was noted before, the stack $\Mod^1_n$ is isomorphic to
the stack $\H^1_n$, but under this isomorphism the maps $\hl$ and
$\hr$ become interchanged. By iterating the definition of the
Hecke property, we obtain that the fact that $\K$ is a Hecke eigensheaf with
respect to $E'$ then implies that
$$(\wt\supp\times\wt\hr)_! \wt\hl{}^*(\K)\otimes \Ql(\frac{n-1}{2})[n-1]
\simeq (E'{}^*)^{\boxtimes d} \boxtimes \K.$$
Hence we obtain:
\begin{equation*}
(\wt\supp\times\wt\hr)_!(\wt\hl{}^*(\K)\otimes \wt\supp^*(E^{\boxtimes
d}))\otimes (\Ql(\frac{1}{2})[1])^{\otimes d\cdot (n-1)}
\simeq (E\otimes E'{}^*)^{\boxtimes d} \boxtimes \K.
\end{equation*}

By taking the direct image of the last isomorphism under
$\on{sym}:X^d\to X^{(d)}$ we obtain
$$(\supp\times\hr)_!(\hl{}^*(\K)\otimes p_! \; \wt\supp^*(E^{\boxtimes
d}))\otimes (\Ql(\frac{1}{2})[1])^{\otimes d\cdot (n-1)}
\simeq \on{sym}_!((E\otimes E'{}^*)^{\boxtimes d})\boxtimes \K.$$
Moreover, this isomorphism is compatible with the $S_d$--action on both
sides.

By passing to the $S_d$--invariants we obtain formula \eqref{strong av
Aut} and hence the statement of the proposition.
\end{proof}

\ssec{Remark}
Let us see what the isomorphism of formula \eqref{strong av Aut} looks
like at the level of fibers over a given point $D\in X^{(d)}$, in
terms of the general Hecke functors $\He^\la_n$ introduced in
\secref{gen hecke}. To simplify notation we take $D=d\cdot x$, for
some $x\in X$.

By formula \eqref{fiber of Laumon's sheaf} below, the stalk of the LHS
of formula \eqref{strong av Aut} at $D$ can be identified with
\begin{equation*}
\underset{\la \in P^+_{n,d}}\oplus\,
{}_x\He_n^{-w_0(\la)}(\K)\otimes E_x^\lambda.
\end{equation*}
Since $\K$ is a Hecke eigensheaf with respect to $E'$, this is
isomorphic to
$$\underset{\la \in P^+_{n,d}}\oplus\, \K\otimes
(E'{}^*)^\la_x\otimes E^\la_{x} \simeq \K\otimes \on{Sym}^d(E_x\otimes
E'{}^*_x),$$ which is the stalk of the RHS of formula \eqref{strong av
Aut} at $d \times x$.

\ssec{Remark}    \label{minus d}
Recall the stack $\Mod^{-d}_k$ introduced in the proof of
\thmref{thesis}. Denote by $\Hav^{-d}_{k,E}$ the functor
\begin{equation*}
\K \mapsto \hr_!(\hl{}^*(\K)\otimes \pi^*(\L^d_E))\otimes
\Ql(\frac{d\cdot n}{2})[d\cdot n]
\end{equation*}
(where $\hl$ and $\hr$ are taken according to the definition of
$\Mod^{-d}_k$). It follows from the definition that the functor
$\Hav^{-d}_{k,E}$ is both left and right adjoint to $\Hav^d_{k,E^*}$.

In the same way as in the proof of \propref{average Aut} we obtain:
\begin{equation}    \label{Hav -d}
\Hav^{-d}_{n,E}(\K) \simeq \K\otimes
\HH^{\bullet}(X^{(d)},(E\otimes E')^{(d)}) \otimes\Ql(\frac{d}{2})[d].
\end{equation}

Now we are ready to prove \thmref{Aut is cusp}.

\begin{lem}    \label{filtr}
For each $d$, $n=n_1+n_2$, a local system $E$ and $\K\in
\on{D}(\Bun_n)$, the object $\on{R}^n_{n_1,n_2}\circ
\Hav^d_{n,E}(\K)\in \on{D}(\Bun_{n_1}\times \Bun_{n_2})$ has a
canonical filtration by the objects $(\Hav^{d_1}_{n_1,E}\times
\Hav^{d_2}_{n_2,E})\circ \on{R}^n_{n_1,n_2}(\K)\otimes
\Ql(\frac{-n_1\cdot d_2}{2})[-n_1\cdot d_2]$ for all possible
partitions $d=d_1+d_2$ with $d_1,d_2\geq 0$.
\end{lem}

\ssec{Proof of \thmref{Aut is cusp}}

\thmref{Aut is cusp} follows from \lemref{filtr}. Indeed, take $d >
2n^2(2g-2)$.  On the one hand, according to \propref{average Aut},
$$\Hav^d_{n,E}(\Aut_E)\simeq \Aut_E\otimes
\HH^{\bullet}(X^{(d)},(E\otimes E^*)^{(d)}),$$ hence
$$\on{R}^n_{n_1,n_2}\circ \Hav^d_{n,E}(\Aut_E)\simeq
\on{R}^n_{n_1,n_2}(\Aut_E)\otimes \HH^{\bullet}(X^{(d)},(E\otimes
E^*)^{(d)}) \otimes \Ql(\frac{d}{2})[d].$$

On the other hand, \conjref{vanishing conjecture} implies that all
$(\Hav^{d_1}_{n_1,E}\times \Hav^{d_2}_{n_2,E})\circ
\on{R}^n_{n_1,n_2}(\Aut_E)$ must vanish, because either $d_1$ or $d_2$
must be greater than $n^2(2g-2)$.

However, since $E\otimes E^*$ contains the trivial rank one local
system, $\HH^{\bullet}(X^{(d)},(E\otimes E^*)^{(d)})\neq 0$ for any
$d$. Hence, $\on{R}^n_{n_1,n_2}(\Aut_E)=0$. This completes the proof
of \thmref{Aut is cusp}.

\ssec{Proof of \lemref{filtr}}

Consider the Cartesian product $\Mod^d_n\underset{\Bun_n}\times
\Bun_{P(n_1,n_2)}$, where we used the projection $\hr:\Mod^d_n\to
\Bun_n$ to form the Cartesian product. Our task is to calculate the
direct image under $$\Mod^d_n\underset{\Bun_n}\times
\Bun_{P(n_1,n_2)}\overset{\hr\times\on{id}}\longrightarrow
\Bun_{P(n_1,n_2)}\overset{\q}\to {\Bun_{n_1}\times \Bun_{n_2}}$$ of
the pull-back under $\Mod^d_n\underset{\Bun_n}\times
\Bun_{P(n_1,n_2)}\overset{\on{id}\times\p}\longrightarrow \Mod_n$ of
the complex $\hl{}^*(\K)\otimes\pi^*(\L^d_E)$.

By definition, the above Cartesian product classifies the data of

\begin{equation}  \label{restr cart product}
\begin{split}
&\M\in \Bun_n,\,\, \M'\in\Bun_n,\,\, \beta:\M'\hookrightarrow \M, \\
&\M_1\in \Bun_{n_1},\,\, \M_2\in \Bun_{n_2},\,\, 0 \to \M_1 \to \M \to
\M_2 \to 0
\end{split}
\end{equation}

First, we decompose $\Mod^d_n\underset{\Bun_n}\times
\Bun_{P(n_1,n_2)}$ into locally closed substacks, which we will denote
by $(\Mod^d_n\underset{\Bun_n}\times \Bun_{P(n_1,n_2)})^{d_1,d_2}$ as
follows:

A point of $\Mod^d_n\underset{\Bun_n}\times \Bun_{P(n_1,n_2)}$ as in
\eqref{restr cart product} belongs to
$(\Mod^d_n\underset{\Bun_n}\times \Bun_{P(n_1,n_2)})^{d_1,d_2}$ if
$\deg(\M'\cap \M_1)=\deg(\M_1)-d_1$.

{}From each $(\Mod^d_n\underset{\Bun_n}\times
\Bun_{P(n_1,n_2)})^{d_1,d_2}$ there is a natural map to
$\Mod^{d_1}_{n_1}\times \Mod^{d_2}_{n_2}$, which sends a point as
above to $(\M'_1:=\M'\cap \M_1\hookrightarrow \M_1,
\M'_2:=\M'/\M'_1\hookrightarrow \M_2)$.

To prove the proposition it suffices to show that the direct image
under this map of the complex that we obtain on
$(\Mod^d_n\underset{\Bun_n}\times \Bun_{P(n_1,n_2)})^{d_1,d_2}$ by
restriction from $\Mod^d_n\underset{\Bun_n}\times \Bun_{P(n_1,n_2)}$
can be canonically identified with $(\hl\times
\hl)^*(\on{R}^n_{n_1,n_2}(\K))\otimes (\pi\times
\pi)^*(\L^{d_1}_{E}\boxtimes \L^{d_2}_{E})\otimes \Ql(\frac{-n_1\cdot
d_2}{2})[-n_1\cdot d_2]$ in the diagram
$$\Bun_{n_1}\times \Bun_{n_2} \overset{\hl\times\hl}\longleftarrow
\Mod^{d_1}_{n_1}\times \Mod^{d_2}_{n_2}\overset{\pi\times
\pi}\longrightarrow \Coh_0^{d_1}\times \Coh_0^{d_2}.$$

For that purpose, we decompose the map
$$(\Mod^d_n\underset{\Bun_n}\times \Bun_{P(n_1,n_2)})^{d_1,d_2}\to
\Mod^{d_1}_{n_1}\times \Mod^{d_2}_{n_2}$$ as a composition of several
ones. First, we introduce the stack $\Y_1$, which classifies the data
of
\begin{equation*}
\begin{split}
& 0\to \M'_1\to \M_1\to \T_1\to 0,\,\,\M_1\in\Bun_{n_1},\,\, \T_1\in
\Coh_0^{d_1}, \\ & 0\to \M'_2\to \M_2\to \T_2\to
0,\,\,\M_2\in\Bun_{n_2},\,\, \T_1\in \Coh_0^{d_2}, \\ & 0\to \T_1\to
\T\to \T_2\to 0,\,\, 0\to \M'_1\to \M'\to \M'_2\to 0.
\end{split}
\end{equation*}

It is easy to see that the natural map
$(\Mod^d_n\underset{\Bun_n}\times \Bun_{P(n_1,n_2)})^{d_1,d_2}\to
\Y_1$ is a fibration into affine spaces, with each fiber being a
principal homogeneous space for $\on{Ext}^1(\T^2,\M'_1)$. Therefore,
by the projection formula, the direct image of our complex to $\Y_1$
is the pull-back under the map $\Y_1\to \Bun_n\times \Coh_0^d$ (which
sends a point as above to $(\M',\T)$) of $\K\boxtimes \L^d_E\otimes
\Ql(\frac{-n_1\cdot d_2}{2})[-n_1\cdot d_2]$.

\medskip

Now, let $\Y_2$ be the stack classifying the data of
\begin{equation*}
\begin{split}
& 0\to \M'_1\to \M_1\to \T_1\to 0,\,\,\M_1\in\Bun_{n_1},\,\, \T_1\in
\Coh_0^{d_1}, \\ & 0\to \M'_2\to \M_2\to \T_2\to
0,\,\,\M_2\in\Bun_{n_2},\,\, \T_2\in \Coh_0^{d_2}, \\ & 0\to \M'_1\to
\M'\to \M'_2\to 0.
\end{split}
\end{equation*}
The projection $\Y_1\to \Y_2$ corresponds to forgetting the class of
the extension $0\to \T_1\to \T\to \T_2\to 0$. Moreover, we have a
Cartesian square:
$$
\begin{CD} \Y_1 @>>> \H\L_0^{d^1} \\ @VVV @VVV \\ \Y_2 @>>>
\Coh_0^{d_1}\times \Coh_0^{d_2}.  \end{CD}
$$
Using the projection formula and the fact that
$\HL_0^{d_1}(\L^d_E)\simeq \L^{d_1}_E\boxtimes \L^{d_2}_E$, we obtain
that direct image under $\Y_1\to\Y_2$ of the pull-back of $\K\boxtimes
\L^d_E$ is the tensor product of the pull-back of $\K$ under the map
$\Y_2\to \Bun_n$, which sends a point as above to $\M'$ and the
pull-back of $\L^{d_1}_E\boxtimes \L^{d_2}_E$ under the natural map
from $\Y_2$ to $\Coh_0^{d_1}\times \Coh_0^{d_2}$.

Finally, note that we have a Cartesian square:
$$
\begin{CD} \Bun_{P(n_1,n_2)} @<<< \Y_2\\ @V{\q}VV @VVV \\ \Bun_{n_1}\times
\Bun_{n_2} @<{\hl\times \hl}<< \Mod^{d_1}_{n_1}\times \Mod^{d_2}_{n_2}
, \end{CD}
$$
where the upper horizontal arrow sends a point of $\Y_2$ as above to
$0\to \M'_1\to \M'\to \M'_2\to 0$.  The assertion follows now by the
projection formula.

\section{Proof of the Vanishing Conjecture over $\Fq$} \label{proof
of vanishing}

In this section we prove \conjref{vanishing conjecture} in the case
when the ground field $\kk$ is a finite field $\Fq$, i.e., that the
functor $\Hav^d_{k,E}:\on{D}(\Bun_k)\to \on{D}(\Bun_k)$ introduced in
\secref{averaging} is identically zero if $E$ is an irreducible local
system of rank $n$, and $k$ and $d$ satisfy the inequalities $k<n$ and
$d > kn(2g-2)$.

Namely, we will prove the following proposition:

\begin{prop}    \label{let E}
Let $E$ be a rank $n$ local system on $X$ over the finite field $\Fq$,
which is

\smallskip

{\em (a)} pure up to a twist by a one-dimensional representation of
the Weil group of $\Fq$,

\smallskip

\noindent and satisfies one of the following conditions:
\begin{enumerate}
\item[(b)] there exists a cuspidal Hecke eigenfunction associated to
the pull-back of $E$ to $X \underset{\Fq}\times {\mathbb F}_{q_1}$ for
any finite extension ${\mathbb F}_{q_1}$ of $\Fq$; or

\item[(b')] the space of unramified cuspidal automorphic functions on
the group $GL_k$ over the ad\`eles is spanned by the Hecke
eigenfunctions corresponding to rank $k$ local systems on $X
\underset{\Fq}\times {\mathbb F}_{q_1}$, for all $k<n$.
\end{enumerate}
Then the Vanishing \conjref{vanishing conjecture} holds for $E$.
\end{prop}

According to \cite{Lf}, Theorem VII.6, any irreducible local system
$E$, such that $\on{det} E$ is of finite order, is pure. Therefore
condition (a) of \propref{let E} is satisfied for any irreducible rank
$n$ local system $E$. Moreover, both statements (b) and (b') hold for
such $E$, according to the main theorem of Lafforgue's work. Hence we
obtain that the Vanishing \conjref{vanishing conjecture} is true for
all irreducible local systems if the ground field $\kk$ is finite.

\medskip

Our proof of \propref{let E} proceeds as follows. We first show
vanishing of $\Hav^d_{k,E}$ at the level of functions.  Using the
purity property conjectured by Deligne and proved by Lafforgue
\cite{Lf}, we will then deduce that $\Hav^d_{k,E}(\K)=0$ for any $\K
\in \on{D}(\Bun_k), k=1,\ldots,n-1, d>kn(2g-2)$.

\ssec{$L$--functions}

Let $\Gamma^{(1)} = (\gamma^{(1)}_x)_{x\in |X|}$ and $\Gamma^{(2)} =
(\gamma^{(2)}_x)_{x\in |X|}$ be two collections of semi-simple
conjugacy classes in $GL_k(\Ql)$ and $GL_n(\Ql)$, respectively.
We attach to it the $L$--function
$$
L(\Gamma^{(1)},\Gamma^{(2)},t) = \prod_{x \in
|X|}\on{det}(\on{Id}_{kn} - (\gamma^{(1)}_x \otimes \gamma^{(2)}_x) \;
t^{\deg x})^{-1},
$$
viewed as a formal power series in $t$.

To an unramified irreducible representation $\pi = \otimes'_{x \in X}
\pi_x$ of $GL_k(\AA)$, where $\AA$ is the ring of ad\`eles of $F=\Fq(X)$,
we attach the collection $\Gamma_\pi = (s_x)_{x\in |X|}$, where $s_x$
is the Satake parameter of $\pi_x$.

If $\pi$ and $\pi'$ are unramified irreducible representations $\pi$
of $GL_k(\AA)$ and $GL_n(\AA)$, respectively, we write
$$
L(\pi \times \pi',t) := L(\Gamma_\pi,\Gamma_{\pi'},t).
$$
If $\pi$ and $\pi'$ are in addition cuspidal automorphic
representations, then $L(\pi \times \pi',t)$ is the Rankin--Selberg
$L$--function of the pair $\pi,\pi'$. The following statement follows
from results of \cite{PS1,CPS} (see \cite{Lf}, Appendice B, for a
review).\footnote{We are grateful to V.~Drinfeld for pointing it out
to us}

\begin{thm}    \label{rs}
If $\pi,\pi'$ are cuspidal automorphic representations and $k<n$, then
$L(\pi \times \pi',t)$ is a polynomial of degree $kn(2g-2)$.
\end{thm}

Next, we attach to a rank $n$ local system $E$ on $X$ the collection
of conjugacy classes $\Gamma_E = (\on{Fr}_x|_{E_x})_{x\in |X|}$.

If $E$ and $E'$ are two local systems on $X$, of ranks $n$ and $k$,
respectively, we write:
$$
L(E' \times E,t) := L(\Gamma_{E'},\Gamma_E,t).
$$

\begin{lem}    \label{E tensor E}
If both $E$ and $E'$ are irreducible and $k<n$, then $L(E' \times
E,t)$ is a polynomial of degree $kn(2g-2)$.
\end{lem}

\begin{proof}
Using the definition of $L(E' \times E,t)$ and the
Grothendieck-Lefschetz formula, we obtain:
\begin{equation*}
L(E' \times E,t) = \sum_{d\geq 0}
\on{Tr}(\on{Fr},\HH^d(X^{(d)},(E'\otimes E)^{(d)})) t^d.
\end{equation*}
Since $E' \otimes E$ is irreducible by our assumptions,
$\HH^0(E'\otimes E) = H^2(E'\otimes E) = 0$. Therefore,
$$\HH^d(X^{(d)},(E'\otimes E)^{(d)}) \simeq
\Lambda^d(\HH^1(X,E'\otimes E)) = 0$$ for all $d>\dim
\HH^1(X^{(d)},E'\otimes E) = kn(2g-2)$.
\end{proof}

If $E$ is a rank $n$ local system on $X$, and $\pi$ is an irreducible
unramified representation of $GL_k(\AA)$, we will write
$$
L(\pi \times E,t) := L(\Gamma_\pi,\Gamma_E,t).
$$

\ssec{Computation of $\Hav^d_{k,E}$ at the level of functions}
\label{comp for functions}

The functor $\Hav^d_{k,E}$ gives rise to a linear map on the space of
functions on the set $\Bun_k(\Fq)$ of $\Fq$--points of $\Bun_k$. We
denote this operator by $\Havv^d_{k,E}$. In this subsection we will
prove that $\Havv^d_{k,E} \equiv 0$ for all $k=1,\ldots,n-1$ and
$d>kn(2g-2)$.

Recall that $\Bun_k(\Fq)$ is naturally identified with the double
quotient $$GL_k(F) \backslash GL_n(\AA)/GL_n(\OO)$$ (see, e.g.,
\cite{FGKV}, Sect.~2). Let $\pi$ be a cuspidal unramified automorphic
representation of $GL_k(\AA)$. Attached to it is a cuspidal
automorphic function on $\Bun_k(\Fq)$, unique up to a non-zero scalar
multiple. We normalize it in some way and denote the result by
$f_\pi$.

In \remref{minus d} we defined the functor $\Hav^{-d}_{k,E}$, which is
left and right adjoint to $\Hav^d_{k,E^*}$. Denote by
$\Havv^{-d}_{k,E}$ the corresponding linear map on the space of
functions on $\Bun_k(\Fq)$. We have the following analogue of formula
\eqref{Hav -d}, which is proved using a calculation similar to the one
presented in the proof of \propref{average Aut}:
\begin{equation}    \label{appearance of L}
\sum_{d \geq 0} \Havv^{-d}_{k,E}(f_\pi) \cdot t^d = L(\pi \times
E,t)\cdot f_\pi.
\end{equation}

It is clear from the definition that
\begin{equation}    \label{dual}
\langle \Havv^d_{k,E}(f),f' \rangle = \langle
f,\Havv^{-d}_{k,E}(f') \rangle,
\end{equation}
where the inner product of two automorphic functions $f_1, f_2$ on
$GL_k(\AA)$ is defined by the formula
$$
\langle f_1,f_2 \rangle = \int_{GL_k(F)\backslash GL_k(\AA)} f_1(g)
f_2(g) dg.
$$
Formula \eqref{appearance of L} then implies:
\begin{equation}    \label{scalar product}
\langle \Havv^d_{k,E}(f),f_\pi \rangle = L(\pi \times E,t)
\langle f,f_\pi \rangle.
\end{equation}

\begin{lem}    \label{poly}
Let $E$ be a rank $n$ local system such that $L(\pi \times E,t)$ is
a polynomial of degree $kn(2g-2)$ for all irreducible cuspidal
automorphic representations $\pi$ of $GL_k(\AA), k=1,\ldots,n-1$. Then
$\Havv^d_{k,E}(f) = 0$ for any function $f$ on $\Bun_k(\Fq)$.
\end{lem}

\begin{proof}
By induction, we may assume that the assertion is known for $k'<k$. In
the same way as in the proof \thmref{Aut is cusp} we then obtain that
$\Havv^d_{k,E}(f)$ is a cuspidal function for any function $f$ on
$\Bun_k(\Fq)$. By formula \eqref{scalar product}, if $L(\pi \times
E^*,t)$ is a polynomial of degree $kn(2g-2)$, then $\langle
\Havv^d_{k,E}(f),f' \rangle = 0$ for any function $f$, any cuspidal
automorphic function $f'$ and all $d>kn(2g-2)$. Therefore
$\Havv^d_{k,E}(f) = 0$.
\end{proof}

Thus, in order to prove vanishing of $\Havv^d_{k,E}$ for
$k=1,\ldots,n-1$ and $d>kn(2g-2)$, we need to show that $L(\pi \times
E,t)$ is a polynomial of degree $kn(2g-2)$ for all irreducible
cuspidal automorphic representations $\pi$ of $GL_k(\AA)$ for
$k=1,\ldots,n-1$. This can be done in two ways: by identifying $L(\pi
\times E,t)$ with $L(\pi \times \pi',t)$ or with $L(E' \times
E,t)$. As the result, we obtain that $\Havv^d_{k,E} \equiv 0,
k=1,\ldots,n-1$, if either of the statements (b) or (b') listed in
\propref{let E} is true.

Indeed, if the statement (b) is true, then there exists an unramified
cuspidal automorphic representation $\pi'$ of $GL_n(\AA)$, such that
$L(\pi \times E,t) = L(\pi \times \pi',t)$. Vanishing of
$\Havv^d_{k,E}$ for $k<n$ and $d>kn(2g-2)$ then follows from
\thmref{rs}.

If the statement (b') is true, then vanishing follows from \lemref{E
tensor E}.

\ssec{Conclusion of the proof of \propref{let E}}

To complete the proof of \propref{let E} we need to show that
vanishing of the operator $\Havv^d_{k,E}$ at the level of functions
implies the vanishing of the operator $\Hav^d_{k,E}$ at the level of
sheaves, provided that $E$ is pure. In order to do that, we proceed as
follows: for each $x: \on{Spec}\Fq \to \Bun_k$, denote by $\delta_x$
the direct image with compact support of the constant sheaf on
$\on{Spec} \Fq$.  Proving that the functor $\Hav^d_{k,E}$ vanishes is
equivalent to showing that $\Hav^d_{k,E}(\delta_x)=0$, for all $x$.

Since $\Bun_k$ is a stack (and not a scheme), $\delta_x$ is not
necessarily an irreducible perverse sheaf, but it is a mixed
complex. Therefore, it suffices to show that $\Hav^d_{k,E}(\K)=0$, for
any mixed complex $\K$. Decomposing $\Hav^d_{k,E}(\K)$ in the derived
category, we obtain that it is enough to show that
$\Hav^d_{k,E}(\K)=0$, when $\K$ is a pure perverse sheaf.

Now let $E$ be a pure irreducible rank $n$ local system on $X$. Then
Laumon's sheaf $\L_E$ is also pure. The pull-back with respect to a
smooth morphism preserves purity, and so does the push-forward with
respect to a proper representable morphism (see \cite{BBD}). But the
morphism $(\hl\times \pi):\Mod^d_n\to \Bun_n\times \Coh_0^d$ is
smooth, and the morphism $\hr: \Mod^d_n \to \Bun_n$ is proper and
representable. Hence $\Hav^d_{k,E}(\K)$ is pure, if $\K$ is pure.

The function $\text{\tt f}_{q_1}(\Hav^d_{k,E}(\K))$ associated to the
sheaf $\Hav^d_{k,E}(\K)$ equals $\Havv^d_{k,E}(\text{\tt
f}_{q_1}(\K))$ for any for $q_1=q^r, r \in \ZZ_{>0}$ (here we use the
notation introduced in \secref{conventions}). But according to the
computation of \secref{comp for functions}, $\Havv^d_{k,E}(\text{\tt
f}_{q_1}(\K)) = 0$ for all $k=1,\ldots,n-1$ and $d>kn(2g-2)$ if either
of the conditions (b) or (b') of \propref{let E} holds for $E$. In
addition, we have:

\begin{lem}
A pure complex $\F$ vanishes if and only if the corresponding function
$\text{\tt f}_{q_1}(\F)$ is zero, for all $q_1=q^r, r\in \ZZ_{>0}$.
\end{lem}

\begin{proof}

Since $\F$ is non-zero, there exists a locally closed subset $U$ such
that $\F|_U$ is locally constant and non-zero. Since $\F$ is pure,
$\F|_U$ is pointwise pure. But for a pointwise pure non-zero locally
constant complex, all functions $\text{\tt f}_{q_1}(\F|_U)$ cannot be
identically equal to zero for all $q_1=q^r, r\in \ZZ_{>0}$, by the
condition on the absolute values of the Frobenius eigenvalues.

\end{proof}

Therefore the statement of \conjref{vanishing conjecture} holds for
$E$. This completes the proof of \propref{let E}.

\ssec{Remark}
Formula \eqref{dual} has a geometric
counterpart:
$$\langle \Hav^d_{k,E}(\K),\K' \rangle \simeq \langle
\K,\Hav^{-d}_{k,E^*}(\K') \rangle, \qquad \K,\K' \in \on{D}(\Bun_k),
$$
where $\langle \K,\K' \rangle := \on{RHom}(\K,\K')$.  Note that a
priori $\on{RHom}(\K,\K')$ makes sense if $\K$ is the $!$--extension
from a substack of $\Bun_k$ of finite type, or $\K'$ is the
$*$--extension from a substack of $\Bun_k$, whose intersection with every
connected component is of finite type.

Let $E$ and $E'$ be two irreducible local systems on $X$, of ranks $n$
and $k$, respectively, where $k<n$.  Let us assume that the vanishing
\conjref{vanishing conjecture} holds for $E'$. In particular,
$\Aut_{E'}$ exists and is cuspidal, according to \thmref{Aut is cusp}.
Therefore, for every $d$, $\Aut_{E'}|_{\Bun_k^d}$ is extended by zero
from an open substack of finite type of $\Bun^d_n$.

\medskip

{}From formula \eqref{Hav -d} we obtain the following analogue of
formula \eqref{scalar product}:
\begin{equation}    \label{scalar product1}
\langle \Hav^d_{k,E}(\K),\Aut_{E'} \rangle \simeq \HH^d(X^{(d)},(E^*
\otimes E')^{(d)}) \otimes \langle \K,\Aut_{E'} \rangle\otimes
\Ql(\frac{d}{2})[d].
\end{equation}
Since $\HH^d(X^{(d)},(E^*\otimes E')^{(d)})=0$ for $d>kn(2g-2)$ (see
the proof of \lemref{E tensor E}), we find that
$$
\langle \Hav^d_{k,E}(\K),\Aut_{E'} \rangle = 0,
$$
for all $\K \in \on{D}(\Bun_k)$, if $d>kn(2g-2)$. Thus we obtain a
geometric analogue of \propref{let E}:

\begin{prop}
Suppose that for $k=1,\ldots,n-1$ the vanishing \conjref{vanishing conjecture}
is true for rank $k$ local systems on
$X$ and in addition the following statement holds:

$\on{(b'')}$ if $\F \in \on{D}(\Bun_k)$ is cuspidal and satisfies
$\langle \F,\Aut_{E'} \rangle = 0$ for all irreducible rank $k$ local
system $E'$ on $X$, then $\F=0$.

Then the Vanishing \conjref{vanishing conjecture} is true for any
irreducible local system on $X$ of rank $n$.
\end{prop}

The above statement $\on{(b'')}$ is known to be true for $k=1$ in the
case when $\on{char}\kk=0$, by the Fourier--Mukai transform
\cite{La:mukai,Roth}.

\appendix

\section{Hecke functors and Whittaker sheaves}

\ssec{General Hecke functors}    \label{gen hecke}

We recall some results from Sect. 5 of \cite{FGV}.

Let $_x\H_n$ be the full Hecke correspondence stack at $x \in |X|$. In
other words, $_x\H_n$ classifies triples $(\M,\M',\beta)$, where $\M$
and $\M'$ are rank $n$ bundles on $X$ and $\beta$ is an isomorphism
$\M|_{X-x} \overset{\sim}\to \M'|_{X-x}$. To a dominant weight
$\lambda$ of $GL_n(\Ql)$ we associate a closed finite-dimensional
substack $_x\Hb_n^\lambda$ of $_x\H_n$, which classifies the triples
$(\M,\M',\beta)$, such that for an algebraic representation $V$ of
$GL_n(\kk)$, whose weights are $\leq \check\nu$, we have the following
embeddings induced by $\beta$ on the entire $X$:
$$
V_{\M'}(\langle w_0(\la),\check\nu \rangle \cdot x) \subset V_{\M}
\subset V_{\M'}(\langle \la,\check\nu \rangle \cdot x),
$$
where $V_{\M}$ is the vector bundle on $X$ associated with $V$ and the
principal $GL_n$-bundle on $X$ corresponding to $\M$ (recall that
$w_0$ stands for the permutation $(d_1,d_2,\ldots,d_n) \mapsto
(d_n,\ldots,d_2,d_1)$).

Using this stack, we define the Hecke functor
$_x\He_n^\lambda:\on{D}(\Bun_n)\to \on{D}(\Bun_n)$ by the formula
\begin{equation*}
\K\mapsto \hr_!(\hl{}^*(\K)\otimes \IC_\lambda)
\otimes\Ql(\frac{\dim(\Bun_n)}{2})[\dim(\Bun_n)],
\end{equation*}
where $\hl$ (resp., $\hr$) sends $(\M,\M',\beta)$ to $\M$ (resp.,
$\M'$), and $\IC_\la$ is the intersection cohomology sheaf on
$_x\Hb_n^\lambda$.

In particular, if $\lambda$ is the $i$-th fundamental weight
$\omega_i$, then the stack ${}_x\Hb_n^{\omega_i}$ is nothing but the
preimage of $x \in X$ in $\H^i_n$ under $\supp:\H^i_n\to X$. Hence
${}_x\He_n^{\omega_i}$ is the composition of $\He^i_n$ followed by the
restriction to $x\times \Bun_n\simeq \Bun_n \subset X\times\Bun_n$.

The results of \cite{Lu2,Gi,MV} imply the following formula:
\begin{equation}    \label{tensor product}
_x\He_n^\la \circ \, {}_x\He_n^\mu
\simeq \underset{\nu \in P^+_n}\oplus {}_x\He_n^\nu \otimes
\on{Hom}_{GL_n}(V^\nu,V^\la \otimes V^\mu),
\end{equation}
where the notation $V^\lambda$ is as in \secref{explicit description}.

Consider the fiber product
$Z^{\lambda,x}:=\ol\Q^0\underset{\Bun^0_n}\times
{}_x\Hb_n^\lambda$. It was proved in \cite{FGV}, Prop.~5.3.4, that
there exists a commutative diagram
\begin{equation*}
\begin{CD}
\ol\Q^0   @<{'\hl}<<  Z^{\lambda,x}  @>{'\hr}>> \ol\Q^{-\lambda,x} \\
@V{q}VV     @V{'q}VV    @V{q}VV \\
\Bun_n  @<{\hl}<< _x\Hb_n^\lambda @>{\hr}>> \Bun_n
\end{CD}
\end{equation*}
where we write $\ol\Q^{-\lambda,x}$ for $\ol\Q^{-\ol\lambda,\ol x}$
when the collection $\ol\lambda$ (resp., $\ol x$) consists of just one
element $\lambda$ (resp., $x$), in the notation of \secref{explicit
description}. According to Theorems 3 and 4 of \cite{FGV}, adapted to
our present notation, we have:
\begin{equation}  \label{Whitaction}
'\hr_!({}'\hl{}^*(\Psi^0)\otimes
{}'q^*(\IC_\lambda))\otimes\Ql(\frac{\dim(\Bun_n)}{2})[\dim(\Bun_n)]
\simeq \Psi^{-\lambda,x}.
\end{equation}

More generally, let $x^1,\ldots,x^m$ be a set of distinct points,
different from $x$, and $\ol\mu = (\mu^0,\mu^1,\ldots,\mu^m)$ be a
collection of dominant weights. Set $Z^{\lambda,\ol\mu,\ol
x}:=\ol\Q^{-\ol\mu,\ol x} \underset{\Bun_n}\times {}_x\Hb_n^\lambda$,
where $\ol x = (x,x^1,\ldots,x^m)$. Denote $\ol\la =
(\la,0,\ldots,0)$. We have a commutative diagram
\begin{equation*}
\begin{CD}
\ol\Q^{-\ol\mu,\ol x} @<{'\hl}<< Z^{\lambda,\ol\mu,\ol x} @>{'\hr}>>
\ol\Q^{-\ol\mu-\ol\la,\ol x} \\ @V{q}VV @V{'q}VV @V{q}VV \\ \Bun_n
@<{\hl}<< _x\Hb_n^\lambda @>{\hr}>> \Bun_n
\end{CD}
\end{equation*}
Denote by $_x'\hspace*{-1mm}\He_n^\la$ the functor
$\on{D}(\ol\Q^{-\ol\mu,\ol x}) \to
\on{D}(\ol\Q^{-\ol\mu-\ol\la,\ol x})$,
$$
\K \mapsto {}'\hr_!({}'\hl{}^*(\K)\otimes
{}'q^*(\IC_\lambda))\otimes\Ql(\frac{\dim(\Bun_n)}{2})[\dim(\Bun_n)].
$$
Then Corollary 5.4.3 of \cite{FGV} gives:
\begin{equation}  \label{Whitaction1}
_x'\hspace*{-1mm}\He_n^\la(\Psi^{-\ol\mu,x}) \simeq \underset{\nu \in
P^+_n}\oplus
\Psi^{-(\nu,\mu^1,\ldots,\mu^m),\ol x} \otimes
\on{Hom}_{GL_n}(V^\nu,V^\la \otimes V^{\mu^0}).
\end{equation}

\ssec{Proof of \propref{cass-shal}}    \label{proof cass-shal}

To simplify the notation, we consider the case when $m=1$, i.e.,
$\ol\mu=\mu$ and $\ol{x}=x$. For $m>1$ the proof is essentially the
same. By construction, $\Q^{\mu,x}$ and $\ol\Q^{\mu,x}$ are substacks
of $\ol\Q^{d\cdot x} \subset \ol\Q^d$.

Observe that the preimage of $d\cdot x \in X^{(d)}$ under $\supp:
\Mod_n^d\to X^{(d)}$ can be identified with the closed substack of
$_x\H_n$, which classifies those triples $(\M,\M',\beta)$, for which
$\beta: \M|_{X-x} \overset{\sim}\to \M'|_{X-x}$ extends to an
embedding $\M \hookrightarrow \M'$ over the entire $X$, and
$\deg(\M')-\deg(\M)=d$.

Recall the morphism $\pi: \on{Mod}_n^d\to \Coh^d_0$. Let us denote
the pull-back $\pi^*(\L_E^d)\otimes \Ql(\frac{d\cdot n}{2})[d\cdot n]$
by $\P^d_E$.

It follows from the
results of \cite{La1}, Sect. 3 and \cite{FGKV}, Sect. 4.2 that the
$*$--restriction of $\P^d_E$
to $\supp^{-1}(d \cdot x) \subset
\on{Mod}_n^d$ can be canonically identified with
\begin{equation}  \label{fiber of Laumon's sheaf}
\underset{\lambda \in P^+_{n,d}}\oplus \,\IC_{-w_0(\laa)}\otimes
\Hom_{GL_n}(V^\lambda,\on{Sym}^d(V\otimes E_x)) \otimes
\Ql(\frac{d}{2})[d].
\end{equation}
Here
\begin{equation}    \label{P+nd}
P^+_{n,d} = \{(d_1,\ldots,d_n)|d_1 \geq d_2 \geq \ldots \geq
d_n \geq 0, \underset{i=1}{\overset{n}\Sigma} d_i = d \},
\end{equation}
and $V = V^{(1,0,\ldots,0)}$ stands for the defining representation of
$GL_n$. Further, we have:
\begin{equation} \label{binom}
\on{Sym}^d(V \otimes E_x)\simeq \underset{\lambda \in
P^+_{n,d}}\oplus V^\lambda\otimes E_x^\lambda.
\end{equation}

By comparing the definition of ${\mathcal W}^d_E$ with formulas
\eqref{Whitaction}, \eqref{fiber of Laumon's sheaf} and \eqref{binom},
we obtain that the restriction of ${\mathcal W}^d_E$ to $\ol\Q^{d\cdot
x}$ can be identified with
\begin{equation}    \label{identify W}
\underset{\la \in P^+_{n,d}}\oplus
\Psi^{w_0(\la),x}\otimes E_{w_0(\la),x}\otimes \Ql(\frac{d}{2})[d],
\end{equation}
which is what we had to prove.

\ssec{Proof of \propref{Hecke on Q}: local computation}

Consider the morphism $\tau: {}\Qpp^d\to X\times X^{(d+1)}$, sending
$(x,\M,(s_i))$ to $(x,D)$, where $D$ is the divisor of zeroes of the
map $s_n: \Omega^{n(n-1)/2} \to (\det \M)(x)$. We start by describing
explicitly the restriction of the complex $\Wp^d$ to $\Qp^d\cap
\tau^{-1}(x\times D)$.

Let us write $D=d^0\cdot x+d^1\cdot x^1+...+d^m\cdot x^m$, where $x^i$
are pairwise distinct and different from $x$, and
$d^0+d^1+...+d^m=d+1$.

It follows from the definitions given in \secref{explicit description}
that the stack $\Qpp^d \cap \tau^{-1}(x\times D)$ is identified with
$\ol\Q^{\ol\nu'',\ol x}$, where $\ol x=(x,x^1,...,x^m)$ and
$\ol\nu'' = (\nu^0{}'',\nu^1,...,\nu^m)$, with $\nu^0{}'' =
(-1,0,\ldots,0,d^0)$ and $\nu^j = (0,\ldots,0,d^j),
j=1,\ldots,m$. Furthermore, $\Qp^d \cap \tau^{-1}(x\times D)$ is
identified with the substack $\ol\Q^{\ol\nu',\ol x}$ of
$\ol\Q^{\ol\nu,\ol x}$, where $\ol\nu' = (\nu^0{}',\nu^1,...,\nu^m)$,
with $\nu^0{}' = (0,-1,0,\ldots,0,d^0)$ and $\nu^j, j=1,\ldots,m$, as
above.

We also have a morphism $\tau: \ol\Q^{d+1} \to X^{(d+1)}$, and we
identify $\tau^{-1}(D)$ with $\ol\Q^{\ol\nu,\ol x}$, where $\ol\nu =
(\nu^0,\nu^1,...,\nu^m)$, with $\nu^0 = (0,\ldots,0,d^0)$ and
$\nu^j, j=1,\ldots,m$, as above.

We have a commutative diagram
\begin{equation*}    \label{diag for Q1}
\begin{CD}
\ol\Q^{\ol\nu,\ol x} @<{'\hl}<< Z^{\omega_1,\mu,\ol x} @>{'\hr}>>
\ol\Q^{\ol\nu'',\ol x} \\
@VVV @VVV @VVV \\
\ol\Q^{d+1} @<{\thl}<< \ol\Q^{d+1} \underset{\Bun_n}\times \H^1_n
@>{\thr}>> \Qpp^d
\end{CD}
\end{equation*}
Moreover, it follows from the definitions that both squares of this
diagram are Cartesian. Therefore we obtain the following formula for
the restriction of the sheaf $\Wp^d = \thr_! \,
\thl{}^*(\W_E^{d+1})|_{\Qp^d}\otimes \Ql(\frac{n-2}{2})[n-2]$
to $\Qp^d \cap \tau^{-1}(x \times D) =
\ol\Q^{\ol\nu',\ol x}$:
\begin{equation}    \label{cart cor}
\Wp^d|_{\ol\Q^{\ol\nu',\ol x}} \simeq {}_x'\hspace*{-1mm}\He_n^{\omega_1}
\left(\W^{d+1}_E|_{\ol\Q^{\ol\nu,\ol x}} \right)|_{\ol\Q^{\ol\nu',\ol
x}}.
\end{equation}
By \propref{cass-shal},
\begin{equation}    \label{restr of WE}
\W^{d+1}_E|_{\ol\Q^{\ol\nu,\ol x}} \simeq \underset{\ol\mu}\oplus
\Psi^{\ol\mu,\ol x} \otimes E_{\mu^0,x} \otimes E_{\mu^1,x^1} \otimes
... \otimes E_{\mu^m,x^m},
\end{equation}
where the summation is over all $\ol\mu=(\mu^0,...,\mu^m)$ with
$\mu^j \in w_0(P^+_{n,d^j}), j=0,\ldots,m$.

Applying formula \eqref{Whitaction1}, we obtain
\begin{multline}   \label{W+descr}
_x'\hspace*{-1mm}\He_n^{\omega_1}
\left(\W^{d+1}_E|_{\ol\Q^{\ol\nu,\ol x}} \right) \simeq \\
\underset{\ol{\mu}'}\oplus \; \Psi^{\ol\mu',\ol x}\otimes
\Hom_{GL_n}(V_{\mu^0{}'},\on{Sym}^{d^0}(V\otimes E_x)\otimes V^*)
\otimes E_{\mu^1,x^1} \otimes ... \otimes E_{\mu^m,x^m}\otimes
\Ql(\frac{d}{2})[d],
\end{multline}
where $\ol{\mu}'=(\mu^0{}',\mu^1,...,\mu^m)$ with $\mu^0{}'$ running
over the set $w_0(P^+_n)$, and $\mu^j$ running over the set
$w_0(P^+_{n,d^j})$ for $j=1,\ldots,m$. Here $V^*$ is the
representation of $GL_n(\Ql)$ dual to $V$.

We have a stratification of the stack $\ol\Q^{\ol\nu'',\ol x}$
analogous to that described in \lemref{strata}. We list only the
strata that can possibly support a sheaf of the form $\Psi^{\ol\mu',\ol
x}$. Those are $\Q^{\ol\mu',\ol x}$, where $\ol\mu'=
(\mu^0{}',\ldots,\mu^m)$ are such that $\mu^0{}' \geq \nu^0{}'', \mu^j \geq
\nu^j, j=1,\ldots,m$. (We recall that the inequality $\la \geq \la'$
means that $\la$ belongs to the set $\la'+R_+$, where $R_+$ is the set
of all linear combinations of simple roots $\al_i, i=1,\ldots,n-1$, of
$GL_n$ with non-negative integer coefficients.)

The stratum $\Q^{\ol\mu',\ol x}$ belongs to the substack
$\ol\Q^{\ol\nu',\ol x}$ if and only if in addition $\mu^0{}' \geq
\nu^0{}' = \nu^0{}'' + \al_1$.

Recall that each sheaf $\Psi^{\ol\mu{}',\ol x}$ is the extension by zero
of its restriction to the stratum $\Q^{\ol\mu{}',\ol x}$. The stratum
with $\ol\mu{}'$ appearing in the summation of the RHS of formula
\eqref{W+descr} belongs to $\ol\Q^{\ol\nu',\ol x}$ if and only if
$d^0>1$, $w_0(\mu^0{}') \in P^+_{n,d^0-1}$, and $w_0(\mu^j) \in
P^+_{n,d^j}$ for all $j=1,\ldots,m$. All of these strata belong to $X
\times \ol\Q^d \subset \Qp^d$. Therefore $\Wp^d$ is supported on
$X\times \ol\Q^d \subset \Qp^d$. This proves the first assertion of
\propref{Hecke on Q}.

Furthermore, we have for $d^0>1$ and any $\mu^0 \in
w_0(P^+_{n,d^0-1})$:
$$\Hom_{GL_n}(V_{\mu^0{}'},\on{Sym}^{d^0}(V\otimes E_x)\otimes V^*)
\simeq E_x\otimes E_{\mu^0{}',x}.$$ Combining this with formulas
\eqref{cart cor}, \eqref{W+descr} and \eqref{restr of WE}, we obtain
the desired isomorphism ${\mathcal W}^d_{E,+}|_{X \times \ol\Q^d}
\simeq E\boxtimes {\mathcal W}^d_E$ over the preimage of each $x\times
D\subset X\times X^{(d+1)}$ in $\Qp^d$.

\ssec{Proof of \propref{Hecke on Q}: global computation}

To complete the proof of \propref{Hecke on Q}, we need to show that
the isomorphism $\Wp^d|_{X \times \ol\Q^d} \simeq E\boxtimes {\mathcal
W}^d_E$ holds globally, and not only on each fiber $\tau^{-1}(x \times
D)$.

In order to show that, we introduce one more stack, $\Qpm^d$. This is
a closed substack of ${}\Qpp^d$ which is the preimage of the incidence
divisor $X\times X^{(d)}\subset X\times X^{(d+1)}$ under the morphism
$\tau:{}\Qpp^d \to X\times X^{(d+1)}$. Equivalently, the stack
$\Qpm^d$ may be defined by the condition that the image of $s_n$ is
contained in $\det\M \subset (\det\M)(x)$.

Let us consider the stack ${}\Mpp$ which classifies the data
$(x,\M_0,\M,\M_0\hookrightarrow \M(x))$, where $x \in X$,
$\deg(\M)-\deg(\M_0)=d$, and the map $\M_0 \to \M(x)$ is such that the
image of $\det \M_0$ is contained in $(\det \M)(x)$ (and not just in
$(\det \M)(n\cdot x)$). Let $\Mpm$ be the closed substack of $\Mpp$,
where the image of $\det \M_0$ is contained in $\det \M$.

Consider the Cartesian product $\Mod_n^{d+1}\underset{\Bun_n}\times
\H^1_n$, which classifies the data
\begin{equation}    \label{data ++}
(x,\M_0,\M,\M',\M_0\hookrightarrow \M',\M\hookrightarrow \M'),
\end{equation}
where $\deg(\M')-\deg(\M_0)=d+1$ and $\M'/\M$ is the simple skyscraper
sheaf supported at $x$. We have a natural pro\-per morphism
\begin{equation*}
\cc: \Mod_n^{d+1}\underset{\Bun_n}\times \H^1_n\to {}\Mpp,
\end{equation*}
which corresponds to "forgetting" $\M'$, and a natural pro\-jec\-tion
$\bb:\Mod_n^{d+1}\underset{\Bun_n}\times\H^1_n\to \Mod_n^{d+1}$, which
corresponds to "forgetting" $\M$. We define the complex $\Ppp$ on
$\Mpp$ as
$$\Ppp:= \cc_! (\bb^*(\P^{d+1}_E))\otimes \Ql(\frac{n}{2})[n]$$
(recall that $\P^d_E:=\pi^*(\L^d_E)\otimes \Ql(\frac{d\cdot n}{2})[d\cdot n]$.)
Let $\Ppm$ be the restriction of $\Ppp$ to
$\Mpm$ tensored with $\Ql(\frac{-1}{2})[-1]$.

Now form a commutative diagram, in which the left square is Cartesian:
\begin{equation*}
\begin{CD} \ol\Q^0 @<{\hlpm}<< Z^d_{+-}
@>{\hrpm}>> \Qpm^d \\ @VqVV @V{'q}VV @VqVV \\ \Bun_n
@<{\hl}<< \Mpm @>{\hr}>> X\times \Bun_n \end{CD}
\end{equation*}

Consider the complex
$$\W^d_{E,+-}:={}\hrpm_!({}\hlpm{}^*(\Psi^0)\otimes
{}'q^*(\Ppm))\otimes \Ql(\frac{-1}{2})[-1].$$
Since we already know that ${\mathcal W}^d_{E,+}$ vanishes on
$\Qp^d-(\Qp^d \cap \Qpm^d)$, it suffices to show that the restriction
of ${\mathcal W}^d_{E,+-}$ to $\Qp^d \cap \Qpm^d$ is supported on
$X\times \ol\Q^d$, where it is isomorphic to $E\boxtimes
{\mathcal W}^d_E$.

Observe that $X\times\Mod^d_n$ is naturally a closed substack in
$\Mpm$. The following result completes the proof of \propref{Hecke on Q}:

\begin{lem}    \label{last lemma}
{\em (1)} The complex $\Ppm$ is a perverse sheaf, and there is a
natural surjection
\begin{equation}    \label{surj map1}
\Ppm \twoheadrightarrow (E\otimes \Ql(\frac{1}{2})[1])\boxtimes \P^d_E.
\end{equation}

{\em (2)} Let $\K_E$ be the kernel of the map \eqref{surj map1}. The
{\em *}--restriction of $\hrpm_!({}\hlpm{}^*(\Psi^0)\otimes
{}'q^*(\K_E))$ from $\Qpm^d$ to $\Qp^d \cap \Qpm^d$ vanishes
identically.
\end{lem}

\ssec{An informal explanation}

Before giving a formal proof of \lemref{last lemma}, we explain the
main idea behind it. In this discussion we will assume that our ground
field $k$ is algebraically closed.

\medskip

Recall the full Hecke correspondence stack $_x\H_n$. Observe that the
fiber $(\hr)^{-1}(\M')$ of $_x\H_n$ over $\M' \in \Bun_n$ under the
map $\hr$ is isomorphic to the affine Grassmannian $\Gr_x$ (see, e.g.,
\cite{FGKV,FGV}). Let $\OO_x$ be the complete local ring at
$x$. The group $GL_n(\OO_x)$ acts naturally on
$\Gr_x$ and we say that a perverse sheaf on $\Gr_x$ is {\it spherical}
if it is $GL_n(\OO_x)$-equivariant. (In particular, every spherical
perverse sheaf is smooth along the stratification of $\Gr_x$
by $\Gr_x^\lambda:=\Gr_x\cap {}_x\H^\la_n$, $\la \in P^+_n$.)
The category of spherical perverse sheaves is known to be semi-simple
and equivalent as a tensor category
to the category of representations of $GL_n(\Ql)$, with the fiber
functor being the functor of (total) global cohomology (see
\cite{Lu2,Gi,MV,BD}).

It is possible to generalize this equivalence to the case when
$x$ and $\M'$ are allowed to ``move'' along $X\times \Bun_n$.
Namely, let $\H_n$ be the stack classifying quadruples
$(\M,\M',\beta,x)$, where $\M$ and $\M'$ are rank $n$ bundles on $X$,
and $\beta$ is an isomorphism $\M|_{X-x}\simeq \M'|_{X-x}$.
Then there is an equivalence between the category of perverse
sheaves on $X$ equipped with $GL_n(\Ql)$--action and a certain
subcategory of the category of perverse sheaves on $\H_n$. For each
$x \in X$ this equivalence ``restricts'' to the equivalence of
the previous paragraph.

\medskip

Moreover, one can generalize this construction to the case of several
points. For any partition ${\mathbf d} = (d^1,\ldots,d^k)$ of $d$,
consider the open subset $\ovc{X}{}^{\mathbf d}$ of $X^{(d^1)} \times
\ldots \times X^{(d^k)}$ consisting of $k$--tuples of divisors
$(D_1,\ldots,D_k)$, such that $\on{supp} D_i \cap \on{supp} D_j =
\emptyset$, if $i\neq j$. Denote the map $\ovc{X}{}^{\mathbf d} \to
X^{(d)}$ by $p_{\mathbf d}$.  We introduce an abelian category ${\mc
A}^d_n$ as follows.  The objects of ${\mc A}^d_n$ are perverse sheaves
$\F$ on $X^{(d)}$ equipped with a $GL_n(\Ql)$--action, together with
the following extra structure: for each partition ${\mathbf d}$, the
sheaf $p_{\mathbf d}^*(\F)$ should carry an action of $k$ copies of
$GL_n(\Ql)$, compatible with the original $GL_n(\Ql)$--action on $\F$
with respect to the diagonal embedding $GL_n(\Ql) \to
(GL_n(\Ql))^{\times k}$. For different partitions, these actions
should be compatible in the obvious sense. In addition, it is required
that whenever $d^i=d^j, i \neq j$, the action of the $i$th and $j$th
copies of $GL_n(\Ql)$ on $p_{\mathbf d}^*(\F)$ commutes with the
corresponding natural $\ZZ_2$--action on $\ovc{X}{}^{\mathbf d}$. The
definition of morphisms in ${\mc A}^d_n$ is clear.

Let now $\H^{\on{BD},d}_n$ be the symmetrized version of the
Beilinson-Drinfeld affine Grassmannian (see \cite{BD}). By definition,
$\H^{\on{BD},d}_n$ is the ind-stack classifying triples
$(\M,\M',D,\beta)$, where $\M,\M$ are as above, $D\in X^{(d)}$ and
$\beta$ is an isomorphism between $\M$ and $\M'$ away from the support
of the divisor $D$.  In particular, $\Mod_n^d$ is naturally a closed
substack of $\H^{\on{BD},d}_n$, corresponding to the condition that
the meromorphic map $\M'\to \M$ defined by $\beta$ is regular.

One can introduce the notion of a spherical perverse sheaf on
$\H^{\on{BD},d}_n$ and construct an equivalence between the above
category ${\mc A}^d_n$ and the category of spherical perverse sheaves
on $\H^{\on{BD},d}_n$. For example, the perverse sheaf
$\on{Sym}^d(V\otimes E)\otimes \Ql(\frac{d}{2})[d]$, which is
naturally on object of ${\mc A}^d_n$, goes to the sheaf $\P^d_E$
(considered as a sheaf on $\H^{\on{BD},d}_n$ supported on $\Mod_n^d$).

One can also define categories analogous to ${\mc A}^d_n$
over partially symmetrized powers of $X$. From this point of view,
the sheaves on $\Mpm$ that we are interested in correspond to
perverse sheaves on $X\times X^{(d)}$ equipped with
$GL_n(\Ql)$--action and an additional structure as above.

In particular, the sheaf $\Ppm$ corresponds to the restriction to $X
\times X^{(d)} \subset X \times X^{(d+1)}$ of the sheaf
$$(V^* \otimes \Ql) \boxtimes \on{Sym}^{d+1}(V \otimes E)\otimes
\Ql(\frac{d+1}{2})[d+1]$$ on $X \times X^{(d+1)}$.

\medskip

The first assertion of \lemref{last lemma}
then translates into the statement that this restriction is perverse,
and that there is a map
\begin{equation}    \label{surj map}
(V^* \otimes \Ql) \boxtimes \on{Sym}^{d+1}(V \otimes E)|_{X \times
X^{(d)}} \rightarrow (\Ql \otimes E) \boxtimes \on{Sym}^d(V\otimes E),
\end{equation}
which becomes after a cohomological shift by $d+1$ a surjection of
perverse sheaves.

The required map is induced by the obvious map
$$\Ql \boxtimes \on{Sym}^{d+1}(V \otimes E)|_{X \times X^{(d)}}
\rightarrow (V \otimes E) \boxtimes \on{Sym}^d(V \otimes E).$$ (In
fact, for any local system $\wt{E}$ on $X$, the map
$\on{Sym}^{d+1}(\wt{E})|_{X \times X^{(d)}}\to \wt{E}\boxtimes
\on{Sym}^d(\wt{E})$, which is an injection of {\it sheaves} becomes
after a cohomological shift a surjection of {\it perverse sheaves}
such that $\wt{E}\boxtimes \on{Sym}^d(\wt{E})[d+1]$ is the cosocle of
$\on{Sym}^{d+1}(\wt{E})|_{X \times X^{(d)}}[d+1]$, see \cite{Dr}).

Moreover, the kernel of the map \eqref{surj map} is supported on the
incidence divisor $X \times X^{(d)} \subset X \times X^{(d+1)}$, and
there it satisfies the following property. Its stalk at a point
$(x,D_1,\ldots,D_k)$ of $X \times X^{(d)}$ (assuming that $\on{supp}
D_i \cap \on{supp} D_j = \emptyset$, $\deg(D_i)=d^i$) is a
$GL_n(\Ql)^{\times k}$--module which decomposes into irreducible
components of the form $V^{\la^1} \otimes \ldots \otimes V^{\la^k}$,
where at least one $\la^i$ does not belong to $P^+_{n,d^i}$. This
proves the second assertion of \lemref{last lemma}.

In the proof of \lemref{last lemma} given below we simply
perform the same manipulations as above directly in the category
of sheaves on $\Mpm$.

\ssec{Proof of \lemref{last lemma}}

First observe that $\Ppp$ is a perverse sheaf on $\Mpp$, by a standard
smallness result in the theory of the affine Grassmannian.  The
assertion about $\Ppm$ follows because $\Ppp$ has no subquotients
supported over the incidence divisor $X\times X^{(d)}\to X\times
X^{(d+1)}$.

To construct the surjection $\Ppm \twoheadrightarrow (E\otimes
\Ql(\frac{-1}{2})[-1])\boxtimes \P^d_E$ we introduce the stack
$'\Mod_n^{d+1} = \Mod^{d+1}_n\underset{X^{(d+1)}}\times (X\times
X^{(d)})$.

We consider two perverse sheaves on it. The first one, denoted by
$\F_1$, is the pull-back of $\P^{d+1}_E$ under $'\Mod_n^{d+1}\to
\Mod^{d+1}_n$. To construct the other sheaf, consider the morphism
${\mathbf a}: \Mod_n^d\underset{\Bun_n}\times \H^1_n\to
{}'\Mod_n^{d+1}$ defined by sending $(\M_0,\M,\M',\M_0\subset \M,\M
\subset \M')$ to $(\M_0\subset \M') \times (x,D)$, where
$D=\tau(\M_0\subset \M)$. The second perverse
sheaf $\F_2$ is by definition $${\mathbf a}_! \; (\on{id} \times
\supp)^*(\P^d_E\boxtimes E)\otimes \Ql(\frac{n}{2})[n],$$ where
$(\on{id} \times \supp): \Mod_n^d\underset{\Bun_n}\times \H^1_n\to
\Mod_n^d\underset{\Bun_n}\times X$. (Thus, $\F_1$ corresponds to the
sheaf $\Ql \boxtimes \on{Sym}^{d+1}(V \otimes E)|_{X \times X^{(d)}}
\otimes \Ql(\frac{d+1}{2})[d+1]$
and $\F_2$ corresponds to the sheaf $(V \otimes E) \boxtimes
\on{Sym}^d(V \otimes E)\otimes \Ql(\frac{d+1}{2})[d+1]$ on $X \times
X^{(d)}$.)

There is a natural surjective map $\F_1 \twoheadrightarrow \F_2$. Now
the desired map $\Ppm \to (E\otimes \Ql(\frac{1}{2})[1])\boxtimes
\P^d_E$ is obtained from the map $\F_1 \to \F_2$ by adjunction. It is
surjective, because $(E\otimes \Ql(\frac{1}{2})[1])\boxtimes \P^d_E$
has no subquotients supported on proper closed substacks. This
completes the proof of part (1) of the lemma.

\medskip

To prove part (2), we choose $x\times D\in X\times X^{(d)}$ and
calculate the restriction of $\K_E$ to its preimage in $\Mpm$. To
simplify notation, assume that $D$ is of the form $d\cdot x$.

The preimage of $x\times d\cdot x\in X\times X^{(d)}$ under
$\Mpm\overset{\supp}\longrightarrow X\times X^{(d)}$ is naturally a
closed substack of $_x\H_n$. Using formula \eqref{fiber of Laumon's
sheaf}, we obtain that the restriction to $\tau^{-1}(x\times d\cdot
x)$ of the surjection $\Ppm \twoheadrightarrow (E\otimes
\Ql(\frac{1}{2})[1])\boxtimes \P^d_E$ can be identified with the map
\begin{align*}
&\underset{\lambda \in P^+_n}\oplus\, \IC_{-w_0(\lambda)}\otimes
\Hom_{GL_n}(V^\lambda, \on{Sym}^{d+1}(V\otimes E_x)\otimes V^*)\otimes
\Ql(\frac{d+1}{2})[d+1]
\twoheadrightarrow \\
&\underset{\lambda \in P^+_n}\oplus\, \IC_{-w_0(\lambda)}\otimes
\Hom_{GL_n}(V^\lambda, E_x\otimes \on{Sym}^{d}(V\otimes E_x))\otimes
\Ql(\frac{d+1}{2})[d+1].
\end{align*}

\medskip

The kernel of this map is nothing but the restriction of $\K_E$ to the
preimage of $x\times d\cdot x\in X\times X^{(d)}$ in $\Mpm$. It is
clear that if the summand corresponding to the sheaf
$\IC_{-w_0(\lambda)}$ appears in $\K_E$, then $\la \not \in
P^+_{n,d}$. Therefore, the required vanishing follows from formula
\eqref{Whitaction} and \lemref{strata}.


\end{document}